\documentclass{svjour3}                     

\def\strokedint{\int\!\!\!\!\!\!-\:}

\usepackage{hyperref}
\usepackage{graphicx}
\usepackage{graphics}
\usepackage{epstopdf}
\usepackage{framed,color}
\usepackage{color,soul}
\usepackage[lined,boxed,commentsnumbered,ruled,vlined]{algorithm2e}

\usepackage{amssymb,amsmath}
\usepackage{calrsfs}
\usepackage{subcaption}
\usepackage{array}
\usepackage{booktabs}
\usepackage{appendix}

\DeclareMathOperator{\tr}{tr} 

\newcommand{\bb}{\bold{b}}

\newcommand{\bp}{\bold{p}}

\newcommand{\bx}{\boldsymbol{x}}

\newcommand{\bpartial}{\boldsymbol{\partial}}
\newcommand{\beps}{\boldsymbol{\eps}}

\newcommand{\bsigma}{\boldsymbol{\sigma}}

\newcommand{\bzero}{\boldsymbol{0}}

\newcommand{\bC}{\bold{C}}

\newcommand{\bI}{\boldsymbol{I}}

\newcommand{\stiffness}{\bold{K}}

\newcommand{\bN}{\bold{N}}

\newcommand{\bQ}{\boldsymbol{Q}}

\newcommand{\eps}{\varepsilon}
\renewcommand{\t}{^\text{t}}
\newcommand{\T}{^\textrm{T}}

\newcommand{\n}{_{n}}

\newcommand{\odd}{\hspace{0.2mm}$\ddot{\mbox{o}}$}

\newcommand{\real}       {\mathbb{R}}

\def\wbox#1;#2;{\vbox{\hrule\hbox{\vrule height#1mm\kern#2mm\vrule height#1mm}\hrule}}

\renewcommand{\phi}{\varphi}

\def\F#1#2{{}_{#1}\kern-0.5ptF_{#2}}
\def\U#1#2{{}_{#1}\kern-0.5ptU_{#2}}

\def\odd#1{\kern-#1pt{}^{\textrm{o}}\,\,}

\newcommand\lshad{{[\kern-0.15em[}}
\newcommand\rshad{{]\kern-0.15em]}}

\newcommand{\vect}[1]{\bold{#1}}

\def\divb{{\bf div}}
\def\dr{{\rm d}}

\def\x{{\bf x}}

\def\t{{\vect{ t}}}
\def\u{{\vect{ u}}}
\def\v{{\vect{ v}}}
\def\S{{\vect{ S}}}

\def\n{{\bf n}}
\def\f{{\vect{ f}}}

\def\F{{\vect{ f}}}

\def\V{{\bf V}}

\def\bD{{\bf D}}

\def\bI{{\bf I}}
\def\tr{\textrm{tr}}
\def\sgrad{\boldsymbol{\varepsilon}}

\def\eps{\boldsymbol{\varepsilon}}

\def\P{{\mathcal P}}  

\def\O{\Omega}

\def\Th{{\mathcal T}_h}
\def\E{E}

\def\VE{{\V_{h|\E}}}

\def\P{\mathsf P}

\definecolor{mygreen}{rgb}{0.2,0.8,0.2}

\begin{document}

\title{Arbitrary order 2D virtual elements for polygonal meshes: Part I, elastic
problem
\thanks{The first author gratefully acknowledges the partial financial support of the Italian Minister of University and Research, MIUR  (Program: Consolidate the Foundations 2015; Project: BIOART; Grant number (CUP): E82F16000850005).\\ 
The third author was partially supported by IMATI-CNR of Pavia, Italy. This support is gratefully acknowledged.}
}

\titlerunning{PART ONE}        

\author{E. Artioli \and
L. Beirao da Veiga \and 
C. Lovadina \and 
E. Sacco
}


\institute{E. Artioli \at
           Department of Civil Engineering and Computer Science, 
           University of Rome Tor Vergata,
           Via del Politecnico 1, 00133 Rome, Italy. 	
           \and
           L. Beirao da Veiga \at
           Dipartimento di Matematica e Applicazioni, Universtit\`a di Milano-Bicocca, Via Cozzi 53, 20125 Milano, Italia.  \email{lourenco.beirao@unimib.it}
           \and
           C. Lovadina \at
           Dipartimento di Matematica, Universit\`a di Milano, Via Saldini 50, 20133 Milano, Italia. \email{carlo.lovadina@unimi.it}
           \and
           E. Sacco \at
           Dipartimento di Ingegneria Civile e Meccanica, Universtit\`a di Cassino e del Lazio Meridionale, Via G. Di Biasio 43, 03043 Cassino, Italia. \email{sacco@unicas.it}
}
\date{Received: date / Accepted: date}

\maketitle

\begin{abstract}
The present work deals with the formulation of a Virtual Element Method (VEM) for two dimensional structural problems.
The contribution is split in two parts: in part I, the elastic problem is discussed, while in part II [\cite{ABLS_part_II}] the method is extended to material nonlinearity, considering different inelastic responses of the material.
In particular, in part I a standardized procedure for the construction of all the terms required for the implementation of the method in a code is explained. The procedure is initially illustrated for the simplest case of quadrilateral virtual elements with linear approximation of displacement variables on the boundary of the element. Then, the case of polygonal elements with quadratic and, even, higher order interpolation is considered.  The construction of the method is detailed, deriving the approximation of the consistent term, the required stabilization term and the loading term for all the considered virtual elements. A wide numerical investigation is performed to assess the performances of the developed virtual elements, considering different number of edges describing the elements and different order of approximations of the unknown field. Numerical results are also compared with the one recovered using the classical finite element method.
\end{abstract}


\section{Introduction} \label{s:intro}

The Virtual Element Method (VEM) was recently introduced in \cite{volley,hitchhikers} as a generalization
 of the Finite Element Method (e.g. \cite{Bathe1996,Hughes:book,Wriggers:book}) with the additional capability of dealing with very general polygonal/polyhedral meshes and the possibility to easily implement
highly regular discrete spaces (see for instance \cite{Brezzi:Marini:plates}). By making use also of non-polynomial shape functions, but
avoiding the explicit construction of the local basis
functions, the VEM can easily handle general polygons/polyhedrons
without complex integrations on the element. 
Polytopal meshes can be very useful for a wide range of
reasons, including meshing of the domain (such as cracks) and data
(such as inclusions) features, automatic use of hanging nodes, use
of moving meshes, adaptivity. Recently, the interest in numerical methods that can make use of general
polytopal meshes has undergone a significant growth in the
mathematical and engineering literature. 
Here, we limit to citing some work in the area of polygonal Finite Elements \cite{ST04} applied to advanced problems in structural mechanics, that is \cite{Paulino-nonlinear-polygonal,TPPM10,Wriggers-2014,Paulino-cracks}, see also \cite{DiPietro-Ern-1}.

The Virtual Element Method has strongly developed and is now able to cope with a wide range of topics and problems, see for instance \cite{CH-VEM} and references therein. In the more specific framework of structural mechanics, VEM has been introduced in \cite{VEM-elasticity} for (possibly incompressible) two dimensional linear elasticity and general ``polynomial'' order, in \cite{Paulino-VEM} for three dimensional linear elasticity and lowest order, in \cite{BeiraoLovaMora} for general two dimensional elastic and inelastic problems under small deformations (lowest order). More specific applications have been also considered, such as  contact problems \cite{wriggers}, topology optimization \cite{Topology-VEM} and geomechanics \cite{Andersen-geo}. 

The present contribution is the first part of a work made of two papers, that take the steps from \cite{BeiraoLovaMora,VEM-elasticity} and extend the VEM to the case of arbitrary order of accuracy (or ``polynomial'' order), in the framework of general elastic and inelastic small deformation 2D problems. In brief, we are able to develop a numerical scheme of arbitrary order of accuracy that works for nonlinear deformation problems, general polygonal meshes and still yields a conforming solution. The key idea of the method is to use a projection of the strains on a polynomial space, apply the constitutive law to this approximated strain at specific points (which may be seen as the analogous of the Gauss points in standard FEM), and finally stabilize the ensuing discrete stiffness operator by standard techniques in the VEM literature.

This first part of the work is focused on linear elasticity and sets all the basic tools necessary for the more advanced problems considered in part two. Namely, we introduce the Virtual Element space, the associated degrees of freedom and projectors, and the approximation of the elasticity stiffness matrix. When compared to the VEM in \cite{VEM-elasticity}, the present scheme has strong similarities (both handle arbitrary order and linear elasticity). Nevertheless, in this paper there are three important additional issues: 
\begin{enumerate}
\item Differently from \cite{VEM-elasticity}, the presented method is already suitable for generalization to nonlinear problems (treated in Part II), without any change in the foundations of the scheme.
\item An important contribution of this paper is that it gives extensive guidelines on the construction and coding of the proposed scheme, starting from the simple case of lowest order on quadrilaterals and going up (step by step) to general order and general polygons. Although it already exists a well written coding guide for VEM \cite{hitchhikers}, that work is more focused on the diffusion problem rather than elasticity; moreover, here we propose a coding approach that is more suitable for engineers. The reader interested in the coding of VEM will surely benefit from reading both works;
\item In this paper we also present, for the first time, a deep numerical testing of a VEM method of general order in the realm of linear elasticity, with highly encouraging results.
\end{enumerate}
The paper is organized as follows. In Section \ref{s:one} we present the continuous problem and the general structure of the Virtual Element formulation. In Section \ref{s:new_form} we thoroughly  describe the construction of the method. To help the reader understand all the details, we start with the lowest order case (linear accuracy) on quadrilaterals in Section \ref{ss:kone}, then consider the quadratic case on general polygons in Section \ref{ss:ktwo}, and only as a last step we consider the fully general case, see Section \ref{ss:kgen}. Finally, in Section \ref{s:num_test} we present a set of numerical tests for different orders and different families of meshes.

\section{General structure of the Virtual Element formulation}
\label{s:one}

In the present section we introduce the continuous problem and give an overview on the proposed Virtual Element discretization. More details and specific explanations will be found in the following sections.

\subsection{The continuous problem}
\label{ss:one}

In this section we present the linear elastic problem in the two-dimensional (2D) framework, under the assumption of infinitesimal strain and small displacements. In the following the Voigt notation is adopted, so that stress and strain tensors are represented as $3-$component vectors, and the fourth-order constitutive tensor is represented as a $3\times3$ matrix.

Let $\O$ be a continuous body occupying a region of the two-dimensional space $\mathcal{R}^2$ in which the Cartesian coordinate system $(\textrm{O},x,y)$ is introduced.
The displacement field is denoted by the vector $\u(x,y)=\left( u \; v \right)^\textrm{T}$ and the associated strain defined as:
\begin{equation}
\label{eq:epsilon}
\beps(\u) = \S \u \qquad \textrm{with} \;
\S = \left [ \begin{matrix} \partial_{x} & 0 \\ 0 & \partial_{y} \\ \partial_{y} & \partial_{x} \end{matrix} \right] \, ,
\end{equation}
where $\partial_{(\bullet)}$ denotes the partial derivative operator with respect to the $(\bullet)$-coordinate.

The linear elastic constitutive law is considered for the body $\O$, so that the stress is given by the relationship:
\begin{equation}
\bsigma = \bC \, \beps 
\label{eq:sigma}
\end{equation}
where $\bC = \bC(x)$ is the $3\times3$ (uniformly positive) elastic matrix, possibly depending on the position vector $\x=( x \; y ) \T \in \O$. 

The body $\O$ is subjected to distributed volume forces $\bb$. Without loosing in generality, for simplicity it is assumed that the displacements are vanishing on the whole boundary of $\O$ (see Remark \ref{BCs} at the end of this section). 
The variational formulation of the elastostatic problem is provided by the virtual work principle:
\begin{equation}
\left\lbrace{
\begin{aligned}
& \textrm{Find $\u\in\V$ such that} \\
& a(\u,\v) = <\bb,\v> \quad \forall \, \v \in \V
\end{aligned}
}\right.
\label{eq:vwork}
\end{equation}
where $\V:=(H^1_0(\O))^2$ is the space of the admissible displacement field and
\begin{equation}
\label{eq:weak_form}
\begin{aligned}
a(\u,\v) & = \int_{\O} \eps(\v)^T \, \bC \,\eps(\u)\,\dr\x
                = \int_{\O} (\S  \v)^T \, \bC \,\S \u\,\dr\x \\
<\bb,\v> & = \int_{\O} \v^T \, \bb \,\dr\x
\end{aligned}
\end{equation}
The form $a(\cdot,\cdot)$ is symmetric, continuous and coercive on $\V$, so that problem \eqref{eq:vwork} is well posed.

\subsection{The Virtual Element formulation}
\label{s:abs}

We start by presenting the discrete (virtual) space of admissible displacements, that is the same one introduced in \cite{VEM-elasticity}. See also \cite{volley,hitchhikers} for a scalar version of the space used for the diffusion problem.

\begin{figure}
\begin{minipage}[b]{.5\linewidth}
\centering
\includegraphics[bb = 5 18 20 300, scale=0.5]
                {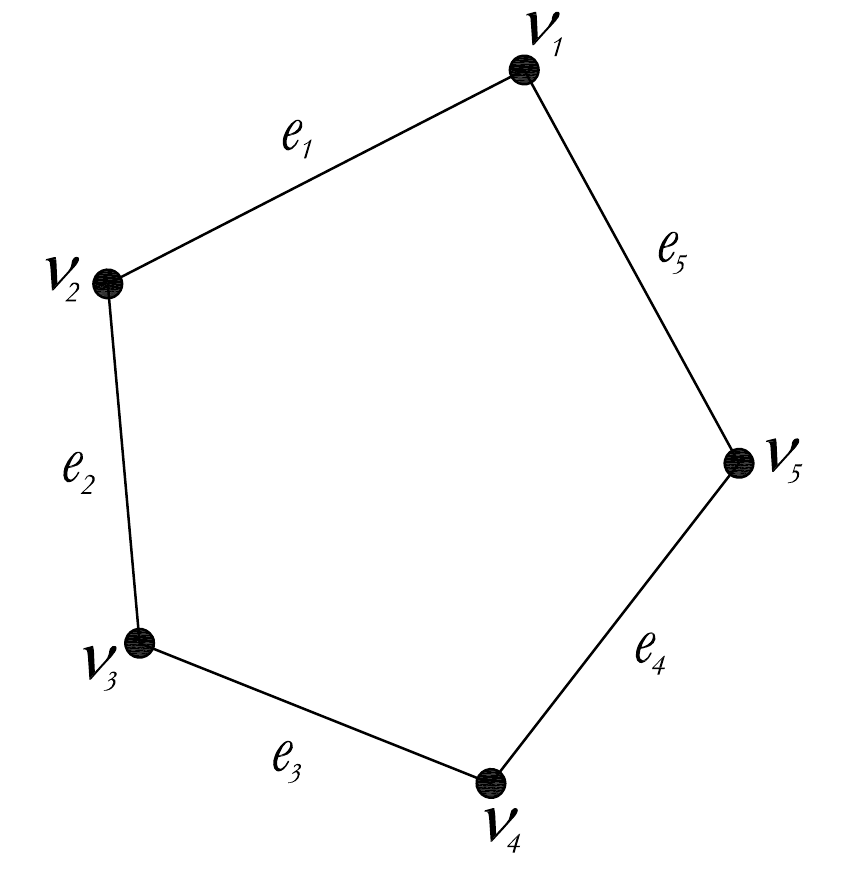}
\end{minipage}
\caption{
A sample element with five edges. The symbols $\nu_i$ and $e_i$, $i=1,2,..,5$, denote vertexes and edges, respectively.}
\label{fig:elem:num}
\end{figure}

Let $\Th$ be a {\it simple polygonal mesh} on $\O$. This can be any decomposition of $\O$ in non overlapping polygons $\E$ with straight edges. Let ${\cal E}_h$ denote the set of the edges of $\Th$. The symbol $m$ represents the number of edges of a polygon $E$,and the typical edge of the polygon $E$ is indicated by $e$, (i.e. $e \in \partial E$). A sample element $E$, with an ordering of vertexes and edges, is depicted in figure \ref{fig:elem:num}. 
The space $\V_h$ will be defined element-wise, by introducing local spaces $\VE$ and the associated local degrees of freedom, as in standard Finite Element (FE) analysis. Instead, differently from standard FE, the definition of the local spaces $\VE$ is not fully explicit. 

Let $k$ be a positive integer, representing the ``degree of accuracy'' of the method. Then, given an element $\E\in\Th$, we start by defining the scalar space:
\begin{equation}
\begin{aligned}
V_{h|E}=\big\{ & v_h \in H^1(\E) \cap C^0(\E) \ : \  \Delta v_h \in \P_{k-2}(\E) , \\
& v_h|_e \in \P_k(e)   \quad \forall e \in \partial\E \big\} ,
\end{aligned}
\end{equation}
where $\P_k(E)$ denotes the space of polynomials of degree (up to) $k$ defined on $E$, with the agreement that $\P_{\!-1}(E) = \{ 0 \}$.
The proposed virtual displacement space is
$$
\VE = \left [ V_{h|E}\right]^2 .
$$

The space $\VE$ is made of vector valued functions $\v_h$ such that \cite{VEM-elasticity,volley}:
\begin{itemize}
\item $\v_h$ is a polynomial of degree $\le k$ on each edge $e$ of $\E$, i.e. $\v_{h|e} \in(\P_k(e))^2$;
\item $\v_h$  is globally continuous on $\partial\E$;
\item $\divb \big( \bC \eps(\v_h) \big)$ is a polynomial of degree $\le k-2$ in $\E$.
\end{itemize}
The following important observations hold:
\begin{itemize}
\item the functions $\v_h\in\VE$ are explicitly known on $\partial\E$;
\item the functions $\v_h \in\VE$ are not explicitly known inside the element $\E$;
\item it holds $(\P_k(\E))^2 \subseteq \VE$ (that is important for approximation).
\end{itemize}
A sample polygon with $m=5$ edges is represented in Figure \ref{fig:elem:num} with indication of vertexes and edges.
Taking into account that the unknown vector $\v_h \in \VE$ has two components, for any $\E\in\Th$, the degrees of freedom for $\VE$ will be:
\begin{itemize}
\item $2m$ pointwise values $\v_h(\nu_i)$ at corners $\nu_i$ of $E$, $i=1,2,..,m$,
\item $2m(k-1)$ pointwise values $\v_h(\bx^{e}_j)$ at edge internal nodes  $\{\bx^{e}_j \}$ (i.e. not comprising extrema of edge $e$), $j=1,...,k-1$, being $m$ the number of edges of a polygon;
\item $2\frac{k(k-1)}{2}$ scalar moments of the unknown field over the element, not associated with a specific location over $E$ (these degrees of freedom will be specified in the following).
\end{itemize}
Sample figures for the degrees of freedom will be shown in the following sections, where more details are given.
The dimension of the space $\VE$ results 
\begin{equation}
\label{dofsnumber}
n = \textrm{dim}(\VE) = 2mk + k (k-1) \: ,
\end{equation}
where $2m$ degrees of freedom are associated with the $m$ corners of $\partial E$, $2m(k-1)$ degrees of freedom are associated with the $m$ edges (with $(k-1)$ internal nodal points per edge), and $k(k-1)$ degrees of freedom are the moments.
A local ``Lagrange-type'' basis for the space $\VE$ is chosen as follows: the $i-th$ function of the local basis is given by the unique function in $\VE$ that has unit value on the $i-th$ degree of freedom and vanishes for the remaining ones.

As in standard FE methods, the global space $\V_h \subseteq \V$ is built by assembling the local spaces $\VE$ as usual:
\begin{equation}
\V_h = \{ \v \in \V \ : \ \v |_{\E} \in \VE \; \forall \E\in\Th \} .
\end{equation}

A virtual element method (VEM) for the above problem is constructed following a procedure that closely resembles the Finite Element approach, i.e. by restricting the original variational formulation \eqref{eq:vwork} to the discrete space $\V_h$ and approximating the ensuing terms:
\begin{equation}\label{gen-vem-str}
\left\lbrace{
\begin{aligned}
&\mbox{find } \u_h\in \V_h \mbox{ such that}\\
& a_h(\u_h,\v_h)=<\bb_h,\v_h>\qquad\forall \,\v_h\in \V_h,
\end{aligned}
} \right.
\end{equation}
where
\begin{itemize}
\item $\V_h \subset \V$ is the virtual space introduced above;
\item $a_h(\cdot,\cdot) \: : \: \V_h \times \V_h \rightarrow {\mathbb R}$ is a discrete bilinear form approximating the continuous form $a(\cdot,\cdot)$;
\item $<\bb_h,\v_h>$ is the term approximating the virtual work of the external load. 
\end{itemize}

The discrete bilinear form is built element by element, assuming:
\begin{equation}
a_h(\u_h,\v_h)=\sum_{{\E\in\Th}} a_h^{\E}(\u_h,\v_h)\quad\forall\, \u_h,\,\v_h\,\in \V_h.
\end{equation}
On the other hand, since the functions of the space $\V_h$ are not know explicitly inside the elements, the bilinear forms $a_h^E(\cdot,\cdot)$ cannot be evaluated by standard Gauss integration. 
The construction of the local discrete form $a_h^E(\cdot,\cdot)$, a key point in the VEM, is performed following the procedure described in the following.

Since neither the function $\v_h$ nor its gradient are explicitly computable in the element interior points, the method proceeds by introducing a projection operator $\Pi$, representing the approximated strain associated with the virtual displacement, defined as:
\begin{align}
\label{uVuP}
\Pi :  \;\;\;\; & \VE \longrightarrow \P_{k-1}(\E)^{2\times 2}_{\textrm{sym}} \\
\nonumber
& \quad \v_h \; \mapsto \quad \Pi(\v_h) ,
\end{align}
where, albeit the adopted notation, the $2 \times 2$ symmetry character of the strain tensor is remarked. It is noted that $\ell =  \textrm{dim}(\P_{k-1}(\E)^{2\times 2}_{\textrm{sym}}) = 2 \frac{k (k+1)}{2}$.

Given $\v_h \in \VE$, such an operator $\Pi$ is defined as the unique function $ \Pi(\v_h) \in \P_{k-1}(\E)^{2\times 2}_{\textrm{sym}}$ that satisfies the condition:
\begin{equation}
\int_{\E} \Pi(\v_h)^T  \sgrad^{P}  = 
\int_{\E}      \eps (\v_h)^T  \sgrad^{P}  ,\quad \forall \sgrad^{P} \in \P_{k-1}(\E)^{2\times 2}_{\textrm{sym}} \,.
\label{eq:def-proj}
\end{equation}
This operator represents the best approximation of the strains (in the square integral norm) in the space of piecewise polynomials of degree $(k-1)$. Although the functions in $\VE$ are virtual, the right hand side in \eqref{eq:def-proj} (and thus the operator $\Pi$) turns out to be computable with simple calculations, as we will show in the next section.

Once computed, the $\Pi$ operator approximates the internal work part associated to the equilibrium weak formulation as follows:
\begin{equation}
\label{eq:approx_bilinear_form}
a^E_h(\u_h , \v_h) = \int_{\E} \left[\Pi (\u_h)\right]^T \bC \Pi (\v_h) + 
 S^E \left( \u_h , \v_h \right)
\end{equation}
with $S^E(\cdot,\cdot)$ a suitable stabilizing term, needed to preserve the coercivity of the system, and  described in the sequel. 
The numerical approximation of the loading term is presented in the next section.

\begin{remark}\label{BCs}
Since on the boundary of the elements the functions of $\V_h$ are piecewise polynomials, the imposition of more general boundary conditions follows exactly the same procedure as in standard finite elements.
For instance, non-homogeneous Dirichlet boundary conditions are enforced by standard interpolation of the prescribed displacements, while Neumann boundary conditions are implemented by adding the usual boundary integral to the formulation.
\end{remark}

\section{Construction of the method}
\label{s:new_form}

In the present section we detail the construction of the proposed virtual scheme. In order to make things clearer, our presentation will be subdivided into three parts. Initially, we start describing the simple case $k=1$ for quadrilateral elements $E$. Afterwards, we describe the case $k=2$ for general polygons $E$ with $m$ edges. Last, we address the case of general polynomial order $k$.

\begin{remark}\label{rem:coord}
In the following, we adopt local (to the current element $E$) scaled space coordinates, such that: (1) the origin of the coordinate axes is the centroid of $E$; (2) all the coordinate values are scaled by the diameter $h_E$ of the element $E$.
For instance, if the position of a point $P$ is $(x \; y)$ in standard cartesian coordinates, then its local scaled coordinates $(\xi \; \eta)$ are given by
\begin{equation}
\label{eq:adimcoor}
\xi = \frac{x - x_c}{h_E} \ , \quad \eta = \frac{ y - y_c}{h_E} \:,
\end{equation}
where $( x_c \; y_c)$ denotes the cartesian coordinates of the centroid of $E$.
\end{remark}

\subsection{The case $k=1$ for quadrilaterals}
\label{ss:kone}
\begin{figure}
\begin{minipage}[b]{.5\linewidth}
\hspace{10mm}
\includegraphics[bb =  10 0 30 280, scale=0.5]{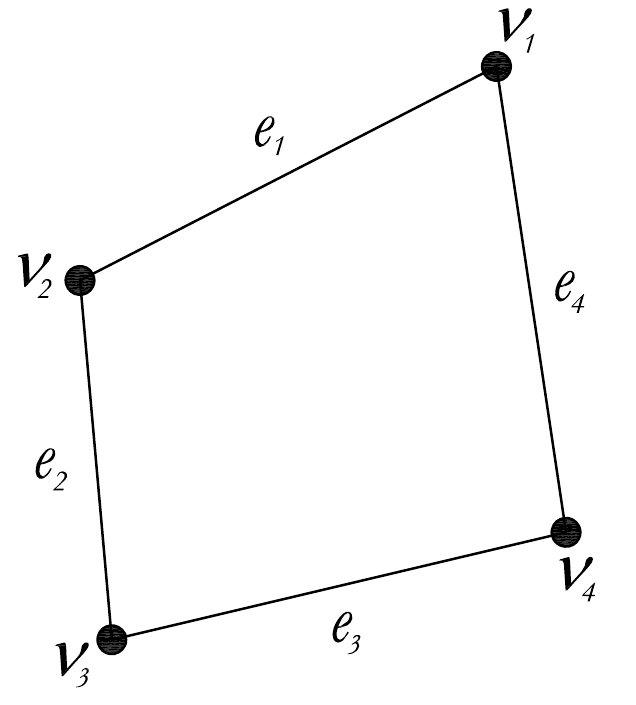}
\subcaption{}
\end{minipage}
\begin{minipage}[b]{.5\linewidth}
\hspace{10mm}
\includegraphics[bb = 30 0 50 280, scale=0.5]{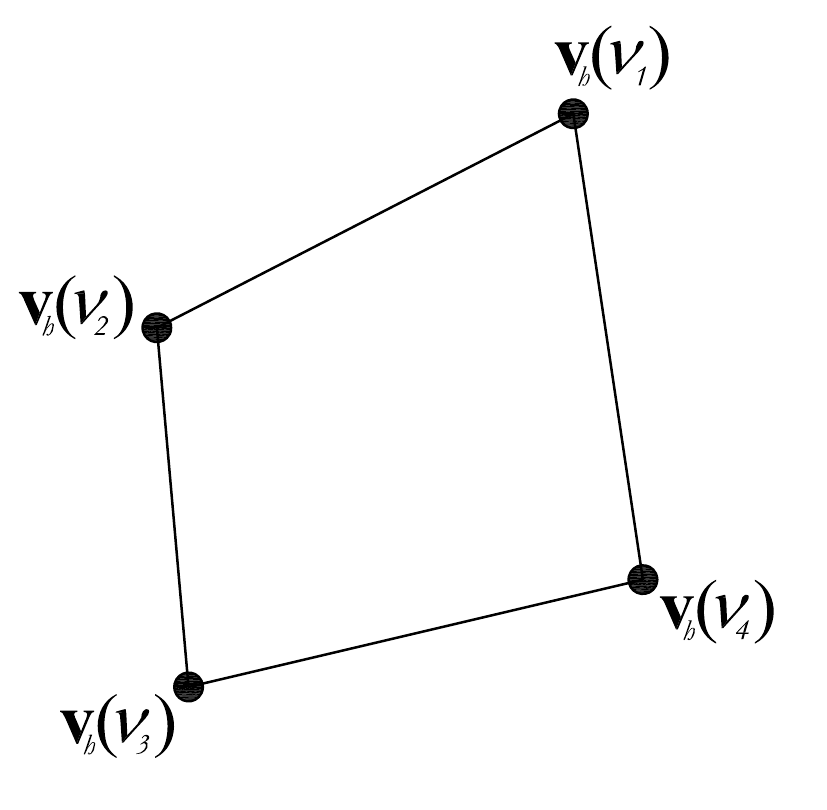}
\subcaption{}
\end{minipage}
\caption{Sample quadrilateral (a) and associated degrees of freedom for the case $k=1$ (b): displacement values at the four corners.}
\label{fig:dof:k1}
\end{figure}

We start considering the case of $E$ being a quadrilateral polygon, i.e. $m=4$, and assuming $k=1$; in this case, the degrees of freedom for $\VE$ are the only displacement field values at the corners of the element, with no edge degrees of freedom and no moments. The degrees of freedom are depicted in figure \ref{fig:dof:k1}.

Therefore, the virtual element space and the polynomial space of approximated strains, on the basis of Eq. \eqref{uVuP}, satisfy the following conditions:
\begin{equation}
\begin{aligned}
n =  \textrm{dim}(\VE) = 4 \times 2 = 8 \\
\ell =  \textrm{dim}(\P_0(\E)^{2\times 2}_{\textrm{sym}}) = 3 \;,
\end{aligned}
\end{equation}
so that it is set
\begin{equation}
\P_0(\E)^{2\times 2}_{\textrm{sym}} = \textrm{span}  \left\{ \begin{pmatrix} 1 \\ 0 \\ 0 \end{pmatrix} ,  \begin{pmatrix} 0 \\ 1 \\ 0 \end{pmatrix}, \begin{pmatrix} 0 \\ 0 \\ 1 \end{pmatrix} \right\} \; .
\end{equation}
Formally, by representing the virtual displacement field and the approximated strains in terms of the basis functions, we can write 
\begin{equation}
\left\{
\begin{aligned}
& \v_h = \bN^V \tilde{\v}_h \;, \\
& \sgrad^{P} = \bN^P \hat{\sgrad} \; ,
\end{aligned}
\right.
\label{eq:u-approx}
\end{equation}
where $\v_h \in \VE$, $\tilde{\v}_h \in \real^8$, $\sgrad^{P} \in \P_0(\E)^{2\times 2}_{\textrm{sym}}$, $\hat{\sgrad} \in \real^3$. The matrix $\bN^P$ is given by
\begin{equation}
\bN^P=\left[
\begin{matrix}
1 & 0 & 0\\
0 & 1 & 0 \\
0 & 0 & 1
\end{matrix}
\right].
\end{equation}
Matrix $\bN^V$ can be viewed as the shape function operator of a typical displacement-based 2D finite element, even if its entries are not analytical over the element interior. In particular, it has the same structure, i.e. odd-index columns are associated with $u_h$ and have zero second entries, whereas even-index columns are associated with $v_h$ and have zero first entries. Since the only degrees of freedom are the vertex values and recalling that the displacements have two components, $\bN^V$ is a $2 \times 8$ matrix operator. Ordering the four vertexes of $E$ with $i = 1,...4$, for $s = 2i-1$, the first entry of the $s^{th}$-column of $\bN^V$ is formally given by the unique function of $V_{h|E}$ that takes the value $1$ at vertex $i$ and vanishes at all other vertexes, the second entry being $0$. Analogously, for $s=2i$, the $s^{th}$-column first entry is zero; the second one is formally given by the unique function of $V_{h|E}$ that takes the value $1$ at vertex $i$ and vanishes at all other vertexes. The operator $\bN^V$ is key for understanding the present construction, but is actually not computable in the interior of the element. On the other hand, since the restrictions of the functions of $\VE$ to $\partial E$ are piecewise linear and continuous, the computation of $\bN^V$ can be carried out easily on the element boundary.

\subsubsection{Consistent term}
\label{sss:cons_term}
In the present section we detail the construction of the consistent part of the bilinear form $a^E_h(\cdot, \cdot)$, i.e the first term appearing in the right hand side of \eqref{eq:approx_bilinear_form}.
We start by the practical construction of the projection operator $\Pi$. 
We denote by $\Pi^{\textrm{m}}$ the operator $\Pi$ expressed as a matrix, in terms of the basis of $\VE$ and 
$\P_0(\E)^{2\times 2}_{\textrm{sym}}$. In other words, the projected strain results in
\begin{equation}
\Pi(\v_h) = \bN^P \Pi^{\textrm{m}} \tilde{\v}_h
\label{eq:pi_matr}
\end{equation}
with $\Pi^{\textrm{m}} \in \real^{3\times 8}$.
The following steps lead to the computation of the matrix operator $\Pi^{\textrm{m}}$.

With the above definitions,   \eqref{eq:def-proj} becomes:
\begin{equation}
\label{eq:X1}
\int_{\E} \left(\bN^P \Pi^{\textrm{m}} \tilde{\v}_h \right)^{\textrm{T}} \bN^P \hat{\sgrad} = 
\int_{\E} \left[ \sgrad(\bN^V \tilde{\v}_h) \right]^{\textrm{T}} \bN^P \hat{\sgrad} 
\quad \forall \hat{\sgrad} \in \real^3 .
\end{equation}
Taking into account Eq. \eqref{eq:epsilon} and noting that the internal contribution vanishes since 
$\sgrad^{P} \!=\! \bN^P \!\! \hat{\sgrad}$ is constant on the element, an integration by parts on the right hand side yields
\begin{equation}
\label{eq:X2}
\int_{\E} \left(\bN^P \Pi^{\textrm{m}} \tilde{\v}_h \right)^{\textrm{T}} \bN^P \hat{\sgrad} = 
\int_{\partial\E} \left(\bN^V \tilde{\v}_h \right)^{\textrm{T}} \left(\bN_E \bN^P \hat{\sgrad} \right) 
\end{equation}
for any $\hat{\sgrad}$ in $\real^3$, where the matrix   
$$
\bN_\E = \left[ \begin{matrix} n_x &  0 & n_y \\ 0 & n_y & n_x\end{matrix} \right] 
$$
is built by the components of the outward normal $\n_\E=(n_x,n_y)^T$ to the element edges.

We can re-write Eq. \eqref{eq:X2} more compactly as:
\begin{equation}
\label{eq:XX3}
\hat{\sgrad}^{\textrm{T}}  \mathcal{G} \Pi^{\textrm{m}} \tilde{\v}_h  = 
\hat{\sgrad}^{\textrm{T}} \mathcal{B} \tilde{\v}_h \quad \forall \hat{\sgrad} \in \real^3 \;,
\end{equation}
where:
\begin{eqnarray}
\label{eq:XX4}
\mathcal{G} = \int_{\E} \left(\bN^P\right)^{\textrm{T}} \bN^P \in \real^{3 \times 3}\\
\mathcal{B} = \int_{\partial\E} \left(\bN_E \bN^P\right)^{\textrm{T}} \bN^V \in \real^{3 \times 8} \;.
\end{eqnarray}

\begin{remark}
Note that, in light of the previous observations, also the matrix $\mathcal{B}$ is computable, since it only requires knowledge of $\bN^V$ on the boundary of the element.
\end{remark}

Eq. \eqref{eq:XX3} provides the required operator $\Pi^{\textrm{m}}$ as the solution of a linear system:
\begin{equation}
\label{eq:Pi_operator}
\Pi^{\textrm{m}} = \mathcal{G}^{-1} \mathcal{B} ,
\end{equation}
where $\mathcal{G}$ is clearly symmetric and positive definite.

Combining eqs. \eqref{eq:pi_matr} and \eqref{eq:Pi_operator}, from \eqref{eq:approx_bilinear_form}, one obtains: 
\begin{align}
\label{eq:approx_bilinear_form_Pi}
\int_{\E} \left( \Pi(\v_h)\right)^{\textrm{T}} \bC  \Pi (\u_h)  = 
\int_{\E} \left( \bN^P \mathcal{G}^{-1} \mathcal{B} \tilde{\v}_h \right)^{\textrm{T}} \bC \left( \bN^P \mathcal{G}^{-1} \mathcal{B} \tilde{\u}_h \right) = 
\tilde{\v}^{\textrm{T}}_h \stiffness_{\textrm{c}} \tilde{\u}_h
\end{align}
which defines explicitly the first part of the stiffness matrix:
\begin{equation}
\label{eq:consistent_stiffness}
\stiffness_{\textrm{c}} = \mathcal{B}^{\textrm{T}} \mathcal{G}^{\textrm{-T}} 
\left[ \int_{\E}  \left(\bN^P \right)^{\textrm{T}} \bC \bN^P \right ] 
\mathcal{G}^{-1} \mathcal{B} .
\end{equation}
As already noted in \eqref{eq:approx_bilinear_form}, the above term is not sufficient for the stability of the method, and needs to be completed with a stabilization term that is discussed in the next section.

\subsubsection{Stabilization term}
\label{sss:stab_term}

The discrete bilinear form cannot be composed of the consistent term alone, otherwise the coercivity of the system may be lost and hourglass modes could arise. The presence of a stability term is therefore standard for the Virtual Element Method, as illustrated for instance \cite{volley,Brezzi:Marini:plates,VEM-elasticity,hitchhikers,BeiraoLovaMora,CH-VEM} for a few examples. 
Herein, we first detail the construction of the stability term, commenting its motivations at the end of the section.

Let a basis for the space $(\P_1(E))^2$ be chosen as:
\begin{equation}
\label{eq:lin_pol_basis}
(\P_1(\E))^2 = \textrm{span}  \left\{ \begin{pmatrix} 1 \\ 0 \end{pmatrix} ,  \begin{pmatrix} 0 \\ 1 \end{pmatrix}, \begin{pmatrix} \xi \\ 0 \end{pmatrix} 
, \begin{pmatrix} 0 \\ \xi \end{pmatrix} , \begin{pmatrix} \eta \\ 0 \end{pmatrix} , \begin{pmatrix} 0 \\ \eta \end{pmatrix} \right\} \;,
\end{equation}
The elements $\bp \in (\P_1(E))^2$ can be written in terms of the above basis as vectors $\hat{\bp} \in \real^6$. Moreover,  
since $(\P_1(E))^2 \subseteq \VE$, the elements $\bp \in (\P_1(E))^2$ can be written also in terms of the basis of $\VE$. 
We call $\bD$ the matrix associated to such change of basis, that can be simply calculated by evaluating the six functions in \eqref{eq:lin_pol_basis} through the degrees of freedom of $\VE$. Specifically, one has
\begin{equation}
\tilde{\bp} = \bD \,\hat{\bp}  
\label{eq:Duu}
\end{equation}
where $\tilde{\bp} \in \real^8$ represents $\bp$ written in terms of the $\VE$ basis. The matrix $\bD$ is given by
\begin{equation}
\bD= \left[ \begin{array}{cccccc}
1   &   0   &   \xi_1 &   0   &   \eta_1 &   0   \\
0   &   1   &   0   &   \xi_1 &   0   &   \eta_1 \\
1   &   0   &   \xi_2 &   0   &   \eta_2 &   0   \\
0   &   1   &   0   &   \xi_2 &   0   &   \eta_2 \\
1   &   0   &   \xi_3 &   0   &   \eta_3 &   0   \\
0   &   1   &   0   &   \xi_3 &   0   &   \eta_3 \\
1   &   0   &   \xi_4 &   0   &   \eta_4 &   0   \\
0   &   1   &   0   &   \xi_4 &   0   &   \eta_4
\end{array}
\right] \,,\label{eq:Dmatrix}
\end{equation}
where the $i-th$ vertex of the polygon has scaled coordinates $(\xi_i,\eta_i)$, $i=1,...,4$.

A stabilization strategy is obtained by considering the form $S^E(\u_h,\v_h) = \tilde{\u}_h^T \stiffness_{S} \tilde{\v}_h$ for all $\u_h \;, \v_h \in \VE$, taking:
\begin{equation}
\label{eq:stab_term}
\stiffness_{S} = \tau \: \tr\left(\stiffness_c\right) \left[ \bI - \bD \left(\bD^{\textrm{T}} \bD \right)^{-1} \bD^{\textrm{T}} \right] \;,
\end{equation}
with $\tau$ positive real number (see below).
The total stiffness operator results
\begin{equation}
\label{eq:final:stiff}
a^E_h(\u_h , \v_h) :=   \tilde{\u}_h^T  \stiffness  \tilde{\v}_h \quad \textrm{with}
\quad \stiffness = \stiffness_c + \stiffness_{S} .
\end{equation}
The role of the stabilization matrix $\stiffness_{S}$ is to keep the positivity of the discrete (local) energy form (up to the constant functions that should still give zero energy). Indeed, it can be checked that, thanks to the contribution of the matrix $\stiffness_{S}$, the discrete bilinear form satisfies the condition:
$$ 
\eps(\v_h) = {\bf 0} ,
$$
for any $\v_h \in \VE$ such that $a^E_h(\v_h , \v_h)  = 0$. Thus, a vanishing discrete elastic energy implies a rigid body motion.

The trace term in the definition of $\stiffness_{S}$ is added in order to guarantee the correct scaling of the energy with respect to the element size and material constants. The factor $\tau \in \mathbb R$, $\tau>0$, can be discarded but is included in order to allow for generalizations. Indeed, a typical good choice is simply $\tau=1$, but other choices can be found in the literature (see for instance \cite{Paulino-VEM}). As shown later in the numerical tests, the method is not much sensible to the parameter $\tau$. 
We refer to the initial paper \cite{volley} for more details on the motivations of the stabilization term, while a comparison with the standard ``single Gauss point'' quadrilateral FEM element and the related hourglass mode control can be found in \cite{hourglass}.

\begin{remark}\label{rem:patch:1} 
It is straightforward to check that, if $\tilde{\bp} \in \VE$ represents a polynomial (in other words satisfies \eqref{eq:Duu} for some $\hat{\bp} \in \real^6$) then $\stiffness_{S} \: \tilde{\bp}$ vanishes. Therefore, from the definition of $\stiffness$, it can be derived that (in the case of constant elastic coefficients) the local patch test is satisfied by the discrete bilinear form, i.e.:
$$
a^E_h(\bp , \v_h) = a^E (\bp , \v_h) \quad \forall \bp \in (\P_1(E))^2, \ \v_h \in \VE ,
$$
a property that in the Virtual Element terminology is called $k$-\emph{consistency} \cite{volley,hitchhikers}.
\end{remark}

\subsubsection{Loading term}
\label{sss:load_term}

The loading term for the case $k=1$ can be approximated by applying an integration rule based on vertexes. As in standard FE technology, the wirtual work exerted by the external loads is the element-wise sum of local contributions:
\begin{equation}\label{load-split}
<\bb_h,\v_h> = \sum_{E \in \Th} <\bb_h,\v_h>_E ,
\end{equation}
where we remind that $<\bb_h,\v_h>$, see also \eqref{gen-vem-str}, represents the VEM approximation of the right hand side, expressed below (for varying $\v_h$) as a linear operator on $\V_h$. 

Let the vector $\overline \bb \in \real^2$ be defined component-wise
$$
\overline b_j := \int_E b_j \quad j=1,2 ,
$$
with $b_1,b_2$ the two components of the volume force $\bb$.
Then, we can approximate the local loading term $\int_E \bb ^{\textrm{T}} \v_h$, for all $\v_h \in \VE$, with the following simple integration rule
\begin{equation}
\label{eq:loadapprox}
<\bb_h,\v_h>_E = \big(\!\!\int_E \!\! \bb \big)^{\textrm{T}} \big( \sum_{i=1}^4 \frac{1}{4} \v_h(\nu_i) \big) = 
\overline b_1 \sum_{i=1}^4 \frac{1}{4} (\tilde{\v}_h)_{2i-1} 
+ \overline b_2 \sum_{i=1}^4 \frac{1}{4} (\tilde{\v}_h)_{2i} .
\end{equation}
where $(\tilde{\v}_h)_{s}$, $s=1,...,8$, is the $s-th$ component of the vector $\tilde{\v}_h$, evaluated at the four corners $\nu_i$, $i=1,...,4$, of $E$. 

The loading term approximation reported in Eq. \eqref{eq:loadapprox} is exact whenever the load $\bb$ is constant over the element, and the displacement field $\v_h$ is a first order polynomial. 

\subsection{The case $k=2$ for general polygons}
\label{ss:ktwo}

The step from $k=1$ to $k=2$ is substantial, since in the case $k=2$ edge degrees of freedom and moments come into play in the approximation of the displacement field. We now consider the case of a general polygon with $m \ge 3$ edges.  

\begin{figure}
\begin{minipage}[b]{.5\linewidth}
\centering
\includegraphics[bb = 10 18 20 300, scale=0.5]{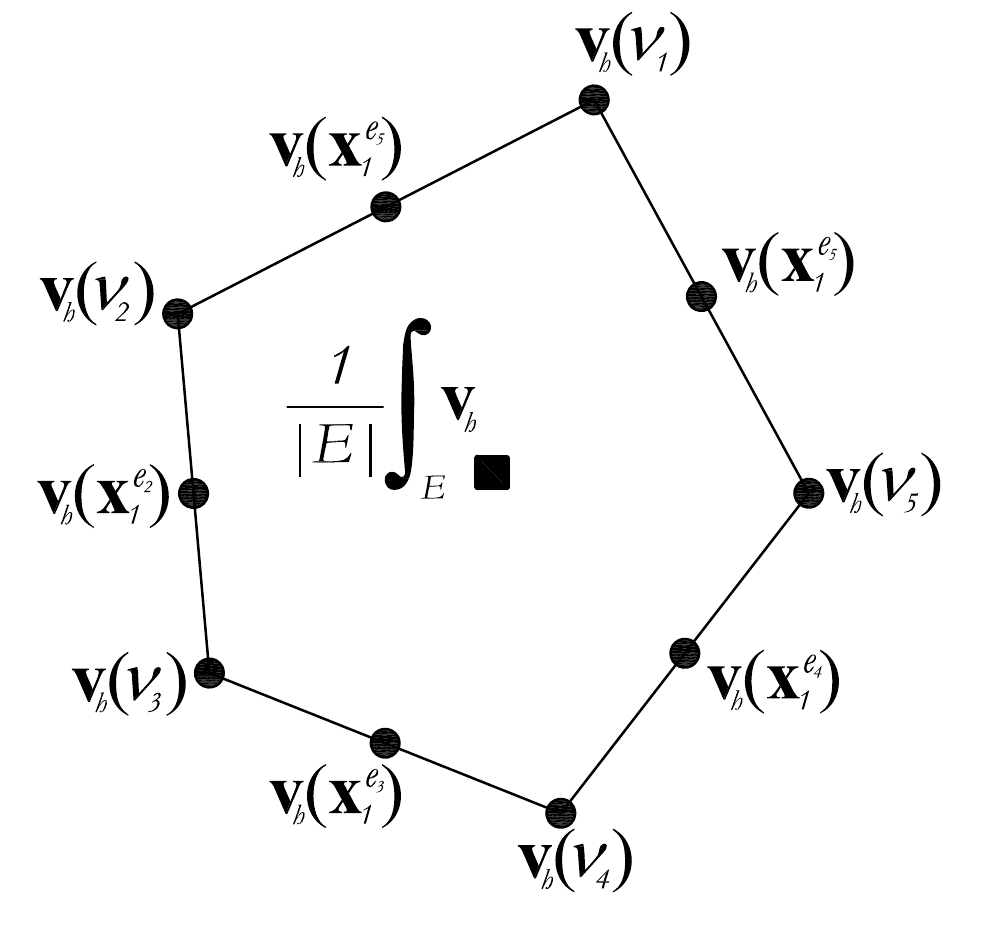}
\end{minipage}
\caption{Degrees of freedom for the case $k=2$, sample case on a pentagon. The degrees of freedom are the displacement values at the vertexes, the displacement values at the edge midpoints and the integral average of the displacement over the element.}
\label{fig:dof:k2}
\end{figure}

By referring to the list of degrees of freedom in Section \ref{s:one}, it is immediate to check that for $k=2$ there are 
\begin{itemize}
\item 2 (since the functions are vector valued) degrees of freedom per vertex, representing the displacement value at the vertex
\item 2 degrees of freedom per edge (representing the displacement value at the midpoint of the edge)
\item 2 scalar moments over the element
\end{itemize}
The degrees of freedom are depicted in figure \ref{fig:dof:k2}. Note that the edge midpoints are denoted by ${\boldsymbol{x}}_1^e$, $e \in \partial\E$.

The moments degrees of freedom are the component-wise average of the function $\v_h$
\begin{equation}
\label{eq:intdof}
\frac{1}{|E|} \int_E \v_h \;.
\end{equation}
The virtual element space and the polynomial space of approximated strains, from \eqref{dofsnumber} with $k=2$, satisfy:
\begin{equation}
\begin{aligned}
n = \textrm{dim}(\VE) = 4 m+2 \\
\ell =  \textrm{dim}(\P_1(\E)^{2\times 2}_{\textrm{sym}}) = 9 \;,
\end{aligned}
\end{equation}
and we set
\begin{equation}
\P_1(\E)^{2\times 2}_{\textrm{sym}} = \textrm{span}  \left\{ 
\begin{pmatrix} 1 \\ 0 \\ 0 \end{pmatrix},  \begin{pmatrix} 0 \\ 1 \\ 0 \end{pmatrix}, 
\begin{pmatrix} 0 \\ 0 \\ 1 \end{pmatrix}, 
\begin{pmatrix} \xi \\ 0 \\ 0 \end{pmatrix},  \begin{pmatrix} 0 \\ \xi \\ 0 \end{pmatrix}, 
\begin{pmatrix} 0 \\ 0 \\ \xi \end{pmatrix}, 
\begin{pmatrix} \eta \\ 0 \\ 0 \end{pmatrix},  \begin{pmatrix} 0 \\ \eta \\ 0 \end{pmatrix}, 
\begin{pmatrix} 0 \\ 0 \\ \eta \end{pmatrix}
\right\} \;.
\end{equation}
Representation \eqref{eq:u-approx} still applies: $\v_h \in \VE$, 
$\tilde{\v}_h \in \real^{4m+2}$, $\sgrad^{P} \in \P_1(\E)^{2\times 2}_{\textrm{sym}}$, $\hat{\sgrad} \in \real^9$. 
The matrix $\bN^P$ is given by
\begin{equation}
\bN^P=\left[
\begin{matrix}
1 & 0 & 0 & \xi & 0 & 0 & \eta & 0 & 0 \\
0 & 1 & 0 & 0 & \xi & 0 & 0 & \eta & 0 \\
0 & 0 & 1 & 0 & 0 & \xi & 0 & 0 & \eta
\end{matrix}
\right]
\end{equation}
The matrix operator $\bN^V$ representing the basis functions of $\VE$ is a $2 \times (4m+2)$ operator. Ordering the $m$ vertexes and the $m$ edge midpoint nodes of $E$ with $i=1,..,2m$, setting $s = 2i-1$, the $s^{th}$-column first entry is given by the unique function of $V_{h|E}$ that takes the value $1$ at $i-th$ location, vanishes at all other vertexes and midside nodes and has null average over the polygon (cf. \eqref{eq:intdof}); the second entry is zero. Analogously, setting $s = 2i$, the $s^{th}$-column has null first entry, being the second entry a basis function with the above mentioned properties. The subsequent $(4m+1)-th , \; (4m+2)-th$ columns contain basis functions associated with moments degrees of freedom, hence they vanish at vertexes and edge nodes, and have unit average over the element. 
Similarly to the $k=1$ case, we note that for $k=2$ the functions of $\VE$ are piecewise \emph{quadratic} and continuous on the boundary of $E$; therefore the ƒƒ $\bN^V$ can be easily computed on the element boundary.
Moreover, also the integral of those functions can be evaluated. Indeed, from the definition of the internal degrees of freedom one gets that
\begin{equation}\label{eq:integral}
\int_E \bN^V = 
\begin{pmatrix} 
0 & 0 & .... & 0 & |E| & 0 \\ 
0 & 0 & .... & 0 & 0 & |E|
\end{pmatrix}
\in \real^{2 \times (4m+2)} \;.
\end{equation}

\subsubsection{Consistent term}
\label{sss:cons_term:k2}
In the present section we describe the consistent part of the bilinear form $a^E_h(\cdot, \cdot)$, starting by
the practical construction of the projection operator $\Pi$ through the definition of the associated matrix 
$\Pi^{\textrm{m}}$
\begin{equation}
\Pi(\v_h) = \bN^P \Pi^{\textrm{m}} \tilde{\v}_h
\label{eq:pi_matr-2}
\end{equation}
with $\Pi^{\textrm{m}} \in \real^{9\times (4m+2)}$.
By inserting the above positions into Eq. \eqref{eq:def-proj} one again obtains \eqref{eq:X2}, that now holds for all 
$\hat{\sgrad}$ in $\real^9$. Integrating by parts, identity \eqref{eq:X2} yields
\begin{eqnarray}
\label{eq:basic_equation-2}
\nonumber
\int_{\E} \left(\bN^P \Pi^{\textrm{m}} \tilde{\v}_h \right)^{\textrm{T}} \bN^P \hat{\sgrad} = 
\int_{\E} \left[ \sgrad(\bN^V \tilde{\v}_h) \right]^{\textrm{T}} \bN^P \hat{\sgrad} = \\
\int_{\partial\E} \left(\bN^V \tilde{\v}_h \right)^{\textrm{T}} \left(\bN_E \bN^P \hat{\sgrad} \right) - \int_{\E} \left( \bN^V \tilde{\v}_h \right)^{\textrm{T}} \bpartial \bN^P \hat{\sgrad}
\end{eqnarray}
with the divergence matrix $\bpartial = \S^{\textrm{T}}$.

Eq. \eqref{eq:basic_equation-2} can be written as
\begin{equation}
\label{eq:compact:2}
\hat{\sgrad}^{\textrm{T}} \Pi^{\textrm{m}} \mathcal{G} \tilde{\v}_h  = 
\hat{\sgrad}^{\textrm{T}} \mathcal{B} \tilde{\v}_h \;,
\end{equation}
where
\begin{eqnarray}
\label{eq:XX3-2}
& \mathcal{G} & = \int_{\E} \left(\bN^P\right)^{\textrm{T}} \bN^P \in \real^{9 \times 9}\\
\label{eq:XX3-2b}
& \mathcal{B} & = \int_{\partial\E} \left(\bN_E \bN^P\right)^{\textrm{T}} \bN^V -
              \int_{\E} \left(\bpartial \bN^P\right)^{\textrm{T}} \bN^V \in \real^{9 \times (4m+2)} \;.
\end{eqnarray}
Therefore, the required operator $\Pi^{\textrm{m}}$ is given by the solution of a linear system:
\begin{equation}
\label{eq:Pi_operator-2}
\Pi^{\textrm{m}} = \mathcal{G}^{-1} \mathcal{B} \;.
\end{equation}

\begin{remark}
Note that matrix $\mathcal{B}$ is also computable. Indeed, the first term in the right hand side of \eqref{eq:XX3-2b} can be calculated: recall that the columns of $\bN^V$, representing the $\VE$ basis functions, are piecewise quadratic functions explicitly known on the boundary, see also Remark \ref{rem:edge-gauss} below. Regarding the second term in the right hand side of \eqref{eq:XX3-2b}, it can be immediately computed from the observation that $\bpartial \bN^P$ has constant entries
$$
\int_{\E} \left(\bpartial \bN^P\right)^{\textrm{T}} \bN^V =  \left(\bpartial \bN^P\right)^{\textrm{T}}  \int_{\E} \bN^V
$$
and using \eqref{eq:integral}.

\end{remark}

Combining eqs. \eqref{eq:pi_matr-2} and \eqref{eq:Pi_operator-2}, from \eqref{eq:approx_bilinear_form}, one obtains:
\begin{align}
\label{eq:approx_bilinear_form_Pi_bis}
\int_{\E} \left( \Pi(\v_h)\right)^{\textrm{T}} \bC  \Pi (\u_h)  = 
\int_{\E} \left( \bN^P \mathcal{G}^{-1} \mathcal{B} \tilde{\v}_h \right)^{\textrm{T}} \bC \left( \bN^P \mathcal{G}^{-1} \mathcal{B} \tilde{\u}_h \right) = 
\tilde{\v}^{\textrm{T}}_h \stiffness_{\textrm{c}} \tilde{\u}_h
\end{align}
which defines the consistent symmetric stiffness matrix:
\begin{equation}
\label{eq:consistent_stiffness-2}
\stiffness_{\textrm{c}} = \mathcal{B}^{\textrm{T}} \mathcal{G}^{\textrm{-T}} 
\left[ \int_{\E}  \left(\bN^P \right)^{\textrm{T}} \bC \bN^P \right ] 
\mathcal{G}^{-1} \mathcal{B}\;.
\end{equation}
\begin{remark}\label{rem:edge-gauss}
The boundary term in \eqref{eq:XX3-2b} is split as a sum of edge integrals involving the product of a polynomial of degree $\le 2$ and a polynomial of degree $\le 1$. Therefore, a Gauss-Lobatto integration rule with three nodes on the edge (that is, the endpoints and the midpoint) is sufficient to compute exactly such integrals. This means that in practice one needs to compute the value of the matrix $\bN^V$ only at the element vertexes and element edge midpoints. Due to our choice of the basis functions, and since such points are associated to the degrees of freedom of the space $\VE$, the matrix $\bN^V$ always takes either the value $0$ or the  value $1$ at the vertexes and the edge midpoints. This observation clearly simplifies the computation of matrix $\mathcal{B}$.
\end{remark}

\subsubsection{Stabilization term}
\label{sss:stab_term-2}

Consider the following basis for the space $(\P_2(E))^2$:
\begin{equation}
\label{eq:lin_pol_basis:2}
\begin{aligned}
(\P_2(\E))^2 = \textrm{span}  \Big\{ & \begin{pmatrix} 1 \\ 0 \end{pmatrix} ,  \begin{pmatrix} 0 \\ 1 \end{pmatrix}, \begin{pmatrix} \xi \\ 0 \end{pmatrix} 
, \begin{pmatrix} 0 \\ \xi \end{pmatrix} , \begin{pmatrix} \eta \\ 0 \end{pmatrix} , \begin{pmatrix} 0 \\ \eta \end{pmatrix}, 
\\
& \begin{pmatrix} \xi^2 \\ 0 \end{pmatrix} 
, \begin{pmatrix} 0 \\ \xi^2 \end{pmatrix} , 
\begin{pmatrix} \xi\eta \\ 0 \end{pmatrix} , 
\begin{pmatrix} 0 \\ \xi\eta \end{pmatrix} ,
\begin{pmatrix} \eta^2 \\ 0 \end{pmatrix} , 
\begin{pmatrix} 0 \\ \eta^2 \end{pmatrix} 
\Big\} \;,
\end{aligned}
\end{equation}
where we recall that we are using local scaled coordinates (cf. Eq. \eqref{eq:adimcoor}).
Similarly to the case $k=1$, the elements $\bp \in (\P_2(E))^2$ can be written in terms of the above basis as vectors $\hat{\bp} \in \real^{12}$ and,  
since $(\P_2(E))^2 \subseteq \VE$, also in terms of the basis of $\VE$. 
We call $\bD$ the matrix associated to such change of basis, that can be simply calculated by evaluating the functions in \eqref{eq:lin_pol_basis:2} through the degrees of freedom of $\VE$. 
One has
\begin{equation}
\tilde{\bp} = \bD \,\hat{\bp}  
\label{eq:Duu:2}
\end{equation}
where $\tilde{\bp} \in \real^{(4m+2)}$ represents $\bp$ written in terms of the $\VE$ basis. 
The matrix $\bD$ is given by 
\begin{equation}
\bD= \left[ \begin{array}{cccccccccccc}
1   &   0   &   \xi_1 &   0   &   \eta_1 &   0   &   \xi_1^2 &   0 & \xi_1 \eta_1 &   0  &   \eta_1^2 &   0     \\
0   &   1   &   0   &   \xi_1 &   0   &   \eta_1 &   0   &   \xi_1^2 & 0   &   \eta_1 \xi_1  &   0   &   \eta_1^2 \\
1   &   0   &   \xi_2 &   0   &   \eta_2 &   0   &   \xi_2^2 &   0  & \xi_2 \eta_2 &   0 &   \eta_2^2 &   0     \\
0   &   1   &   0   &   \xi_2 &   0   &   \eta_2 &   0   &   \xi_2^2 & 0   &   \eta_2 \xi_2 &   0   &   \eta_2^2  \\
\vdots   &   \vdots   &   \vdots &   \vdots   &   \vdots &   \vdots   &   \vdots &   \vdots  &   \vdots &   \vdots   & \vdots &   \vdots  \\
1   &   0   &   \xi_m &   0   &   \eta_m &   0   &   \xi_m^2 &   0 & \xi_m \eta_m &   0  &   \eta_m^2 &   0     \\
0   &   1   &   0   &   \xi_m &   0   &   \eta_m &   0   &   \xi_m^2 & 0   &   \eta_m \xi_m &   0   &   \eta_m^2  \\
1   &   0   &   \zeta_1 &   0   &   \omega_1 &   0   &   \zeta_1^2 &   0   & \zeta_1 \omega_1 &   0 &   \omega_1^2 &   0     \\
0   &   1   &   0   &   \zeta_1 &   0   &   \omega_1 &   0   &   \zeta_1^2 & 0   &   \omega_1 \zeta_1 &   0   &   \omega_1^2  \\
1   &   0   &   \zeta_2 &   0   &   \omega_2 &   0   &   \zeta_2^2 &   0  & \zeta_2 \omega_2 &   0 &   \omega_2^2 &   0     \\
0   &   1   &   0   &   \zeta_2 &   0   &   \omega_2 &   0   &   \zeta_2^2 & 0   &   \omega_2 \zeta_2 &   0   &   \omega_2^2  \\
\vdots   &   \vdots   &   \vdots &   \vdots   &   \vdots &   \vdots   &   \vdots &   \vdots  &   \vdots &   \vdots   & \vdots &   \vdots  \\
1   &   0   &   \zeta_m &   0   &   \omega_m &   0   &   \zeta_m^2 &   0   & \zeta_m \omega_m &   0 &   \omega_m^2 &   0     \\
0   &   1   &   0   &   \zeta_m &   0   &   \omega_m &   0   &   \zeta_m^2 & 0   &   \omega_m \zeta_m &   0   &   \omega_m^2  \\
1   &   0   &   \strokedint \xi &   0    &   \strokedint \eta &   0      &   \strokedint \xi^2 &   0     &   \strokedint \xi\eta &   0     &   \strokedint \eta^2 &   0       \\
0   &   1   &  0       & \strokedint \xi &   0      &   \strokedint \eta &   0          &   \strokedint \xi^2 &   0        &   \strokedint \xi\eta &   0     &   \strokedint \eta^2  
\end{array}
\right] \,,
\label{eq:Dmatrix:2}
\end{equation}
where the $i-th$ vertex of the polygon is supposed to have scaled coordinates $(\xi_i\; \eta_i)$, $i=1,...,m$, the $i-th$ edge midpoint is assumed to have coordinates $(\zeta_i\;\omega_i)$, $i=1,...,m$, and we denote with the short symbol $\strokedint p$ the integral average of any scalar polynomial $p$ on the polygon $E$. 

Once the matrix $\bD$ is built, the rest of the procedure follows exactly the same line of the case  $k=1$, namely Eq. \eqref{eq:stab_term} for the construction of the stabilizing matrix, and Eq. \eqref{eq:final:stiff} for the final stiffness matrix.
\begin{remark}
Analogously to Remark \ref{rem:patch:1}, it is easy to check that in this case the following ``patch test'' condition holds (in the case of constant elastic coefficients)
$$
a^E_h(\bp , \v_h) = a^E (\bp , \v_h) \quad \forall \bp \in (\P_2(E))^2, \ \forall \v_h \in \VE .
$$
\end{remark}

\subsubsection{Loading term}
\label{sss:load_term-2}

The loading term for the case $k=2$ is simpler than $k=1$, due to the presence of the moment degrees of freedom. Indeed, after the usual element-wise splitting \eqref{load-split}, we approximate the local loading term by
$$
<\bb_h,\v_h>_E := \frac{1}{|E|} \left( \int_E \bb \right) ^{\textrm{T}} \left( \int_E \v_h \right) ,
$$
that corresponds to an approximation of the local load by its average.
Following the notation of Section \ref{sss:load_term}, this yields
$$
<\bb_h,\v_h>_E = \frac{1}{|E|} \: \overline\bb^{\textrm{T}} \left( \int_E \bN^V \right) \tilde{\v}_h ,
$$
that can be immediately computed using \eqref{eq:integral}.

\subsection{The case of general order $k$}
\label{ss:kgen}

The degrees of freedom for the case of general order $k$ are given in Section \ref{s:abs}. In particular, for any edge of a polygon, edge degrees of freedom are interpolatory displacements at points $\{\bx^e_i \}$, $i=1,...,k-1$, chosen as the internal nodes of the Gauss-Lobatto integration rule with $k+1$ nodes. For instance, for $k=3$ the Gauss-Lobatto integration rule leads to $k+1=4$ nodes, i.e. $\{\bx^e_i \}$, $i=1,2$ internal nodes, plus $2$ extrema. In passing, we note that, for $k=2$, this scheme provides the edge midpoint of the previous section. 
The definition of the $k(k-1)$ internal degrees of freedom is set by considering the following basis for $\P_{k-2}(E)$
\begin{equation}
\label{eq:Pkm2basis}
\left\{ \xi^{\alpha} \eta^{\beta} \ \textrm{ with } \ \alpha,\beta \in {\mathbb N}, \: \alpha + \beta \le k-2 
\right\} 
\end{equation}
The monomials of basis \eqref{eq:Pkm2basis} are ordered as the exponent couples $(\alpha \; \beta)$ of Pascal's triangle and named accordingly $q_1,q_2,...,q_r$, with $r=k(k-1)/2$ the dimension of $\P_{k-2}(E)$. This ordered basis results: 
\begin{equation}
\label{eq:orderedbasis}
\{ q_1,q_2,...,q_r \} \:  = \: \Big\{ 1 , \xi , \eta , \: \xi^2, \xi\eta, \xi \eta^2, \: \xi^3, \xi^2 \eta, \xi \eta^2, \eta^3 , ..... , \xi^{k-2}, \xi^{k-3} \eta, ...., \eta^{k-2} \Big\} .
\end{equation}
The moments of the displacement components, playing the role of internal degrees of freedom, are the $2r$ linear operators from $\VE$ into $\real$ defined as:
\begin{equation}\label{int-dofs}
\begin{aligned}
& \Xi^\textrm{int}_{2j-1} (\v_h) = \frac{1}{|E|} \int_E q_j \: u_h \quad  j = 1,2,...,r \: , \\
& \Xi^\textrm{int}_{2j} (\v_h) = \frac{1}{|E|} \int_E q_j \: v_h \quad j = 1,2,...,r
\end{aligned}
\end{equation}
for all $\v_h = (u_h\;v_h)^T \in \VE$. The internal degrees of freedom are simply the scaled moments of the function $\v_h$ (for the first and second component) up to order $k-2$. 
For example, for $k=3$, the $2r=6$ internal degrees of freedom are:
$$
\Big\{ 
\frac{1}{|E|} \int_E u_h , 
\frac{1}{|E|} \int_E v_h , 
\frac{1}{|E|} \int_E \xi \: u_h ,
\frac{1}{|E|} \int_E \xi \: v_h ,
\frac{1}{|E|} \int_E \eta \: u_h ,
\frac{1}{|E|} \int_E \eta \: v_h 
\Big\}.
$$

The dimensions of the virtual element space and of the polynomial space of approximated strains (cf. Eqs. \eqref{dofsnumber}, \eqref{uVuP}) are, respectively:
\begin{equation}
\begin{aligned}
& n := \textrm{dim}(\VE) = 2km + 2r \\
& \ell : = \textrm{dim}(\P_{k-1}(\E)^{2\times 2}_{\textrm{sym}}) = 3k(k+1)/2 \; .
\end{aligned}
\end{equation}
We order the basis for $\P_{k-1}(\E)^{2\times 2}_{\textrm{sym}}$ coherently with the ordering of $q_i$'s:
\begin{equation}\label{k-1:basis}
\begin{aligned}
\P_{k-1}(\E)^{2\times 2}_{\textrm{sym}} = \textrm{span}  & \left\{ 
\begin{pmatrix} 1 \\ 0 \\ 0 \end{pmatrix},  \begin{pmatrix} 0 \\ 1 \\ 0 \end{pmatrix}, 
\begin{pmatrix} 0 \\ 0 \\ 1 \end{pmatrix}, 
\begin{pmatrix} \xi \\ 0 \\ 0 \end{pmatrix},  \begin{pmatrix} 0 \\ \xi \\ 0 \end{pmatrix}, 
\begin{pmatrix} 0 \\ 0 \\ \xi \end{pmatrix}, 
\begin{pmatrix} \eta \\ 0 \\ 0 \end{pmatrix},  \begin{pmatrix} 0 \\ \eta \\ 0 \end{pmatrix}, 
\begin{pmatrix} 0 \\ 0 \\ \eta \end{pmatrix} , \right. \\
&  
\begin{pmatrix} \xi^2 \\ 0 \\ 0 \end{pmatrix},  \begin{pmatrix} 0 \\ \xi^2 \\ 0 \end{pmatrix}, 
\begin{pmatrix} 0 \\ 0 \\ \xi^2 \end{pmatrix}, 
\begin{pmatrix} \xi\eta \\ 0 \\ 0 \end{pmatrix},  \begin{pmatrix} 0 \\ \xi\eta \\ 0 \end{pmatrix}, 
\begin{pmatrix} 0 \\ 0 \\ \xi\eta \end{pmatrix}, 
\begin{pmatrix} \eta^2 \\ 0 \\ 0 \end{pmatrix},  \begin{pmatrix} 0 \\ \eta^2 \\ 0 \end{pmatrix}, 
\begin{pmatrix} 0 \\ 0 \\ \eta^2 \end{pmatrix} , \\
& \left.
\begin{pmatrix} \xi^3 \\ 0 \\ 0 \end{pmatrix},  \begin{pmatrix} 0 \\ \xi^3 \\ 0 \end{pmatrix}, 
\begin{pmatrix} 0 \\ 0 \\ \xi^3 \end{pmatrix}, 
............. , 
\begin{pmatrix} \eta^{k-1} \\ 0 \\ 0 \end{pmatrix},  \begin{pmatrix} 0 \\ \eta^{k-1} \\ 0 \end{pmatrix}, 
\begin{pmatrix} 0 \\ 0 \\ \eta^{k-1} \end{pmatrix}
\right\} \;.
\end{aligned}
\end{equation}
Representation \eqref{eq:u-approx} still applies, where $\v_h \in \VE$, $\tilde{\v}_h \in \real^n$, $\sgrad^{P} \in \P_{k-1}(\E)^{2\times 2}_{\textrm{sym}}$, $\hat{\sgrad} \in \real^\ell$. 
The matrix $\bN^P \in \real^{3 \times \ell}$ is immediately derived from the basis in \eqref{k-1:basis}
\begin{equation}
\bN^P=\left[
\begin{array}{ccccccccccccc}
1 & 0 & 0 & \xi & 0 & 0 & \eta & 0 & 0 & ...... & \eta^{k-1} & 0  & 0  \\
0 & 1 & 0 & 0 & \xi & 0 & 0 & \eta & 0 & ...... & 0 & \eta^{k-1} & 0   \\
0 & 0 & 1 & 0 & 0 & \xi & 0 & 0 & \eta & ...... & 0 & 0           & \eta^{k-1} \\
\end{array}
\right] .
\end{equation}

$\bN^V$ is a $2 \times n$ matrix operator representing the basis function of $\VE$; odd columns have zero second entry, even columns have the reverse property. Each shape function takes the value $1$ for the corresponding degree of freedom and vanishes for all the others. The ordering of the degrees of freedom in vector $\tilde{\v}_h$ (and associated basis functions) is the same one for the case $k=2$, for $\tilde{u}_h$ and $\tilde{v}_h$ displacement components, respectively. Such indexing depends on the usual indexing of vertexes, edge nodes, and on the ordering of basis functions \eqref{eq:orderedbasis}. It is reported in Table \ref{T:dofs} for compactness. 



\begin{table}[!htbp]
\centering
\caption{Indexing $s$ of vector $\tilde{\v}_h$ components $(\tilde{u}_{h,s}$\;$\tilde{v}_{h,s})$ for the case with general order $k$, and polygon with $m$ edges.}
\begin{tabular}{lccccc}
\toprule
&  \multicolumn{2}{c}{$\tilde{u}_{h,s}$} & \multicolumn{2}{c}{$\tilde{v}_{h,s}$}\\
&  \multicolumn{2}{c}{$(s=2i-1)$} & \multicolumn{2}{c}{$(s=2i)$}\\
\midrule
Vertex $(i)$    &  &  \multicolumn{2}{l}{$i=1,...,m$}    &  \\
Edge $(i)$      &  &  \multicolumn{2}{l}{$i=m+1,...,km$} &  \\
Internal $(i)$  &  &  \multicolumn{2}{l}{$i=km+1,...,km+r$} &  \\
\bottomrule
\end{tabular}
\label{T:dofs} 

\end{table}

\begin{remark}\label{rem:k:int}
The basis functions associated to the internal degrees of freedom clearly vanish on the boundary, since by definition these functions take zero value at all the boundary nodes. 
The internal basis functions are ordered in the same way as the polynomials $\{q_1,q_2,...,q_r \}$ introduced above. Namely, for $j=1,...,r$ and $i=km+1,...,km+r$ (that is the index range associated to the internal basis functions accordingly to the ordering in Table \ref{T:dofs}) we have for the first component
\begin{equation}
\label{eq:piropo1}
\begin{aligned}
& \frac{1}{|E|} \int_E q_j \bN^V_{1,2i-1} =  \delta_{i,j} \\
& \frac{1}{|E|} \int_E q_j \bN^V_{2,2i-1} =  0 ,
\end{aligned}
\end{equation}
with $\delta_{ij}$ representing the Kronecker symbol. The internal basis functions associated to the second component follow an analogous rule. For $j=1,...,r$ and $i=km + 1,...,km+r$, it holds:
\begin{equation}
\label{eq:piropo2}
\begin{aligned}
& \frac{1}{|E|} \int_E q_j \bN^V_{1,2i} =  0\\
& \frac{1}{|E|} \int_E q_j \bN^V_{2,2i} =  \delta_{i,j} ,
\end{aligned}
\end{equation} 
\end{remark}

Similarly to the previous cases, the functions of $\VE$ on the boundary of $E$ are continuous and piecewise polynomials of degree $k$. Therefore the evaluation of $\bN^V$ can be easily computed on the element boundary.
Moreover, the integral on $E$ of the basis functions multiplied by any polynomial $q_j$ is computable immediately by using Eqs. \eqref{eq:piropo1}-\eqref{eq:piropo2}. This allows to deduce that, for all $j=1,2,..,r$, the integral matrix
\begin{equation}\label{eq:intk:fund}
\int_E q_j \: \bN^V \: \in \real^{2 \times n}
\end{equation}
has only two non-vanishing entries, with both value $|E|$, at position $(1,2km+2j-1)$ and $(2,2km+2j)$, respectively. Note that the first $2km$ columns of zeros are associated with the boundary degrees of freedom.
For instance, in the case $k=3$ (recalling that $q_1=1, q_2 = \xi, q_3 = \eta$) one gets
$$
\begin{aligned}
& \int_E \bN^V = 
\begin{pmatrix} 
0 & 0 & .... & 0 & |E| & 0 & 0 & 0 & 0 & 0 \\ 
0 & 0 & .... & 0 & 0 & |E| & 0 & 0 & 0 & 0
\end{pmatrix}
\in \real^{2 \times n} \; , \\
& \int_E \xi \: \bN^V = 
\begin{pmatrix} 
0 & 0 & .... & 0 & 0 &  0  &|E|& 0 & 0   & 0\\ 
0 & 0 & .... & 0 & 0 & 0   & 0 &|E| & 0  & 0
\end{pmatrix}
\in \real^{2 \times n} \; , \\
& \int_E \eta \: \bN^V = 
\begin{pmatrix} 
0 & 0 & .... & 0 & 0 & 0 &  0 & 0 & |E| & 0 \\ 
0 & 0 & .... & 0 & 0 & 0 & 0   & 0 & 0 & |E|
\end{pmatrix}
\in \real^{2 \times n} \; .
\end{aligned}
$$

\subsubsection{Consistent term}
\label{sss:cons_term:k_gen}
The practical construction of the projection operator $\Pi$, through the definition of the associated matrix 
$\Pi^{\textrm{m}}$
\begin{equation}
\Pi(\v_h) = \bN^P \Pi^{\textrm{m}} \tilde{\v}_h , \ \Pi^{\textrm{m}} \in \real^{\ell \times n} ,
\label{eq:pi_matr-k}
\end{equation}
follows the same steps shown for the $k=2$ case in Section \ref{ss:ktwo}.
Indeed, by inserting the above definitions into Eq. \eqref{eq:def-proj}, one again obtains \eqref{eq:X2}, that now holds for all $\hat{\sgrad}$ in $\real^\ell$. Integrating by parts, identity \eqref{eq:X2} again gives \eqref{eq:basic_equation-2}, that in turn still can be written as in \eqref{eq:compact:2} and thus yields the fundamental Eq. \eqref{eq:Pi_operator-2}. Clearly, now the dimensions involved are different:
\begin{eqnarray}
&& \mathcal{G} = \int_{\E} \left(\bN^P\right)^{\textrm{T}} \bN^P \in \real^{\ell \times \ell} \\
\label{eq:XX3-2b-k}
&& \mathcal{B} = \int_{\partial\E} \left(\bN_E \bN^P\right)^{\textrm{T}} \bN^V -
              \int_{\E} \left(\bpartial \bN^P\right)^{\textrm{T}} \bN^V \in \real^{\ell \times n} \; .
\end{eqnarray}
Finally, the consistent stiffness matrix $\stiffness_{\textrm{c}}$ can be computed as in \eqref{eq:consistent_stiffness-2}.

We close this section with important observations regarding the construction of matrix $\mathcal{B}$. 

The first term on the right hand side of \eqref{eq:XX3-2b-k} can be calculated recalling that the columns of $\bN^V$ (representing the basis functions of $\VE$) are piecewise polynomial functions of degree $\le k$ that are explicitly known on the boundary. Moreover, such calculation can be greatly simplified if one adopts as edge degrees of freedom the $(k-1)$ internal points associated to the Gauss-Lobatto integration rule, as already suggested.
Indeed, the boundary term in \eqref{eq:XX3-2b-k} is split as a sum of edge integrals involving the product of a polynomial of degree $\le k$ and a polynomial of degree $(k-1)$. Therefore, a Gauss-Lobatto integration rule with $(k+1)$ nodes on the edge (that is, the vertex nodes and the $(k-1)$ internal edge nodes) is sufficient to compute exactly such integrals. This means that in practice one needs to compute the value of the matrix $\bN^V$ only at the element vertexes and element edge nodes $\{\bx^e_i \}$, $i=1,...,(k-1)$. Due to our choice of the basis functions, and since such points are associated with the degrees of freedom of the space $\VE$, the matrix $\bN^V$ always takes either the value $0$ or the  value $1$ at the vertexes and the edge nodes. This observation clearly simplifies the computation of matrix $\mathcal{B}$.

Regarding the second term in the right hand side of \eqref{eq:XX3-2b-k}, it can be computed by noting that the entries of $\bpartial \bN^P$ are polynomials of degree $\le (k-2)$ and using \eqref{eq:intk:fund}.  
First of all, one needs to express the components of the polynomial matrix $\bpartial \bN^P$ in terms of the $q_j$ basis, introduced previously in the present section. Indeed, one can write
\begin{equation}\label{oneabove}
\bpartial \bN^P = \sum_{j=1}^r M^j q_j
\end{equation}
where the matrices $M^j \in \real^{2\times\ell}$ have constant entries. Note that, given $k$, the matrices $M_j$ can be computed offline once and for all with symbolic calculus tools, since \eqref{oneabove} is a polynomial identity and does not depend on the element $E$ (see also Remark \ref{comp:Mj} below).
Having the above matrices, one can simply write
$$
\int_{\E} \left(\bpartial \bN^P \right)^{\textrm{T}} \bN^V = \sum_{j=1}^r (M^j)^T \int_E q_j \bN^V , 
$$
and make use of \eqref{eq:intk:fund}. 

\begin{remark}\label{comp:Mj}
Note that the matrix $\int_E q_j \bN^V$ in \eqref{eq:intk:fund} is a sparse matrix (only two components are different from zero) and therefore the calculations above can be done at minimal computational cost. Moreover, the matrix $ \int_E q_j \bN^V$ can be immediately written, for every element $E$, as a sparse matrix with only two unit entries, multiplied by the factor $|E|$. Therefore, the product matrices $(M^j)^T \int_E q_j \bN^V$ have the following properties: (1) only two columns are different from zero; (2) they can be computed and stored offline very easily, since the only factor depending on the element is the scalar factor $|E|$, which can be included later on, during each element stiffness matrix construction.
\end{remark}

\subsubsection{Stabilization term}
\label{sss:stab_term-k}

The guidelines of the procedure are the same as in the cases $k=1$ and $k=2$.
Therefore, since the construction of the stabilization term has been already explained with full details in Section \ref{sss:stab_term-2}, we here limit to a compact description of the involved terms.

We start by introducing a basis $\{ {\mathsf p}_j \}$, $j=1,...,(k+1)(k+2)$, for the space $(\P_{k}(E))^2$. We set
$$
\left\{
\begin{aligned}
& {\mathsf p}_{2j-1} = \begin{pmatrix} q_j \\ 0 \end{pmatrix} \ \textrm{ for } j=1,2,...,(k+1)(k+2)/2 , \\
& {\mathsf p}_{2j} = \begin{pmatrix} 0 \\ q_j \end{pmatrix} \ \textrm{ for } j=1,2,...,(k+1)(k+2)/2 ,
\end{aligned}
\right.
$$
where the scaled monomials $q_j$ have been introduced at the beginning of Section \ref{ss:kgen} (cf. \eqref{eq:Pkm2basis}). Note that for $k=2$ we obtain the same identical basis shown in \eqref{eq:lin_pol_basis:2}.

Similarly to the previous cases, the elements $\bp \in (\P_k(E))^2$ can be written in terms of the above basis as vectors $\hat{\bp} \in \real^{(k+1)(k+2)}$ and, since $(\P_k(E))^2 \subseteq \VE$, also in terms of the basis of $\VE$. 
As usual we call $\bD$ the matrix associated to such change of basis, that can be simply computed by evaluating the polynomials $\{ {\mathsf p}_j \}$, $j=1,...,(k+1)(k+2)$, through the degrees of freedom of $\VE$. 
It holds:
\begin{equation}
\tilde{\bp} = \bD \,\hat{\bp}  \ , \quad \bD \in \real^{n \times (k+1)(k+2)}
\label{eq:Duu:k}
\end{equation}
where $\tilde{\bp} \in \real^{n}$ represents $\bp$ written in terms of the $\VE$ basis. 
The components of the matrix $\bD$ are given by
$$
\bD_{ij} = `` \textrm{evaluation of }  {\bf p}_j \textrm{ on the } i^{th} \textrm{ degree of freedom of } \VE''
$$
and can be calculated as follows.

Let $\{ \nu_i \}$, $i=1,...,m$ be the set of the polygon vertexes, and $\{ \bx^e_i \}$, $i=(m+1),...,km$ the edge nodes, associated with the space $\VE$. These are assumed to be ordered as described in Section \ref{ss:kgen}, i.e. first the vertex nodes and then the edge nodes (given some ordering of the edges and the local ordering of the Gauss-Lobatto nodes along each edge). Thus, the first $2km$ rows of the matrix $\bD$ are given by
$$
\left\{
\begin{aligned}
& \bD_{2i-1,2j-1} = {q}_j (\nu_i), \quad i=1,2,...,m , \ j =1,2,..,(k+1)(k+2)/2 \: , \\
& \bD_{2i,2j} = {q}_j (\nu_i), \quad i=1,2,...,m , \ j =1,2,..,(k+1)(k+2)/2 \: , \\
& \bD_{2i-1,2j-1} = {q}_j ( \bx^e_i), \quad i=m+1,m+2,...,km , \ j =1,2,..,(k+1)(k+2)/2 \: , \\
& \bD_{2i,2j} = {q}_j ( \bx^e_i), \quad i=m+1,m+2,...,km , \ j =1,2,..,(k+1)(k+2)/2 \: , 
\end{aligned}
\right.
$$
and all the remaining components of the first $2km$ rows are zero.
It may be helpful for the reader to give a look at the explicit matrix representation given in \eqref{eq:Dmatrix:2} for the case $k=2$. 

The remaining $2r=k(k-1)$ rows of $\bD$ are given by
$$
\left\{
\begin{aligned}
& \bD_{2i-1, 2j-1} = |E|^{-1} \int_E q_{i} {q}_j ,\quad i=km+1,2,...,n , \ j =1,2,..,(k+1)(k+2)/2 \: , \\
& \bD_{2i  , 2j} = |E|^{-1} \int_E q_{i} {q}_j ,\quad i=km+1,2,...,n , \ j =1,2,..,(k+1)(k+2)/2 \: , \\
\end{aligned}
\right.
$$
and all the remaining components of the last $2r$ rows are zero.

Once the matrix $\bD$ is built, the rest of the construction is identical to the previous cases, namely Eq. \eqref{eq:stab_term} is used for the construction of the stabilizing matrix, and Eq. \eqref{eq:final:stiff} for the final stiffness matrix.
Analogously to Remark \ref{rem:patch:1}, it is easy to check that in this case the following ``patch test'' condition holds (for constant elastic coefficients)
$$
a^E_h(\bp , \v_h) = a^E (\bp , \v_h) \quad \forall \bp \in (\P_k(E))^2, \ \forall \v_h \in \VE .
$$
\subsubsection{Loading term}
\label{sss:load_term-k}

As usual, we start by the element-wise splitting \eqref{load-split}, and we approximate the local loading term $\int_{E} \v^T \bb$. 
The first step is to introduce $\bb_h \in (\P_{k-2}(E))^2$ as the $L^2$ projection of the load $\bb|_E$ onto the polynomials of degree $(k-2)$. The second step is simply to define
\begin{equation}\label{load-k-basic}
<\bb_h,\v_h>_E \: = \int_E \v_h^T \bb_h \quad \forall \v_h \in \VE.
\end{equation}
We now describe such steps more in detail. Since $\bb_h\in(\P_{k-2}(E))^2$,  $\bb_h$ can be expressed in terms of the basis \eqref{eq:orderedbasis}:
\begin{equation}\label{qi_expr}
\bb_h = \sum_{i=1}^r \: \widetilde{\bb}_h^i \: q_i \ , \quad \widetilde{\bb}_h^i \in \real^2 .
\end{equation}
The collection of the above real numbers leads to the following vector $\widetilde{\bb}_{h} \in \real^{2r}$:
\begin{equation}
\label{eq:loadcoeff}
\widetilde{\bb}_{h} = \big( \widetilde{b}_{h,1}^1  \widetilde{b}_{h,2}^1 ...  \widetilde{b}_{h,1}^r   \widetilde{b}_{h,2}^r \big)  ,
\end{equation}
By definition of $L^2$ projection, and using \eqref{qi_expr}, it is easy to check that the coefficient vectors in \eqref{qi_expr} must satisfy
\begin{equation}
\label{eq:loadcoeffsystem}
\sum_i \int_E \widetilde{\bb}_{h}^i \: q_i \: q_j = \int_E \bb \: q_j , \quad j=1,2,..., r.
\end{equation}
The above linear system can be explicitly written for the coefficient vector \eqref{eq:loadcoeff}, as
\begin{equation}
\label{eq:loadsystem}
\bQ \: \widetilde{\bb}_{h} = \f \ , \quad \bQ \in \real^{2r \times 2r} , \ \f \in \real^{2r}\; ,
\end{equation}
where
\begin{equation}
\begin{aligned}
& \bQ_{2i-1,2j-1} = \bQ_{2i,2j} = \int_E q_i \: q_j \ , \quad i,j =1,2,...,r \; , \\
& \f_{2i-1} =   \int_E b_1 \: q_i   \ , \quad i=1,2,...,r \; ,  \\
& \f_{2i} = \int_E b_2 \: q_i \ , \quad i=1,2,...,r \; .
\end{aligned}
\end{equation}

Finally, once the $\widetilde{\bb}_h$ vectors is computed as shown above, using \eqref{load-k-basic} and \eqref{qi_expr}, we have
$$
<\bb_h,\v_h>_E \: = \sum_{i=1}^r \widetilde{\bb}_h^i \int_E q_i \bN^V ,
$$
that is immediately computed recalling \eqref{eq:intk:fund}. The same observations of Remark \ref{comp:Mj} still apply.

\section{Numerical tests}
\label{s:num_test}
The proposed VEM methods are implemented into in-house Matlab \cite{matlab:2015} codes for validation and analysis. In this section we assess accuracy and robustness through a number of representative boundary value problems on a unit square domain plus a classical 2D problem taken from the literature. They are organized as follows:
\begin{itemize}
\item[$\bullet$]
Tests $1a$, $1b$: 2D plane strain elasticity linear patch tests on the unit square $\Omega = [0,1]^2$;
\item[$\bullet$]
Tests $2a$, $2b$: 2D plane strain convergence tests with known analytical solution on the unit square $\Omega = [0,1]^2$;
\item[$\bullet$]
Test $3$: assessment of the influence of the stabilization parameter;
\item[$\bullet$]
Test $4$: Cook's Membrane.
\end{itemize}
For the sake of simplicity, the geometric and physical quantities are set omitting the system of measurement; it is nonetheless assumed that the relevant units are taken consistently.
In all the numerical tests we use the choice $\tau=1/2$ for the stabilization parameter, see \eqref{eq:stab_term}. In section \ref{ss:test_3} we show that taking other choices has limited influence on the results.
 
\subsection{Test $1a$ and $1b$}
\label{ss:test_1}
\begin{figure}
\begin{minipage}[b]{.5\linewidth}
\centering
\includegraphics[bb = 150 18 185 320, scale = 0.5]
                {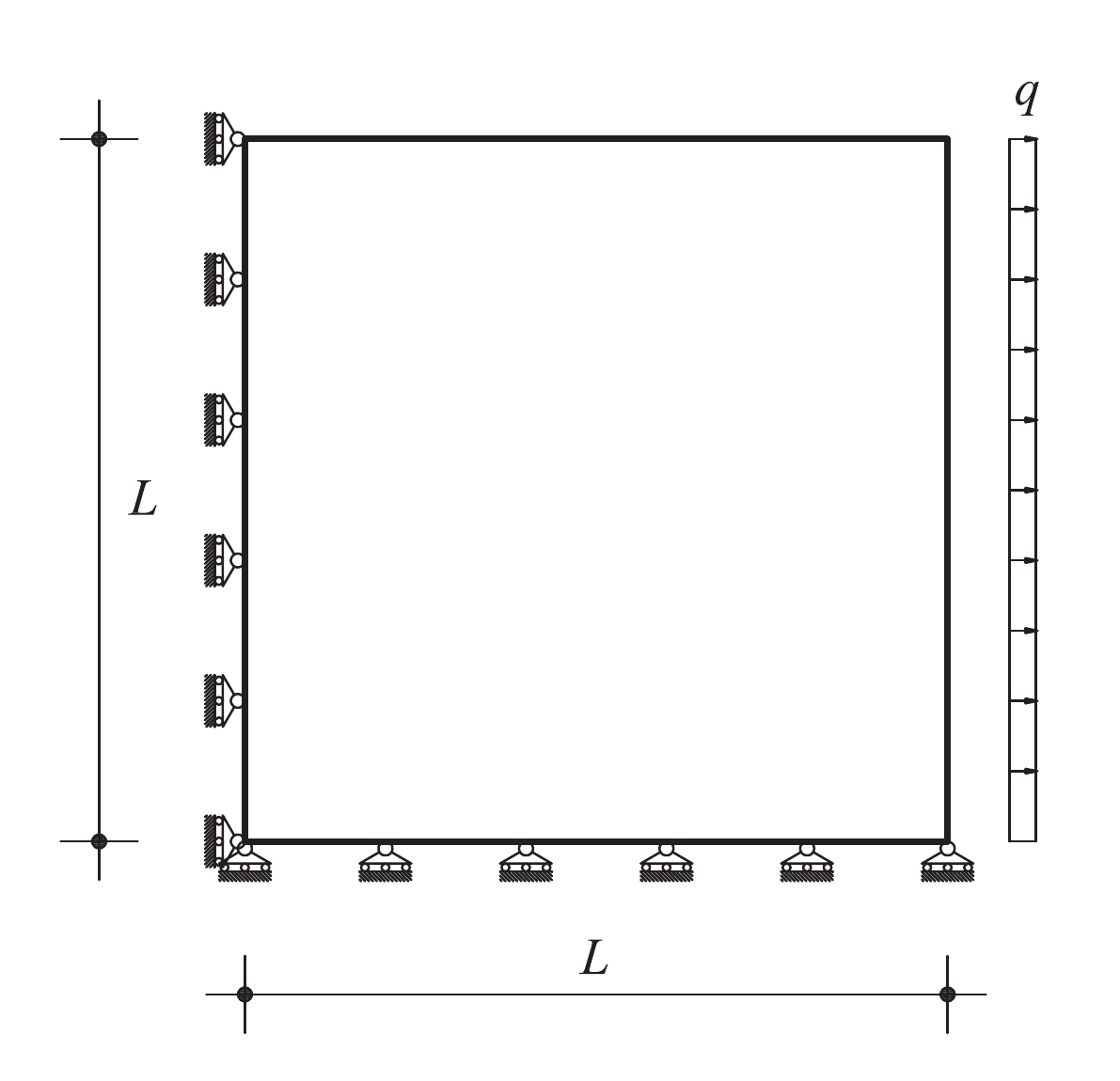}
\subcaption{}
\end{minipage}
\begin{minipage}[b]{.5\linewidth}
\centering
\includegraphics[bb = 150 18 180 320, scale = 0.5]
                {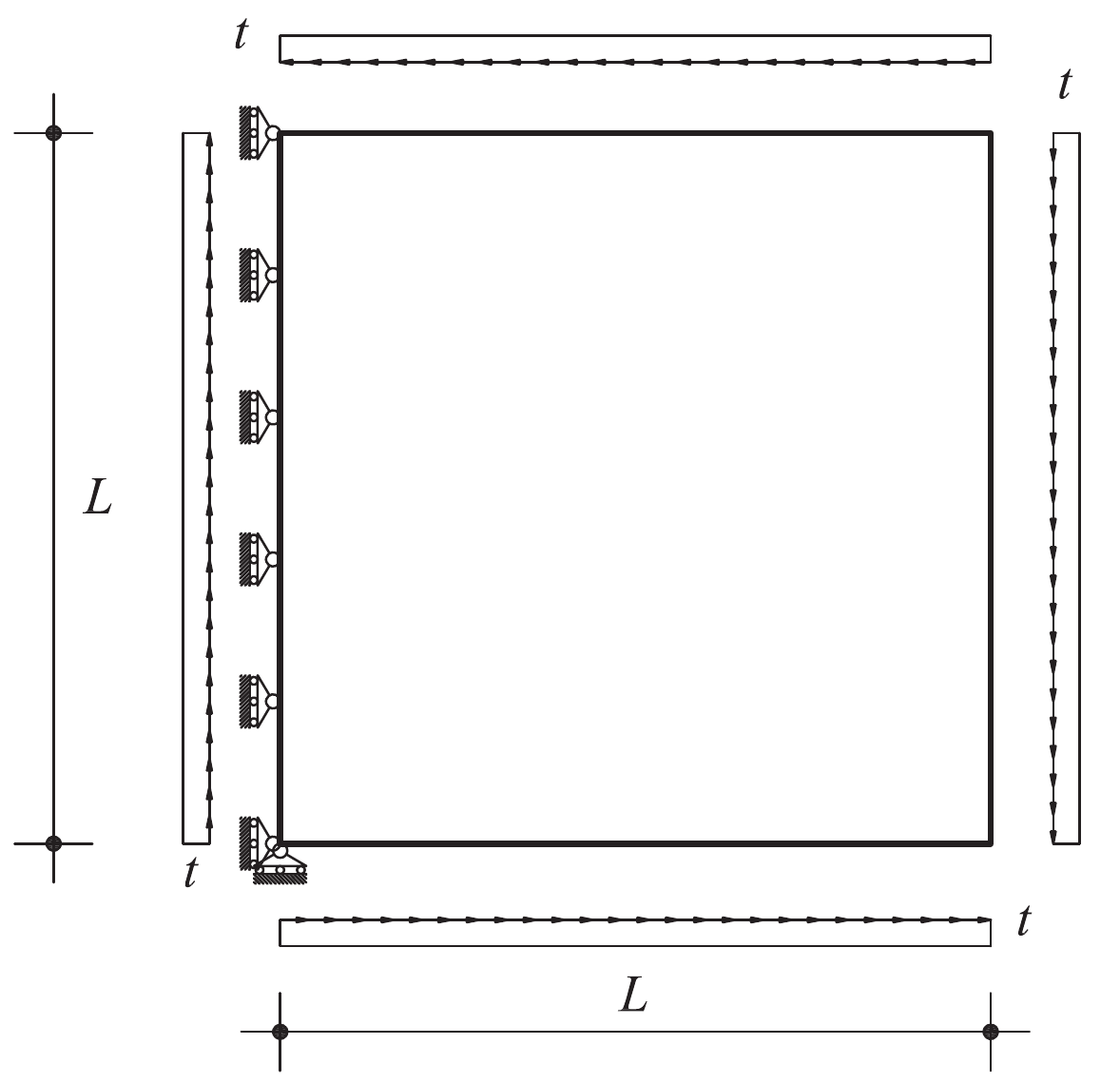}
\subcaption{}
\end{minipage}
\caption{Test 1: Linear patch tests. Geometry, loading and boundary conditions for tensile (a) and shear (b) constant stress states. }
\label{fig:Test_1_geom}
\end{figure}

We consider two boundary value problems on the unit square domain $\Omega = [0, 1]^2$ which induce constant stress states examined in \cite{Bathe1996}, under plain strain assumption. The material parameters are introduced setting the Young modulus $E = 7000$ and the Poisson ratio $\nu = 0.3$. Geometry, boundary conditions and loading cases for Test $1a$ [resp. Test $1b$] are reported in Fig. \ref{fig:Test_1_geom}(a) [resp. Fig. \ref{fig:Test_1_geom}(b)]. The former case represents a constant tensile stress state induced by the normal traction $q = 2000$ on the right edge of the unit square, the latter represents a constant shear stress state induced by the tangential traction $t = 400 $ assigned on the whole boundary. 

\begin{figure}
\begin{minipage}[b]{.5\linewidth}
\hspace{20mm}
\includegraphics[bb = 200 -20 600 480, scale=0.3]
                {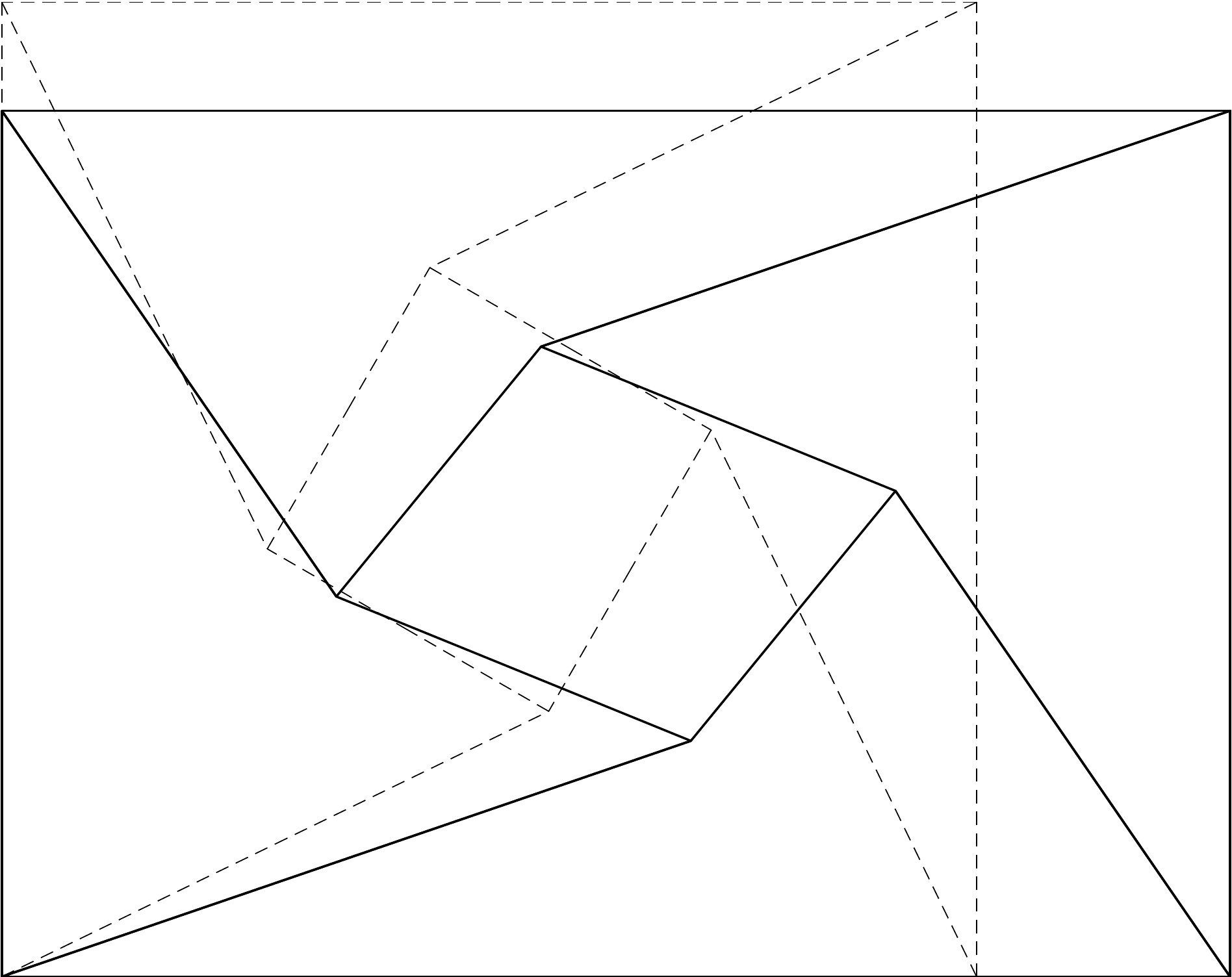}
\subcaption{}
\end{minipage}
\begin{minipage}[b]{.5\linewidth}
\hspace{40mm}
\includegraphics[bb = 220 18 400 620, scale=0.42]
                {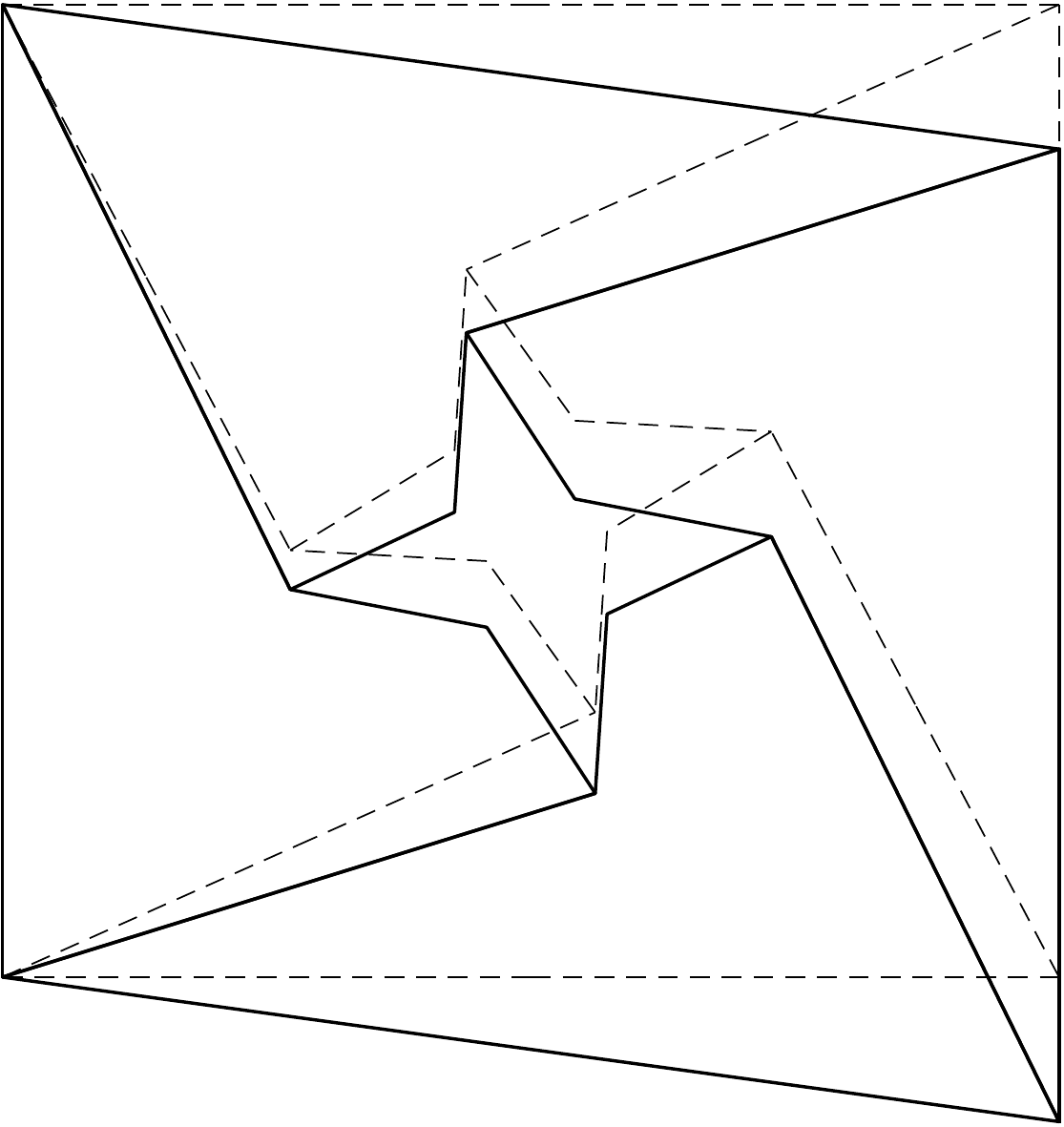}
\subcaption{}
\end{minipage}
\caption{
Test 1: Linear patch tests. Adopted mesh, reference (dash line) and deformed (continuous line) configurations for (a) tensile  and (b) shear constant stress states.}
\label{fig:Test_1_mesh}
\end{figure}

Linear VEMs are sufficient for the problem under investigation. Test $1a$ is analyzed with a very coarse mesh made of $4$ equal concave quadrilaterals and one square; the latter is centered in the center of $\Omega$ and has a side length $1/3$, and is rotated counter-clockwise of $\pi/3$ as can be seen in Fig. \ref{fig:Test_1_mesh}(a), where the meshes for the undeformed (dotted line) and deformed (continuous line) configurations are reported. Test $1b$ is analyzed with a very coarse mesh made of $4$ equal concave pentagons and one star-shaped octagon; the latter is centered in the center of $\Omega$ and has is derived by the square element of Test $1a$ by punching each midside point through the inside of a quarter of its side length, plus a counter-clockwise of $\pi/3$ as can be seen in Fig. \ref{fig:Test_1_mesh}(b), where the meshes for the undeformed (dotted line) and deformed (continuous line) configurations are reported.

\begin{figure}[htb]
\begin{minipage}[b]{.5\linewidth}
\hspace{10mm}
\includegraphics[bb = 80 5 400 400, scale=0.41]
                {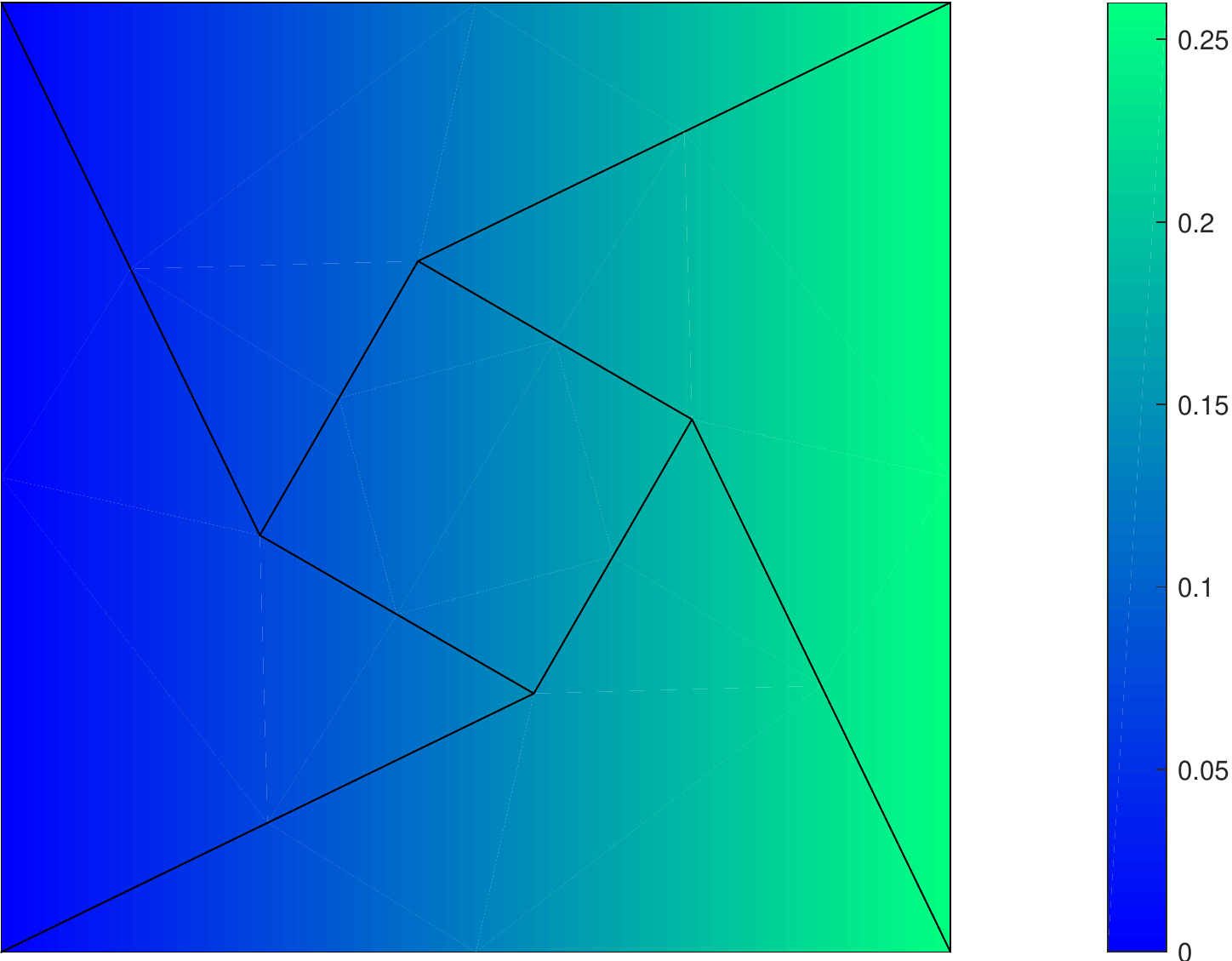}
\subcaption{}
\end{minipage}
\begin{minipage}[b]{.5\linewidth}
\hspace{13mm}
\includegraphics[bb = 80 5 400 400, scale=0.41]
                {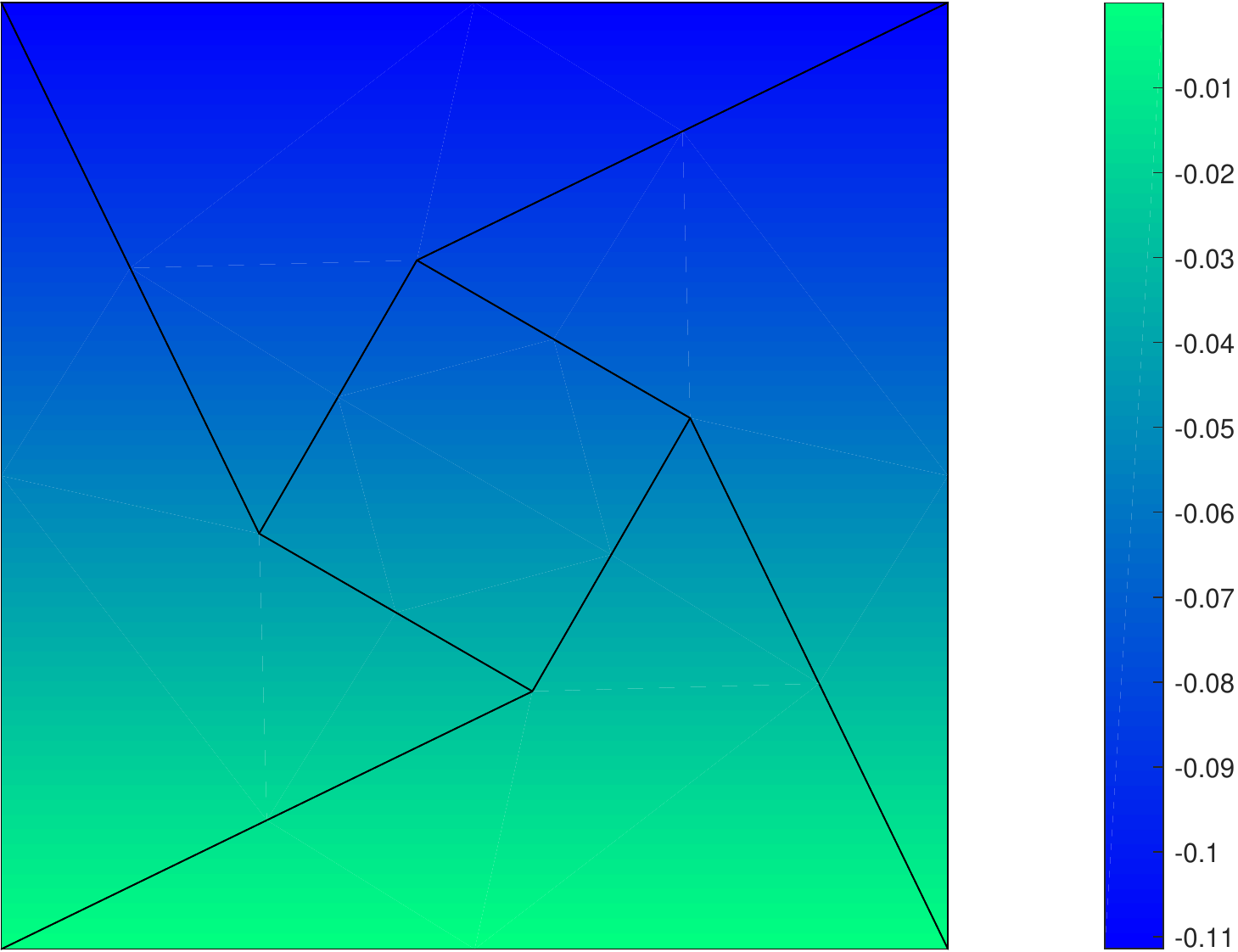}
\subcaption{}
\end{minipage}
\begin{minipage}[b]{.5\linewidth}
\hspace{10mm}
\includegraphics[bb = 80 5 400 400, scale=0.41]
                {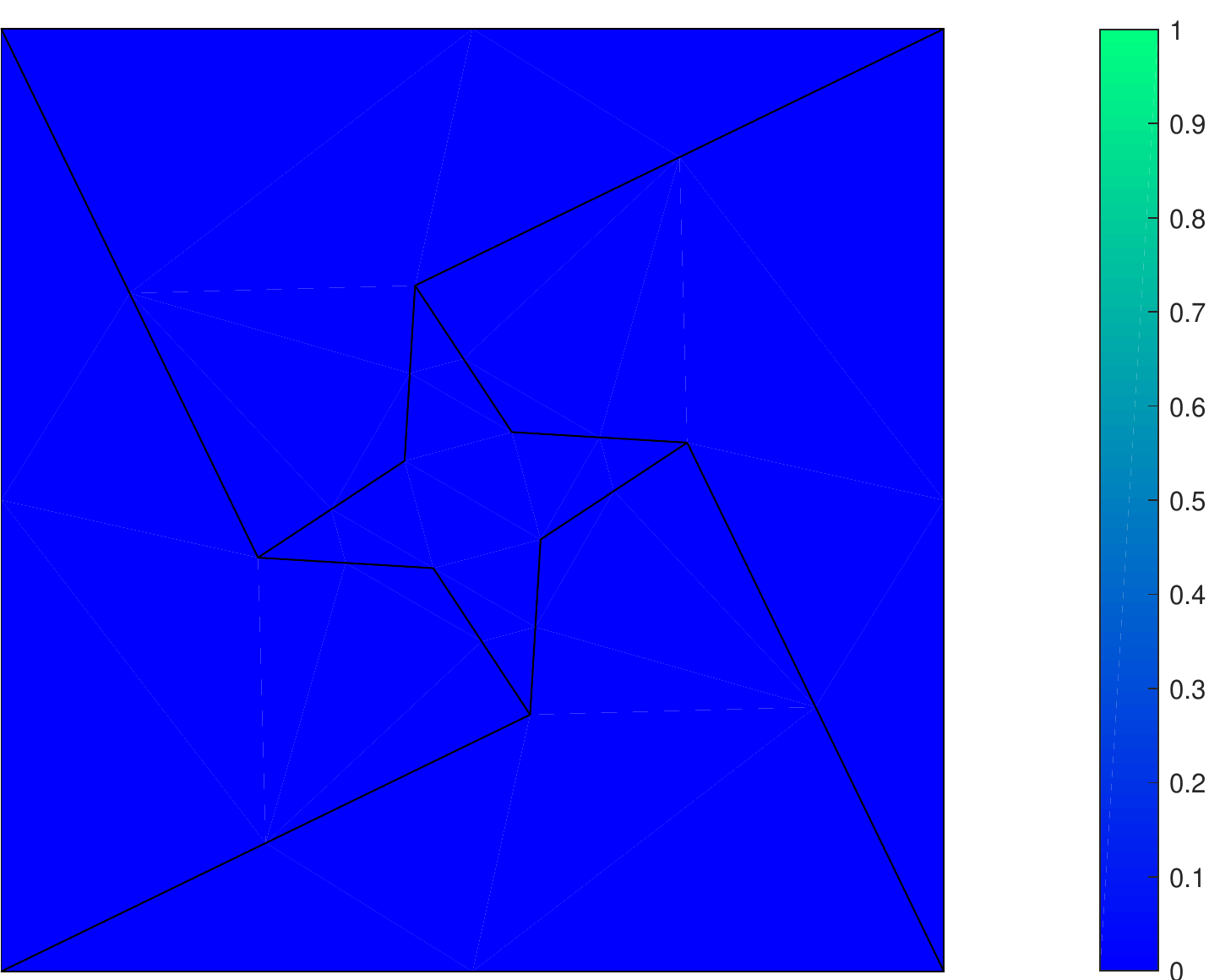}
\subcaption{}
\end{minipage}
\begin{minipage}[b]{.5\linewidth}
\hspace{13mm}
\includegraphics[bb = 80 5 400 400, scale=0.41]
                {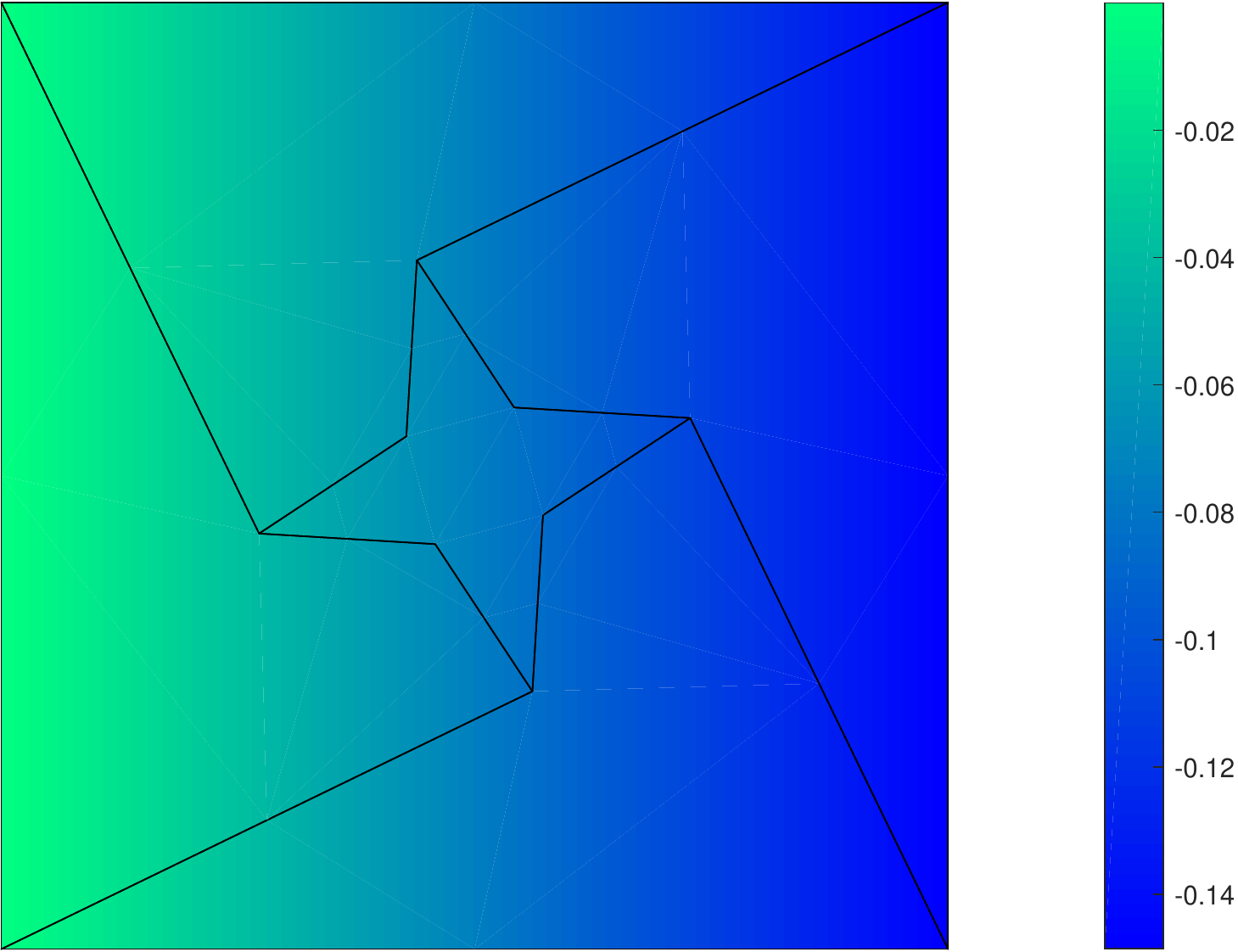}
\subcaption{}
\end{minipage}
\caption{
Test 1: Linear patch tests. Color plots of displacement field components $\left\{ u,v \right\}$ on undeformed meshes. Tensile stress state: $\left\{ \right.$(a),(b)$\left. \right \}$. Shear stress state: $\left\{ \right.$(c),(d)$\left. \right \}$.}
\label{fig:Test_1_disp}
\end{figure}

Qualitative color plots\footnote{For brevity, contour plots are interpolated through vertex values using command {\it fill} \cite{matlab:2015}.} for the displacement components stemming from the two test solutions are shown in Fig. \ref{fig:Test_1_disp}: as expected a linear pattern can be appreciated in both cases, being horizontal displacement for Test $1b$ obviously zero. The order of magnitude of the maximum euclidean norm of the difference between the analytical and the numerical stresses over a $100 \times 100$ uniform grid, and over the whole domain, is $10^{-16}$ in both cases, ascertaining that the method is capable of exactly representing a simple stress state.
 
\subsection{Test $2a$ and $2b$}
\label{ss:test_2a2b}
We consider two boundary value problems on the unit square domain $\Omega = [0, 1]^2$, with known analytical solution, discussed in \cite{BeiraoLovaMora}. In this case, the material parameters are assigned in terms of Lam\'e constants $\lambda = 1$, $\mu = 1$, and plain strain regime is still invoked. The tests are defined by choosing a required solution and deriving the corresponding load $\bb$, as synthetically indicated in the following:
\begin{itemize}
\item[$\bullet$]
Test $2a$
\begin{eqnarray}
\label{eq:test_2a_sol}
\left\{
\begin{array}{l}
u = x^3 - 3 x y^2 \\
v = y^3 - 3 x^2 y \\
\bb = \bzero
\end{array}
\right .
\end{eqnarray}
\item[$\bullet$]
Test $2b$
\begin{eqnarray}
\label{eq:test_2a_sol}
\left\{
\begin{array}{l}
u = v = \sin(\pi x)  \sin(\pi y) \\
b_1 = b_2 = -\pi^2 \left[ -(3 \mu + \lambda ) \sin(\pi x) \sin ( \pi y) + ( \mu + \lambda ) \cos ( \pi x) \cos ( \pi y ) \right]
\end{array}
\right .
\end{eqnarray}
\end{itemize}
As it can be immediately observed, Test $2a$ is a problem with Dirichlet non-homogeneous boundary conditions, zero loading and a polynomial solution; whereas Test $2b$ has homogeneous Dirichlet boundary conditions, trigonometric distributed loads with a trigonometric solution.

For the purpose of comparing the accuracy level and the convergence rate of the proposed VEM to standard FEM, meshes with only quadrilaterals will be considered for Test $2a$. Linear and quadratic VEMs are compared to Lagrangian linear and quadratic quadrilateral finite elements, indicated in the following as $Q4$ and $Q9$, respectively \cite{Zienckiewicz_Taylor2000}. Test $2b$ is instead analyzed with general polygonal meshes only testing VEM solutions.

\begin{figure}
\hspace{45mm}
\begin{minipage}[b]{.5\linewidth}
\centering
\includegraphics[bb = 80 5 400 400, scale=0.42]
                {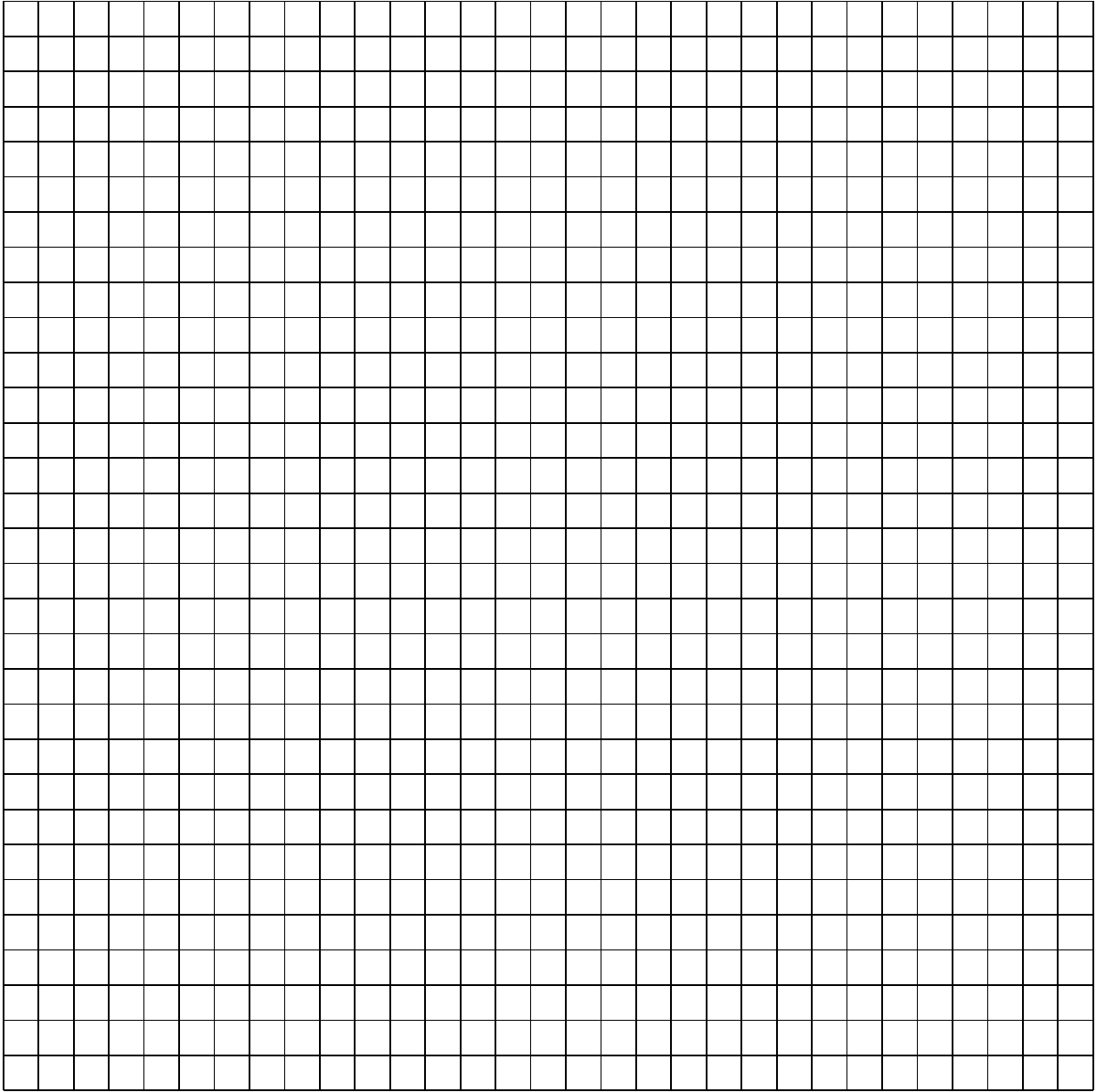}
\subcaption{}
\end{minipage}%
\\

\hspace{45mm}
\begin{minipage}[b]{.5\linewidth}
\centering
\includegraphics[bb = 80 5 400 400, scale=0.42]
                {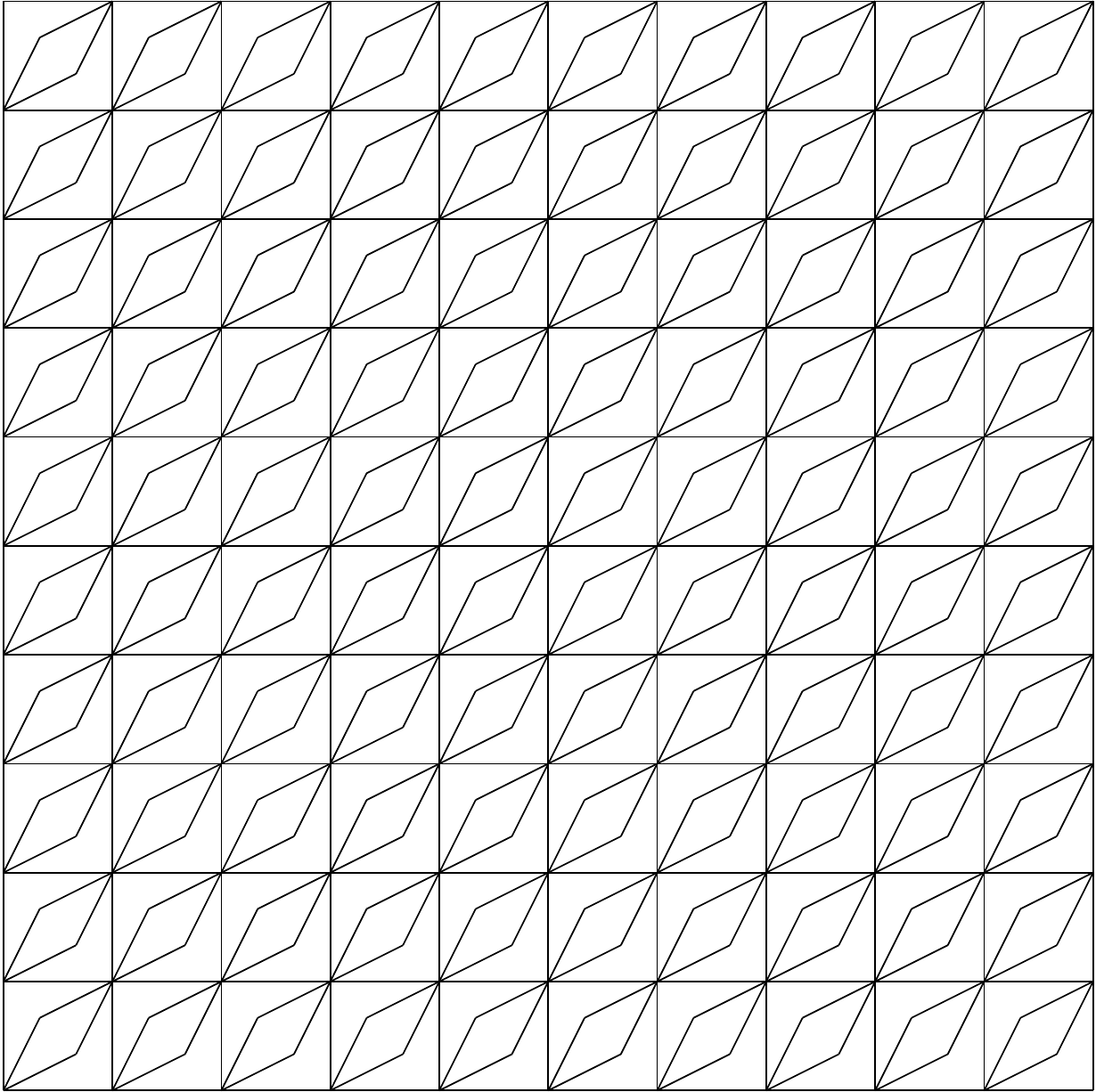}
\subcaption{}
\end{minipage}%
\\

\hspace{45mm}
\begin{minipage}[b]{.5\linewidth}
\centering
\includegraphics[bb = 80 5 400 400, scale=0.42]
                {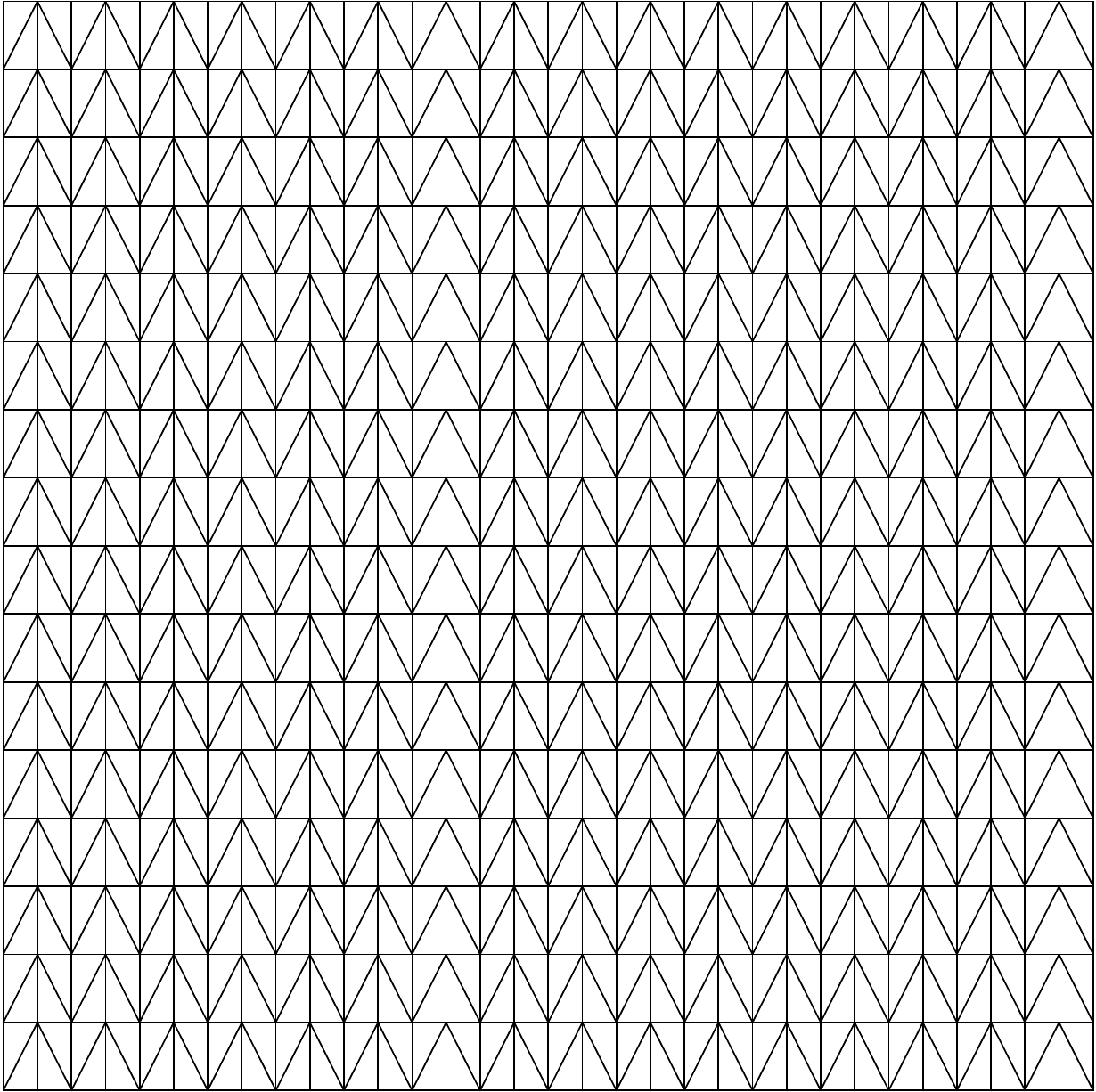}
\subcaption{}
\end{minipage}
\caption{
Test 2a: Sample meshes. Uniform mesh of squares (a); Distorted concave rhombic quadrilaterals (b); Rectangle trapezoids almost collapsing into triangles (c).}
\label{fig:Test_2_a_mesh}
\end{figure}

\begin{figure}
\vspace{20mm}
\begin{minipage}[b]{.5\linewidth}
\centering
\includegraphics[bb = 80 5 400 400, scale=0.41]
                {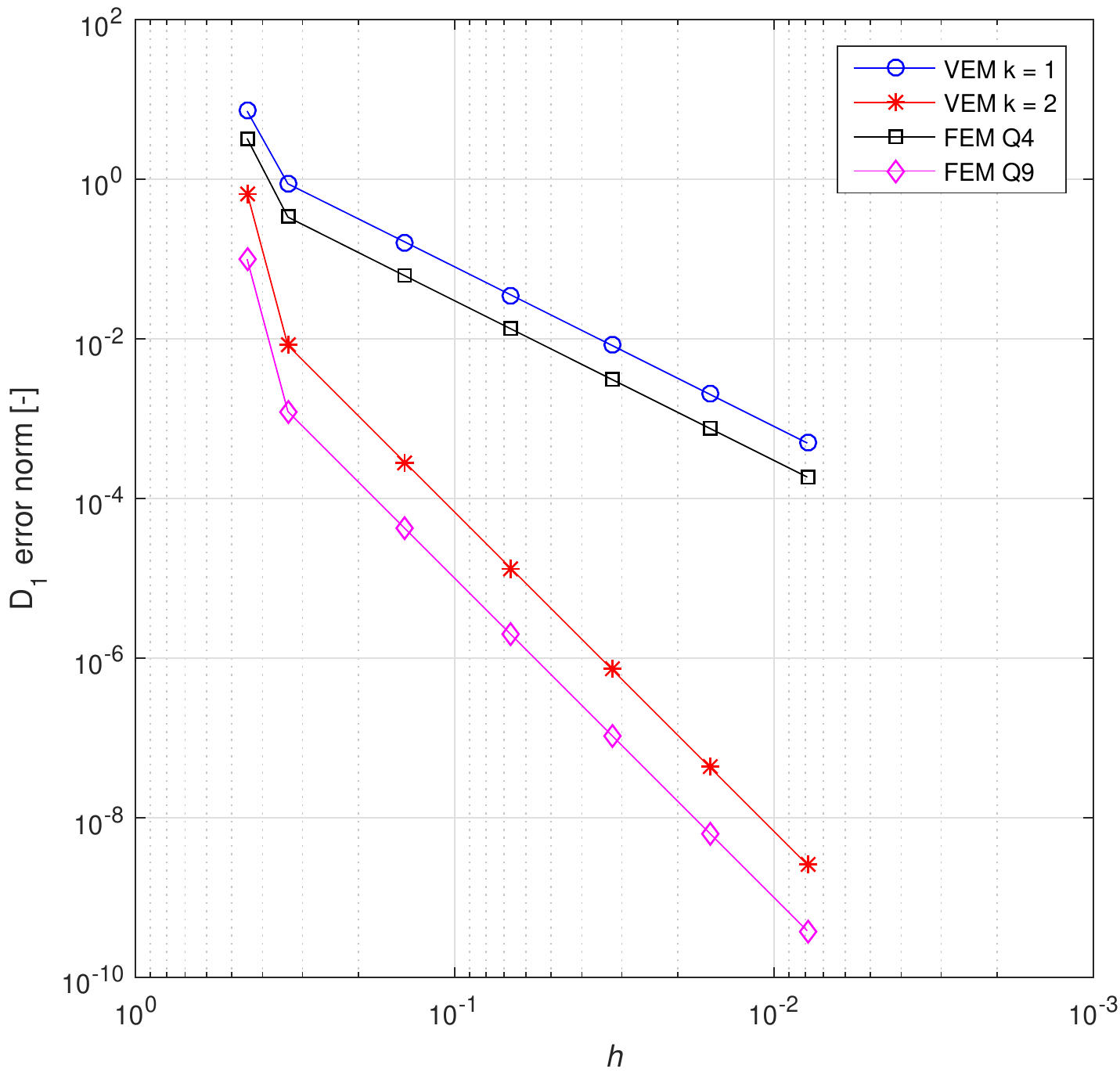}
\subcaption{}
\end{minipage}%
\begin{minipage}[b]{.5\linewidth}
\centering
\includegraphics[bb = 80 5 400 400, scale=0.41]
                {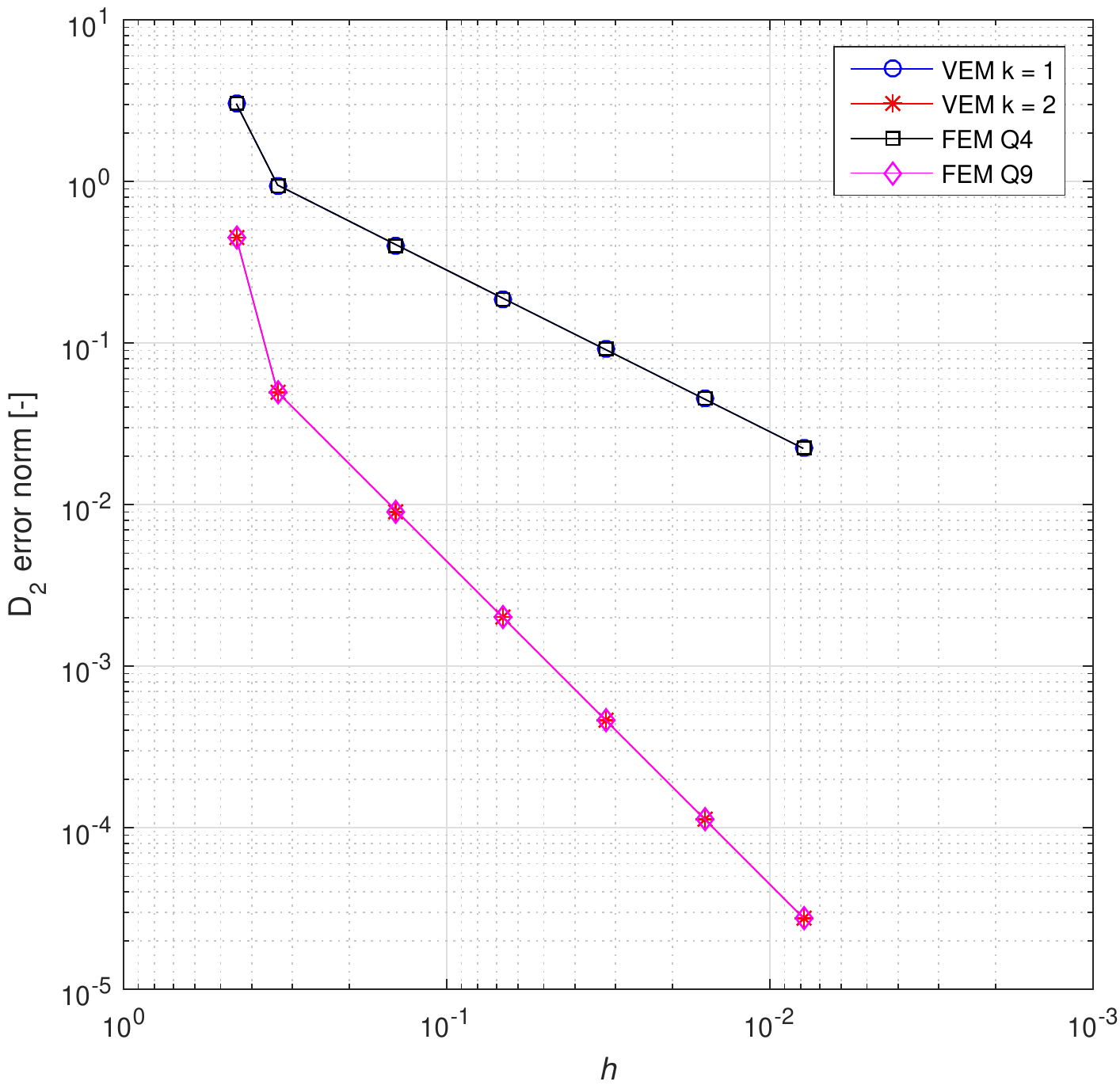}
\subcaption{}
\end{minipage}
\begin{minipage}[b]{.5\linewidth}
\centering
\includegraphics[bb = 80 5 400 400, scale=0.41]
                {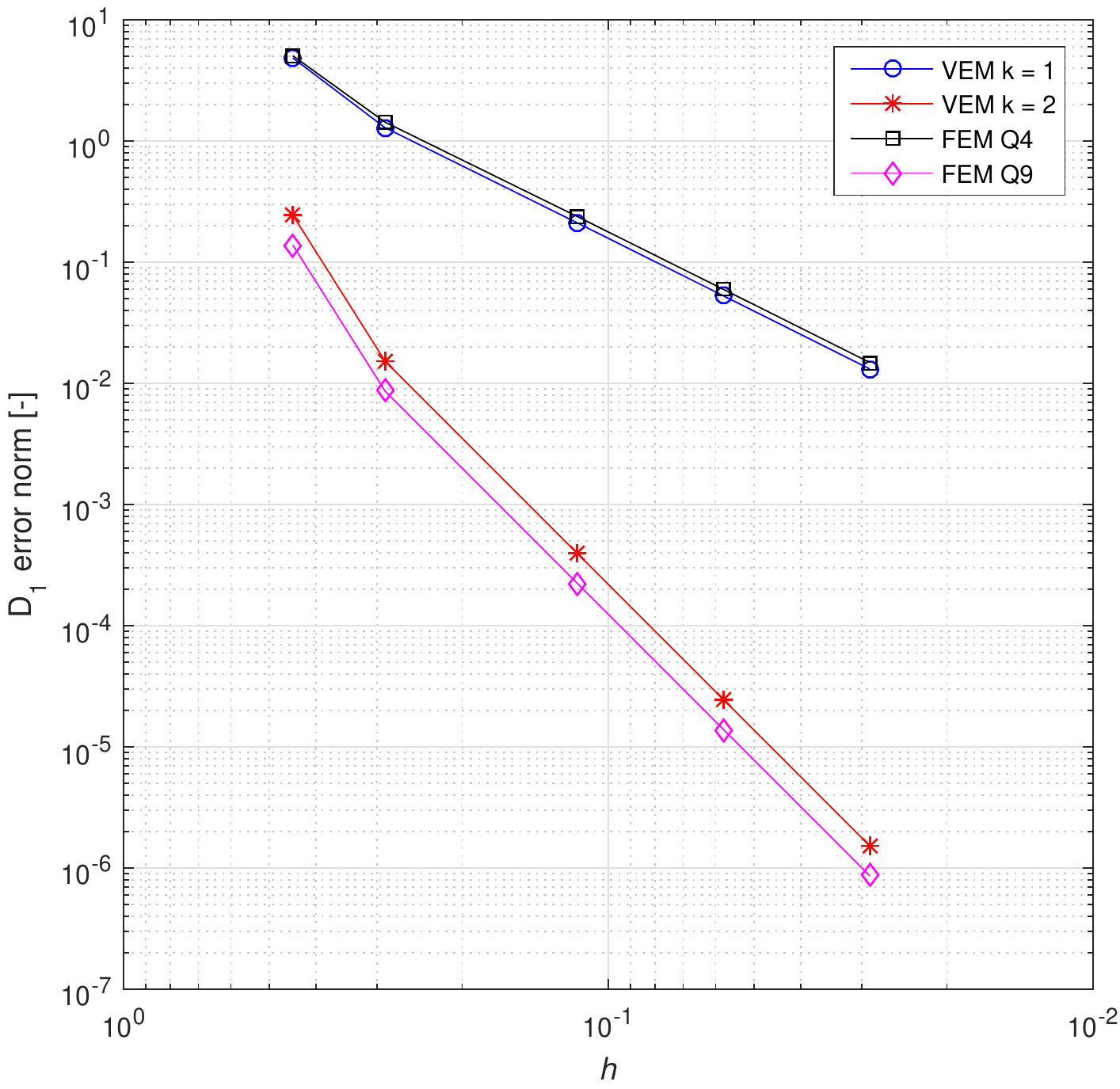}
\subcaption{}
\end{minipage}%
\begin{minipage}[b]{.5\linewidth}
\centering
\includegraphics[bb = 80 5 400 400, scale=0.41]
                {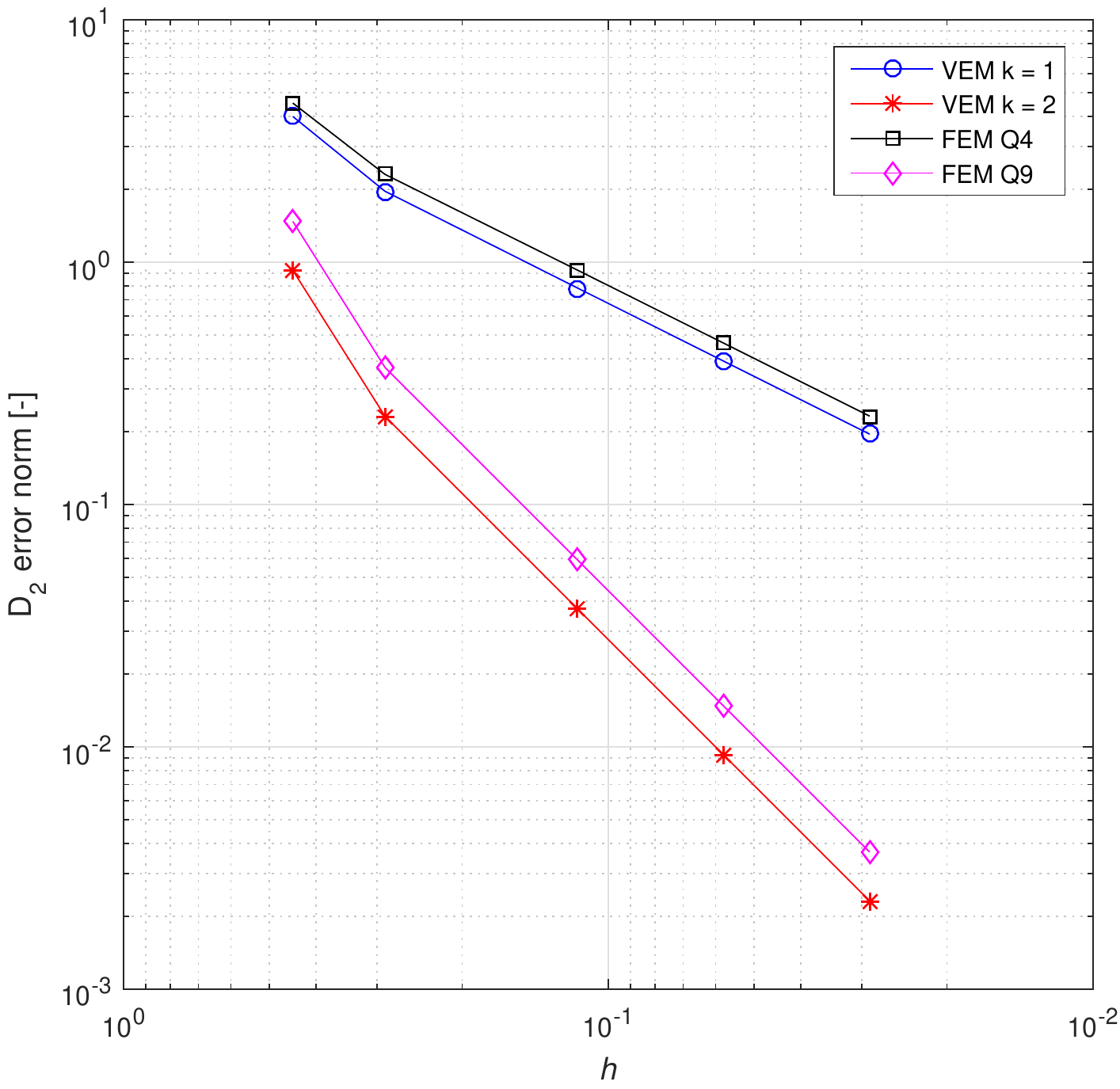}
\subcaption{}
\end{minipage}
\begin{minipage}[b]{.5\linewidth}
\centering
\includegraphics[bb = 80 5 400 400, scale=0.41]
                {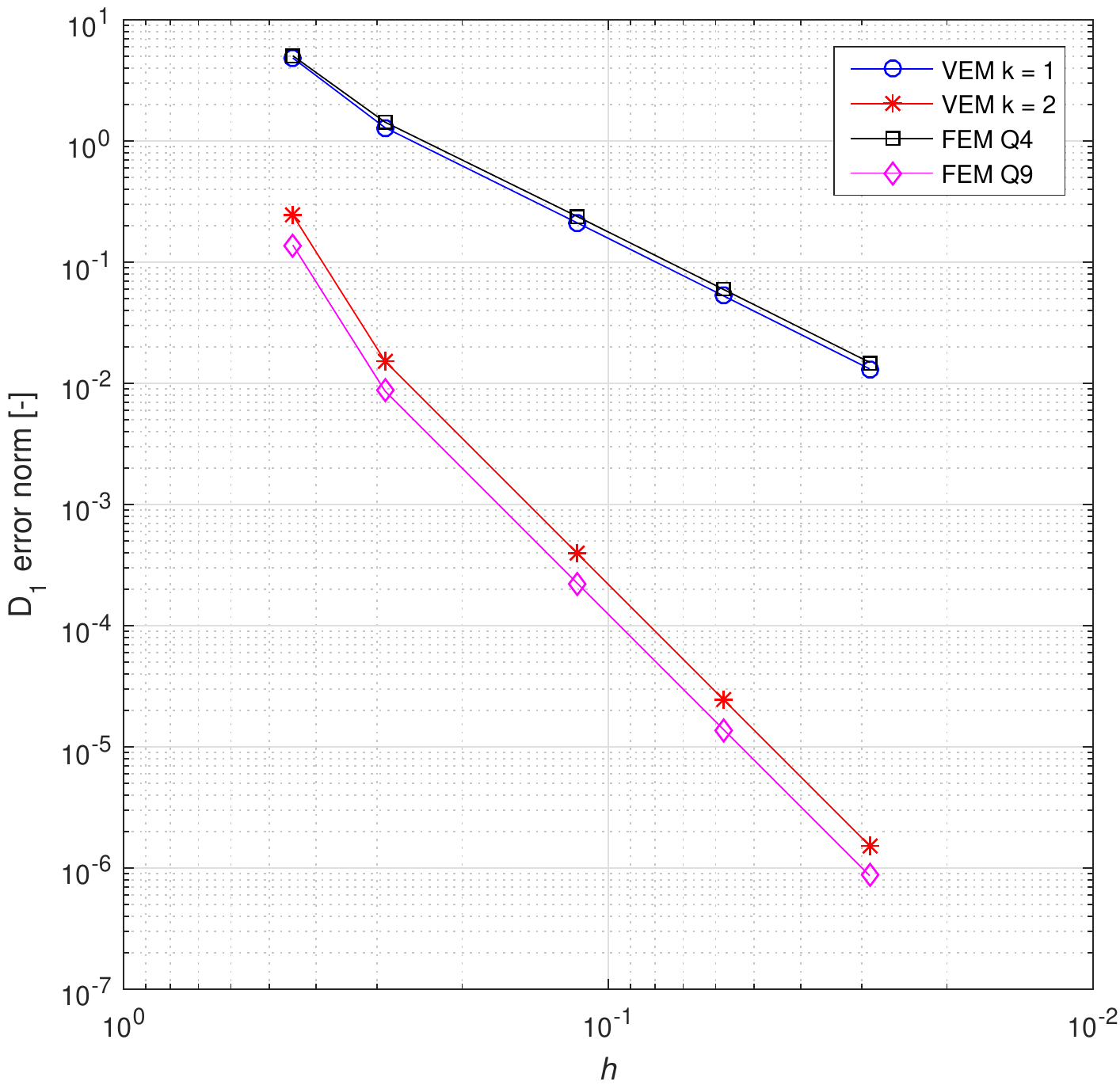}
\subcaption{}
\end{minipage}
\begin{minipage}[b]{.5\linewidth}
\centering
\includegraphics[bb = 80 5 400 400, scale=0.41]
                {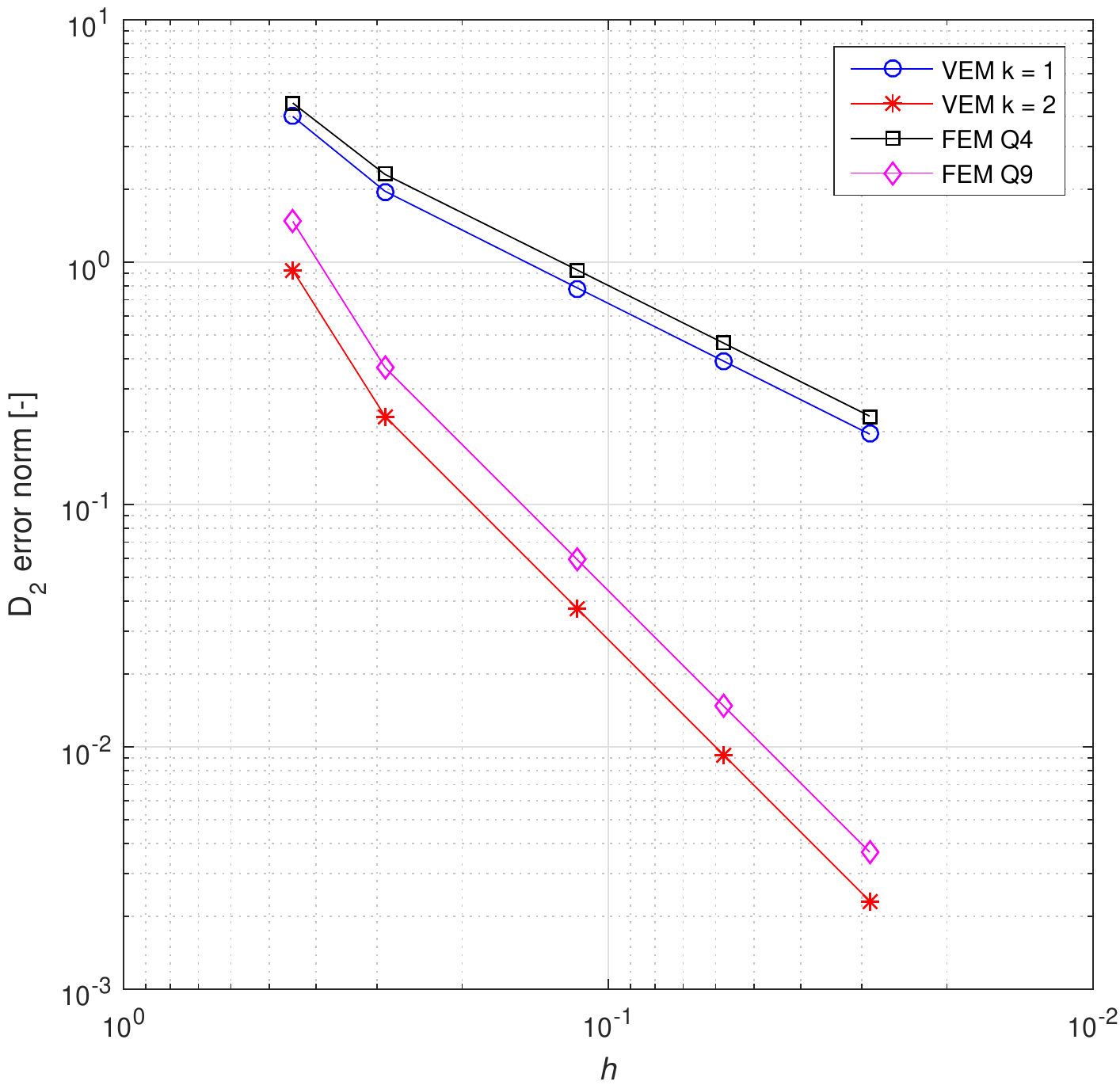}
\subcaption{}
\end{minipage}
\caption{
Test 2a: convergence plot. Error norms $D_1$-$D_2$ {\it vs} element size $h$. FEM $Q4$ and $Q9$; VEM $k=1$ and $k=2$, $m=4$. Uniform meshes of squares (a)-(b); Distorted concave rhombic quadrilaterals (c)-(d); Non-uniform meshes of collapsing quadrilaterals (e)-(f).}
\label{fig:Test_2_a_error_D1_D2}
\end{figure}

In the following, such accuracy and convergence rate assessment is carried out using the following error norms:
\begin{itemize}
\item[$\bullet$] Energy-type error norm:
\begin{equation}
\label{eq:energy_err_norm}
D_1 = ||| \v_{\textrm{ex}} - \v_h |||_{1,2} := \sum_{E \in \Th} \int_{E} \left\| \beps_{\textrm{ex}} - \beps(\v_h) \right \|^2
\end{equation}
where $\beps_{\textrm{ex}}$ is the exact strain, and $\beps(\v_h)$ is the numerical strain according to either a FEM or a VEM solution. In particular, for the VEM solution, such a numerical strain is computed by means of \eqref{uVuP} and \eqref{eq:def-proj} while for the FEM solution it is obtained by standard derivation of the displacements.
\item[$\bullet$] Discrete $H^1$-type error norm:
\begin{equation}
\label{eq:H1_discrete_err_norm}
D_2 = |||\v|||_{1,2}:=\left( \sum_{f \in {\cal E}_h} \int_{f} h_f \left\| \frac{\partial \v}{\partial \t_f} \right\|_{0,e}^2 \right)^{1/2}
\end{equation}
where $h_f$ denotes the length of any edge $f$ of the reticulation ${\cal T}_h$, and $\t_f$ denotes the unit tangent vector to the edge $f$ chosen once and for all. 
\end{itemize}
Sample meshes for Test $2a$ are represented in Fig. \ref{fig:Test_2_a_mesh}, comprising evenly distributed squares (a), rhombic and concave quadrilaterals (b), and rectangular trapezoids almost collapsing into triangles (c). Correspondingly, Figure \ref{fig:Test_2_a_error_D1_D2}(a)-(f) show the $h-$convergence plots in terms of error $D_1$ and $D_2$ for the $Q4$, $Q9$ FEM, and $k=1$, $k=2$ VEM solutions. It is observed that, the tested methods present the expected convergence rates, namely linear and quadratic, respectively. The comparison in terms of error $D_1$ indicates that FEM methods slightly outperform the VEM of same order; while the opposite patterns is encountered in terms of error $D_2$. 
This seems to indicate that the VEM solution on the mesh skeleton is slightly better than the FEM solution in the presence of distorted meshes (see $D_2$ errors). On the other hand the VEM strain, since it is computer through the ``filter'' of the projection operator instead of being directly derived from the local displacements, suffers from a slight loss of accuracy (see $D_1$ errors). Possibly making use of more sofisticated post-processing procedures to compute the VEM strains could yield results that are more accurate than FEM. Nevertheless, it is important to remark that all the errors shown in Figure \ref{fig:Test_2_a_error_D1_D2} are comparable; this result fits exactly into the scope of the present test $2a$, that is to show than even for conforming quadrilaterals (that is, playing on a ground where FEM are very strong) the VEM is able to obtain competitive results. 

\begin{figure}
\vspace{20mm}
\begin{minipage}[b]{.5\linewidth}
\hspace{10mm}
\includegraphics[bb = 80 5 400 400, scale=0.42]
                {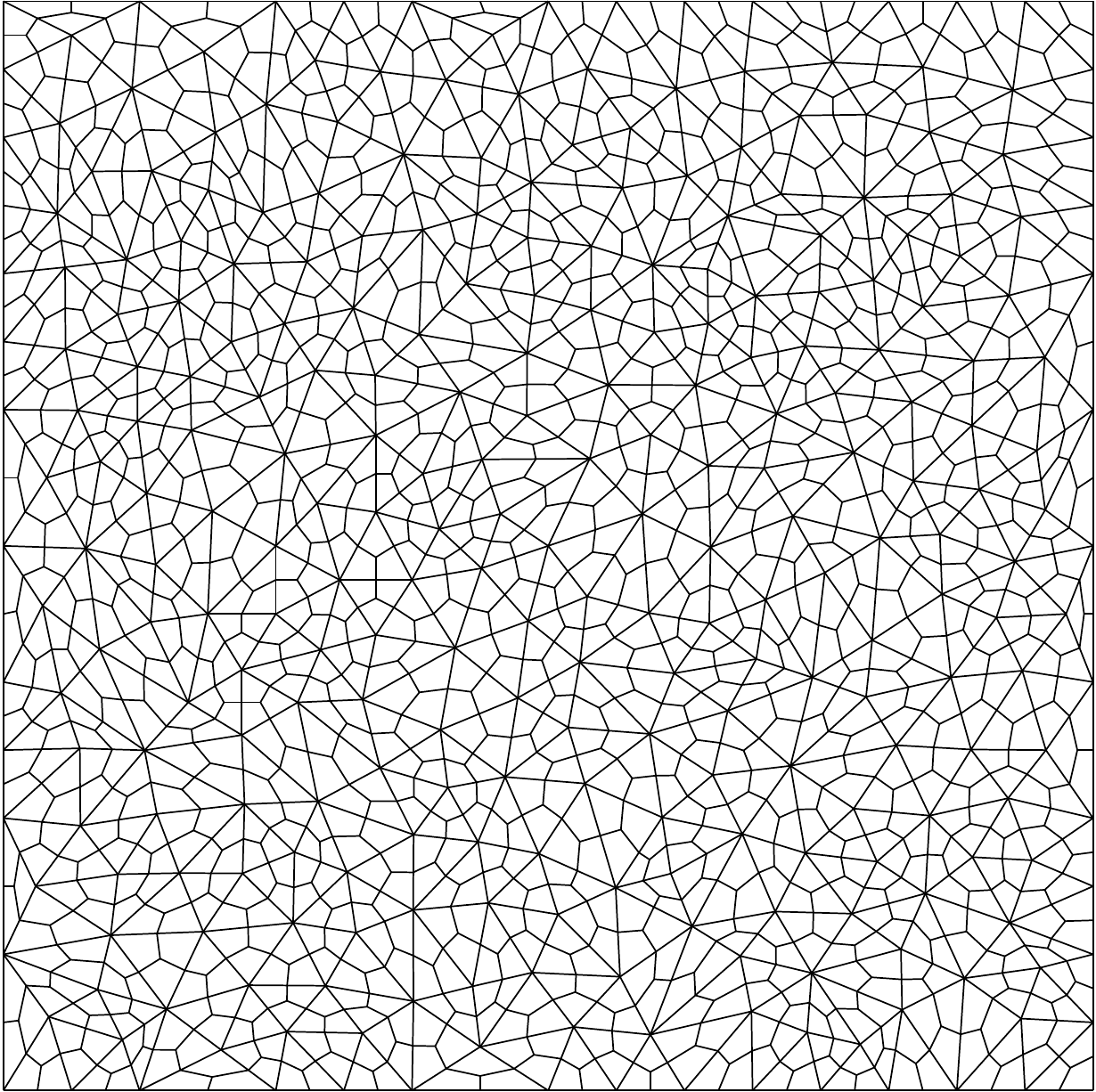}
\subcaption{}
\end{minipage}%
\begin{minipage}[b]{.5\linewidth}
\hspace{10mm}
\includegraphics[bb = 80 5 400 400, scale=0.42]
                {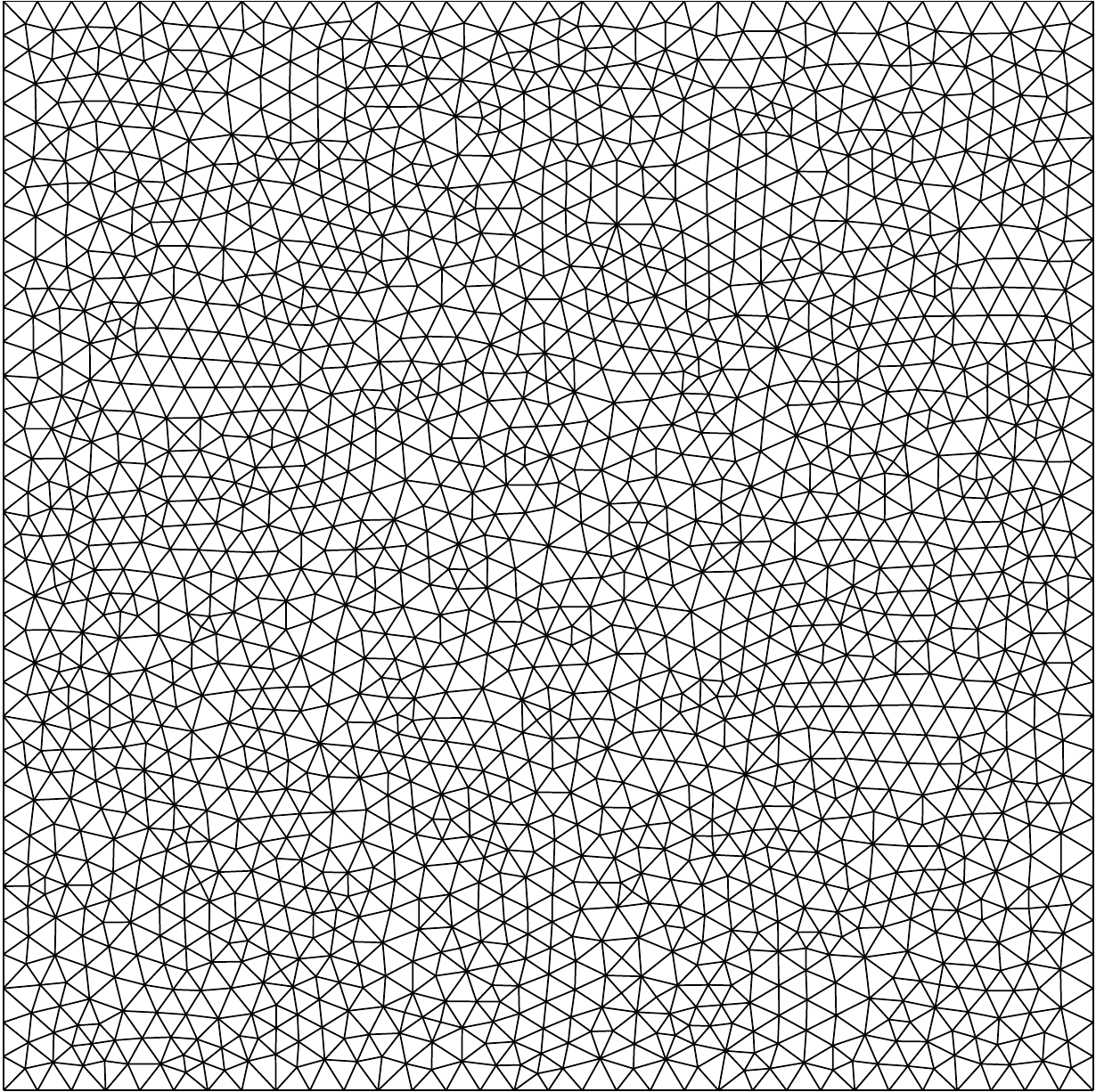}
\subcaption{}
\end{minipage}
\begin{minipage}[b]{.5\linewidth}
\hspace{10mm}
\includegraphics[bb = 80 5 400 400, scale=0.42]
                {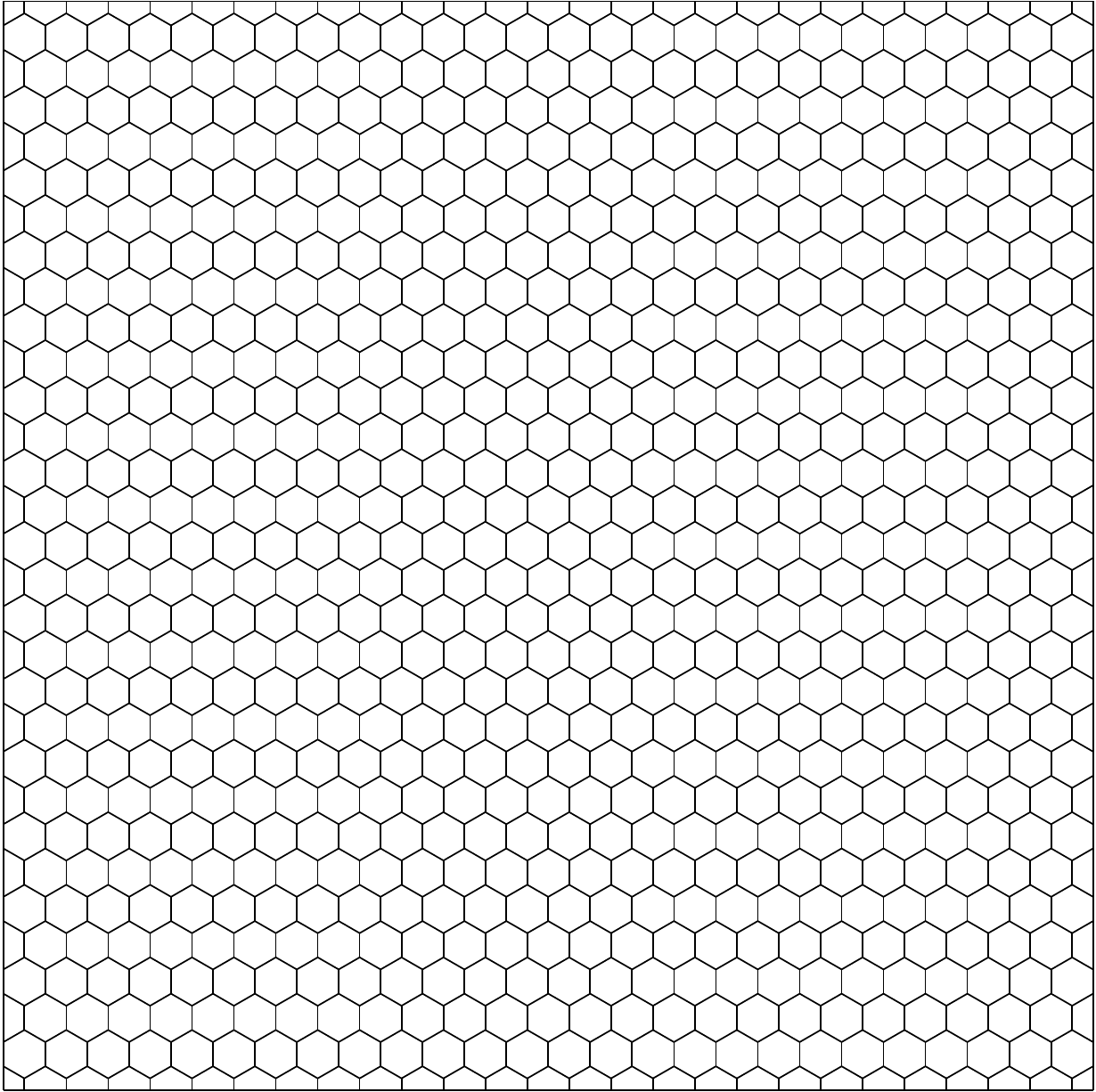}
\subcaption{}
\end{minipage}%
\begin{minipage}[b]{.5\linewidth}
\hspace{10mm}
\includegraphics[bb = 80 5 400 400, scale=0.42]
                {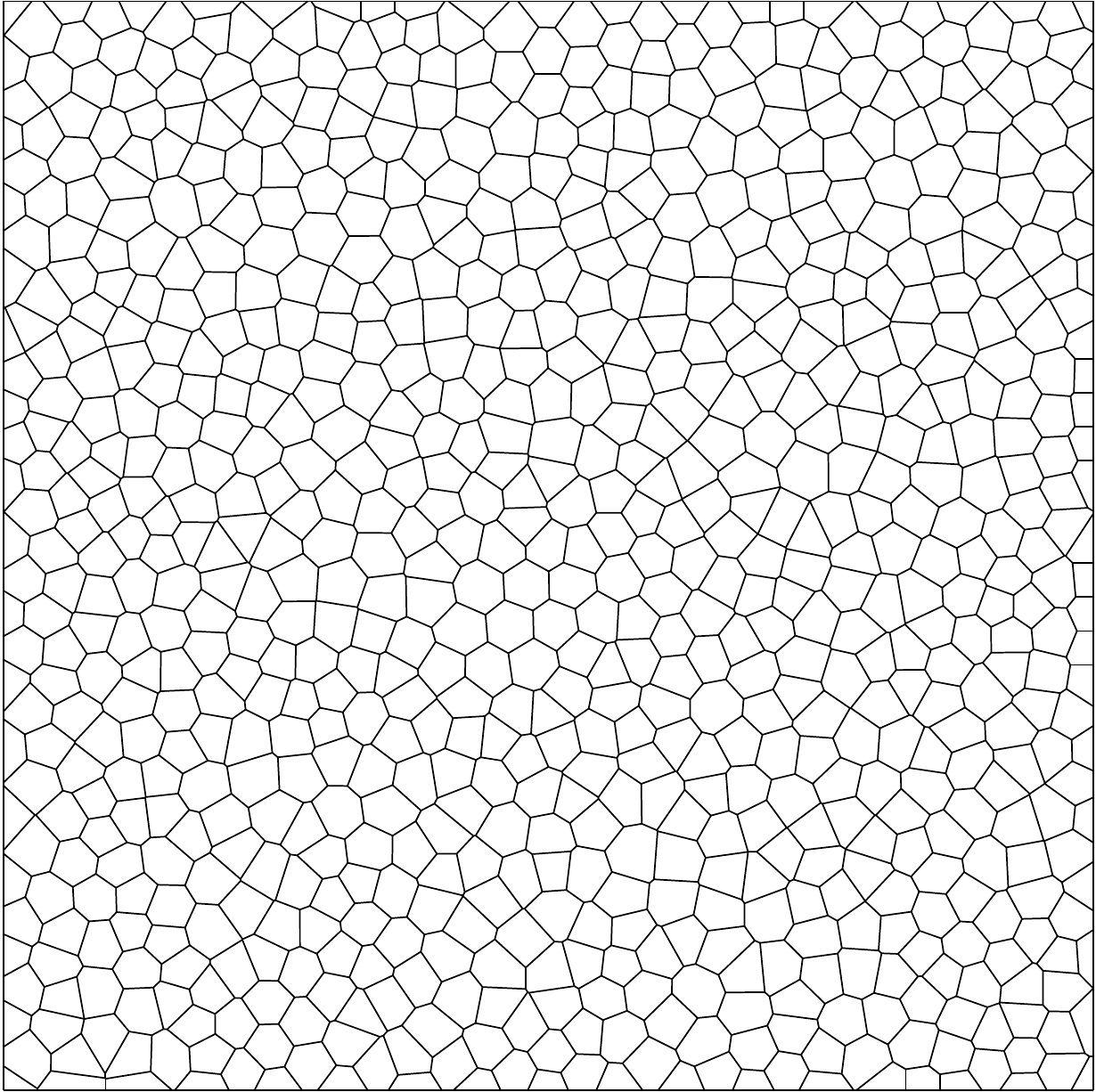}
\subcaption{}
\end{minipage}
\caption{
Test 2b: Sample meshes. Non-uniform mesh of convex quadrilaterals (a); Non-uniform mesh of triangles (b); Uniform mesh of hexagons (c); Voronoi tessellation (d).}
\label{fig:Test_2_b_mesh}
\end{figure}

\begin{figure}
\vspace{20mm}
\begin{minipage}[b]{.5\linewidth}
\centering
\includegraphics[bb = 80 5 400 400, scale=0.41]
                {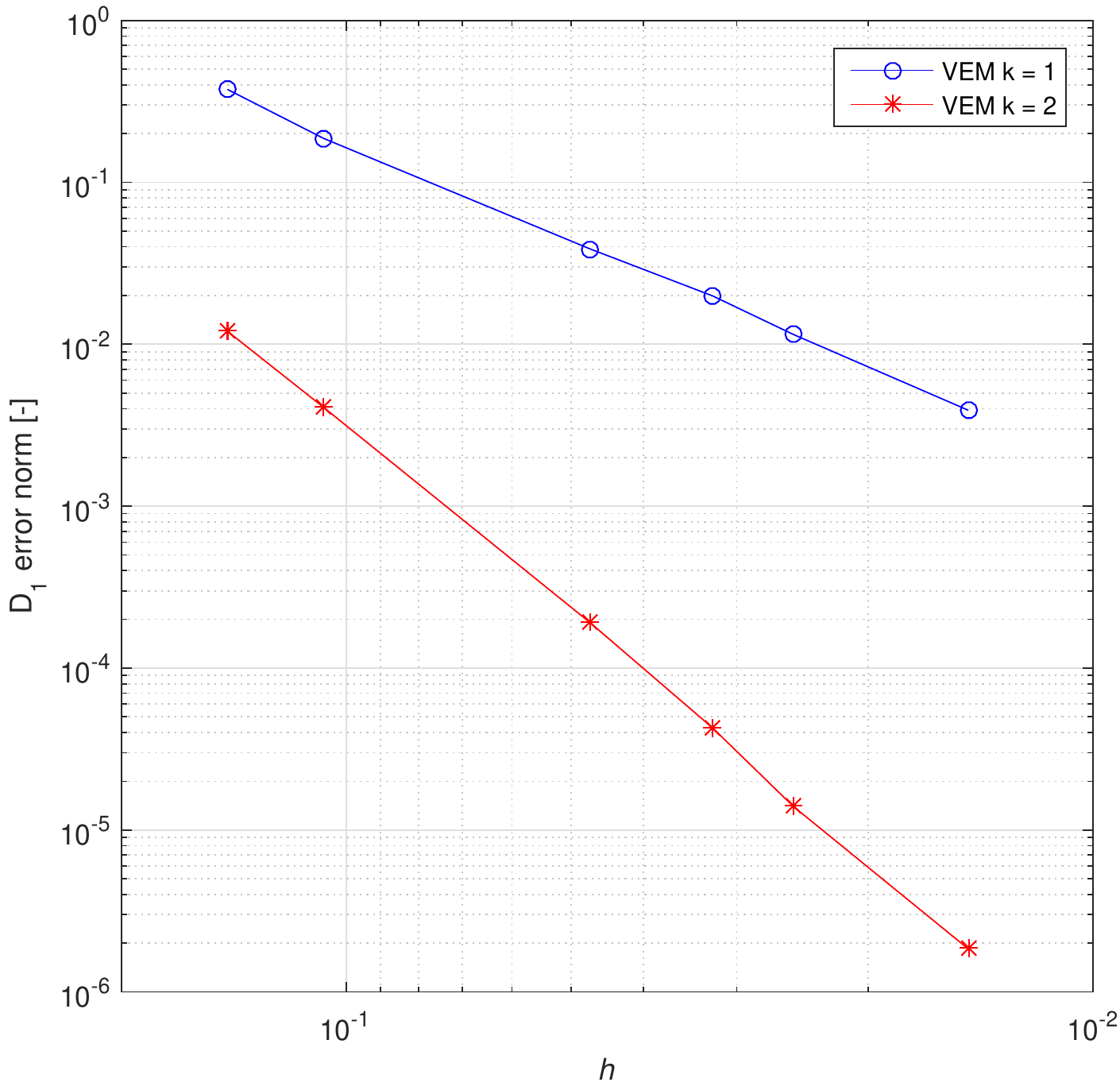}
\subcaption{}
\end{minipage}%
\begin{minipage}[b]{.5\linewidth}
\centering
\includegraphics[bb = 80 5 400 400, scale=0.41]
                {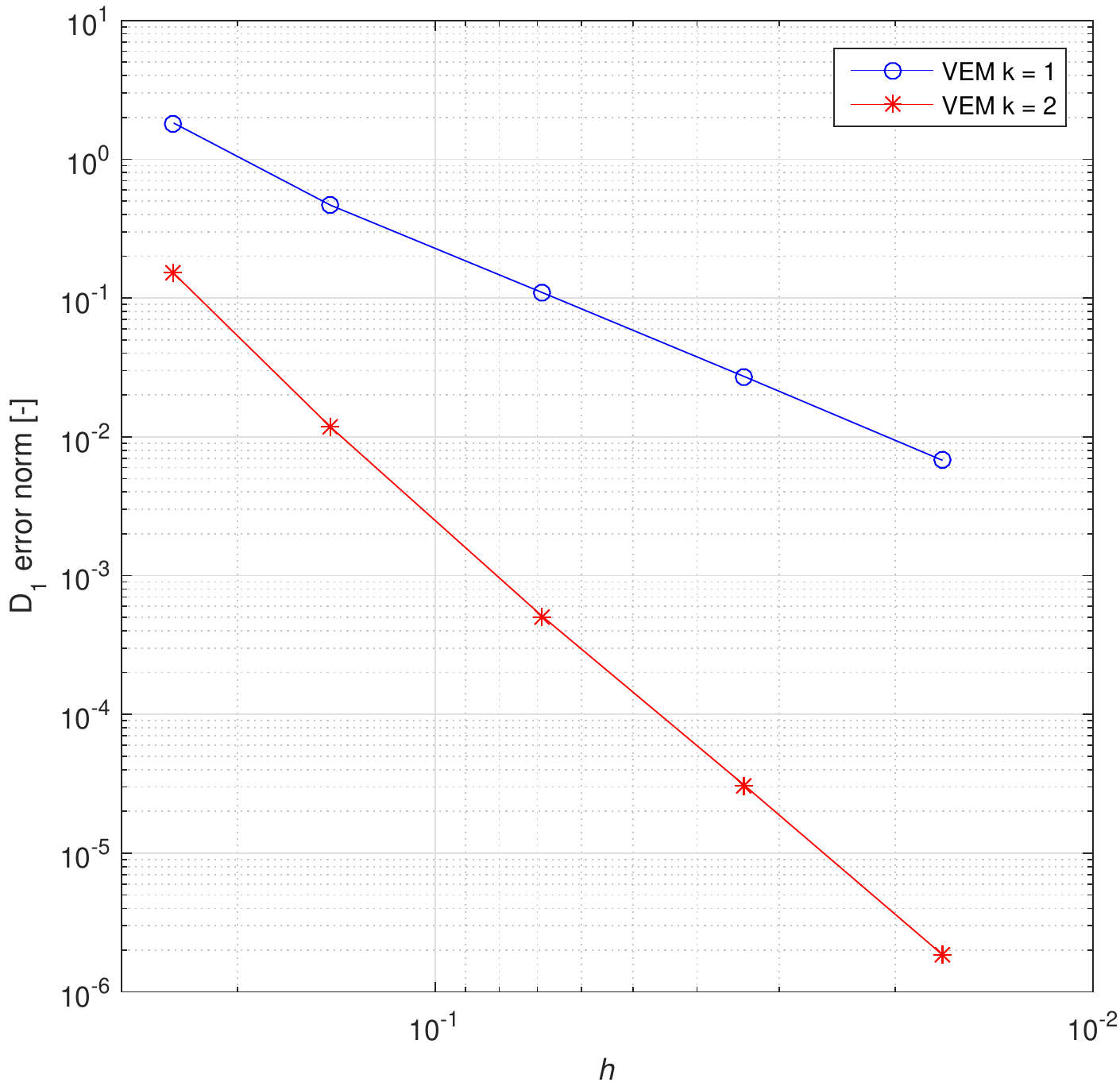}
\subcaption{}
\end{minipage}
\begin{minipage}[b]{.5\linewidth}
\centering
\includegraphics[bb = 80 5 400 400, scale=0.41]
                {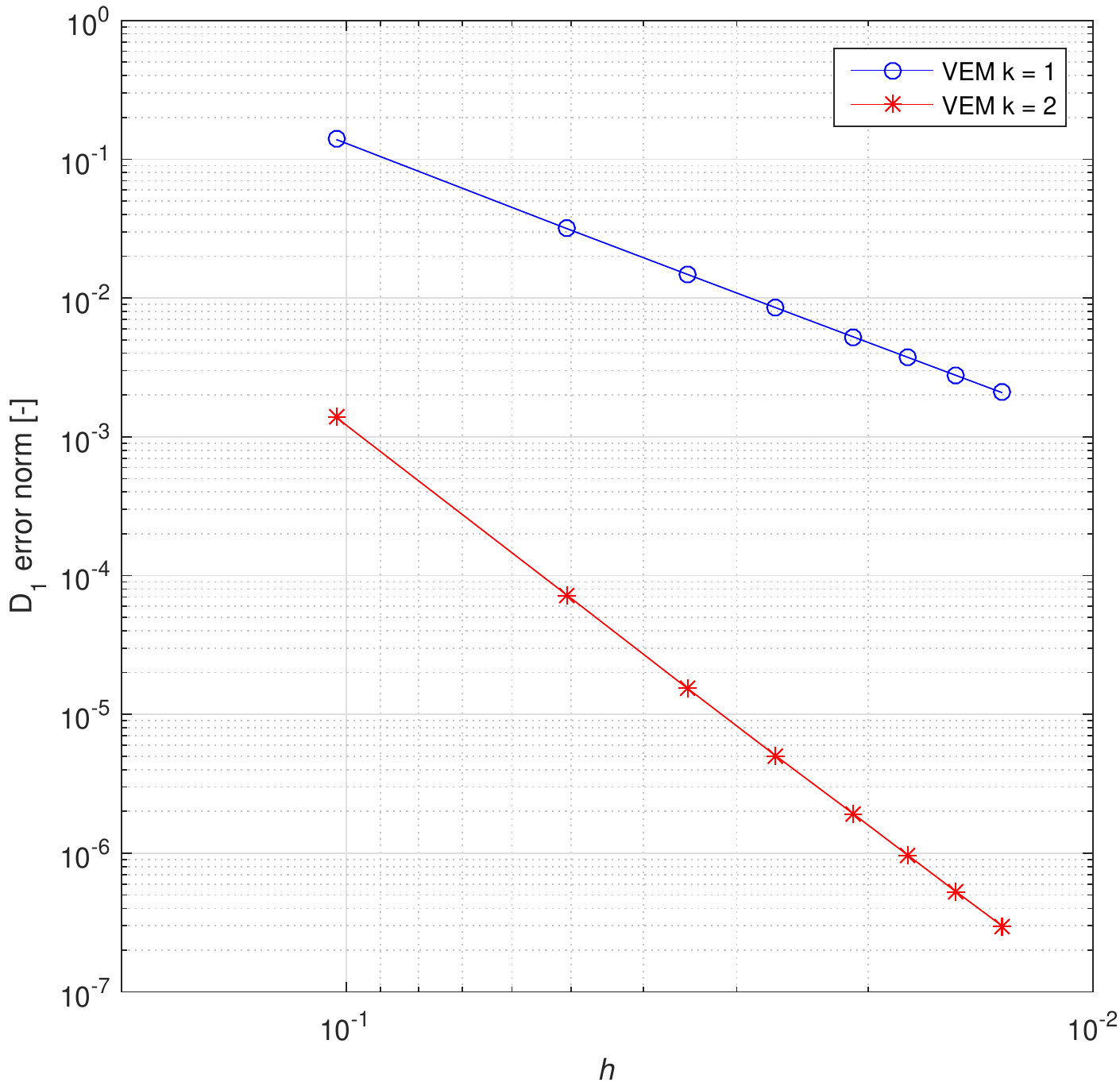}
\subcaption{}
\end{minipage}%
\begin{minipage}[b]{.5\linewidth}
\centering
\includegraphics[bb = 80 5 400 400, scale=0.41]
                {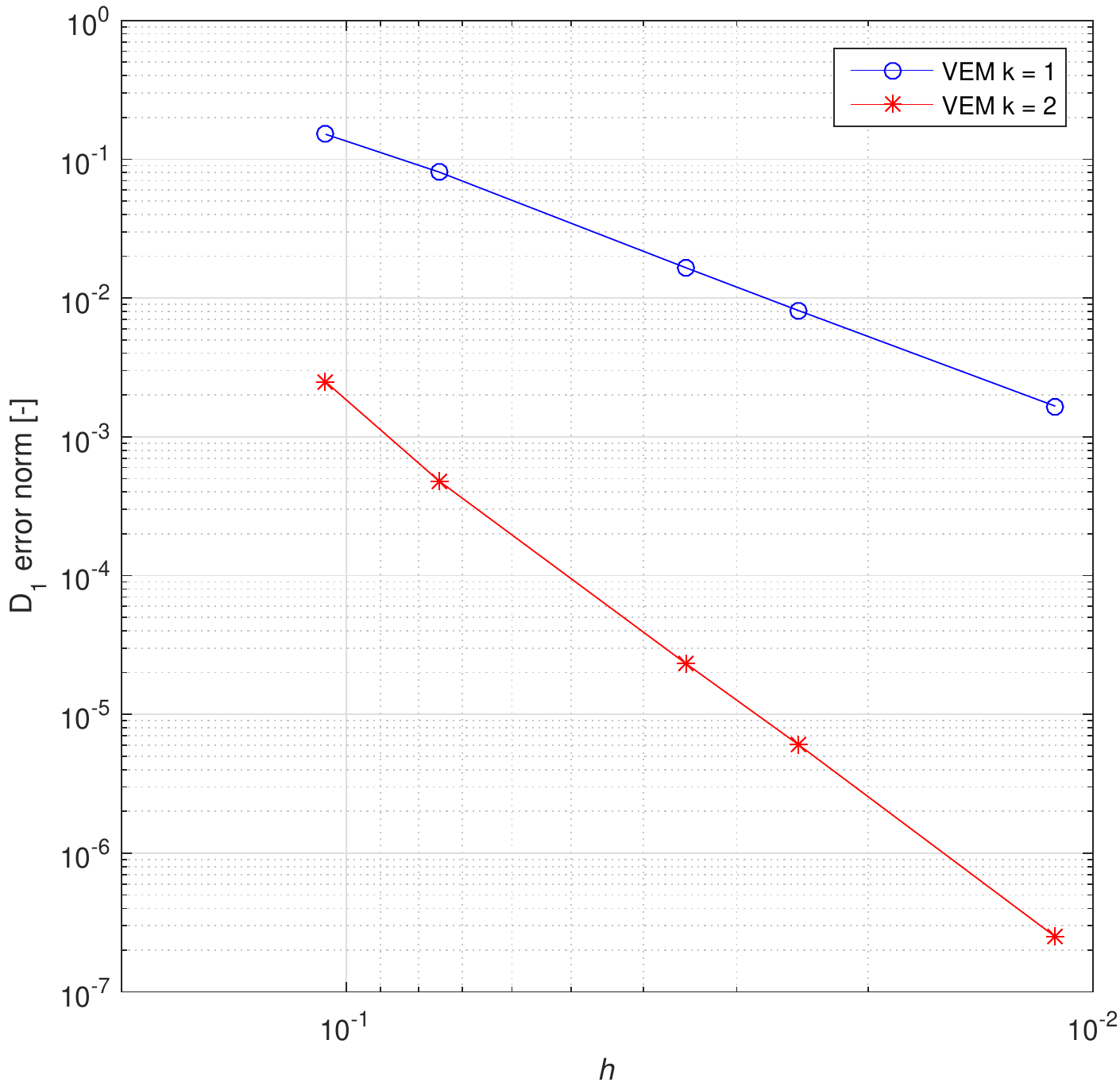}
\subcaption{}
\end{minipage}
\caption{
Test 2b: convergence plot. Error $D_1$ {\it vs} element size $h$. VEM $k=1$ and $k=2$. Non-uniform meshes of convex quadrilaterals (a); non-uniform meshes of triangles (b); uniform meshes of hexagons (c); Voronoi tessellations (d).}
\label{fig:Test_2_b_error_D1}
\end{figure}

\begin{figure}
\vspace{20mm}
\begin{minipage}[b]{.5\linewidth}
\centering
\includegraphics[bb = 80 5 400 400, scale=0.41]
                {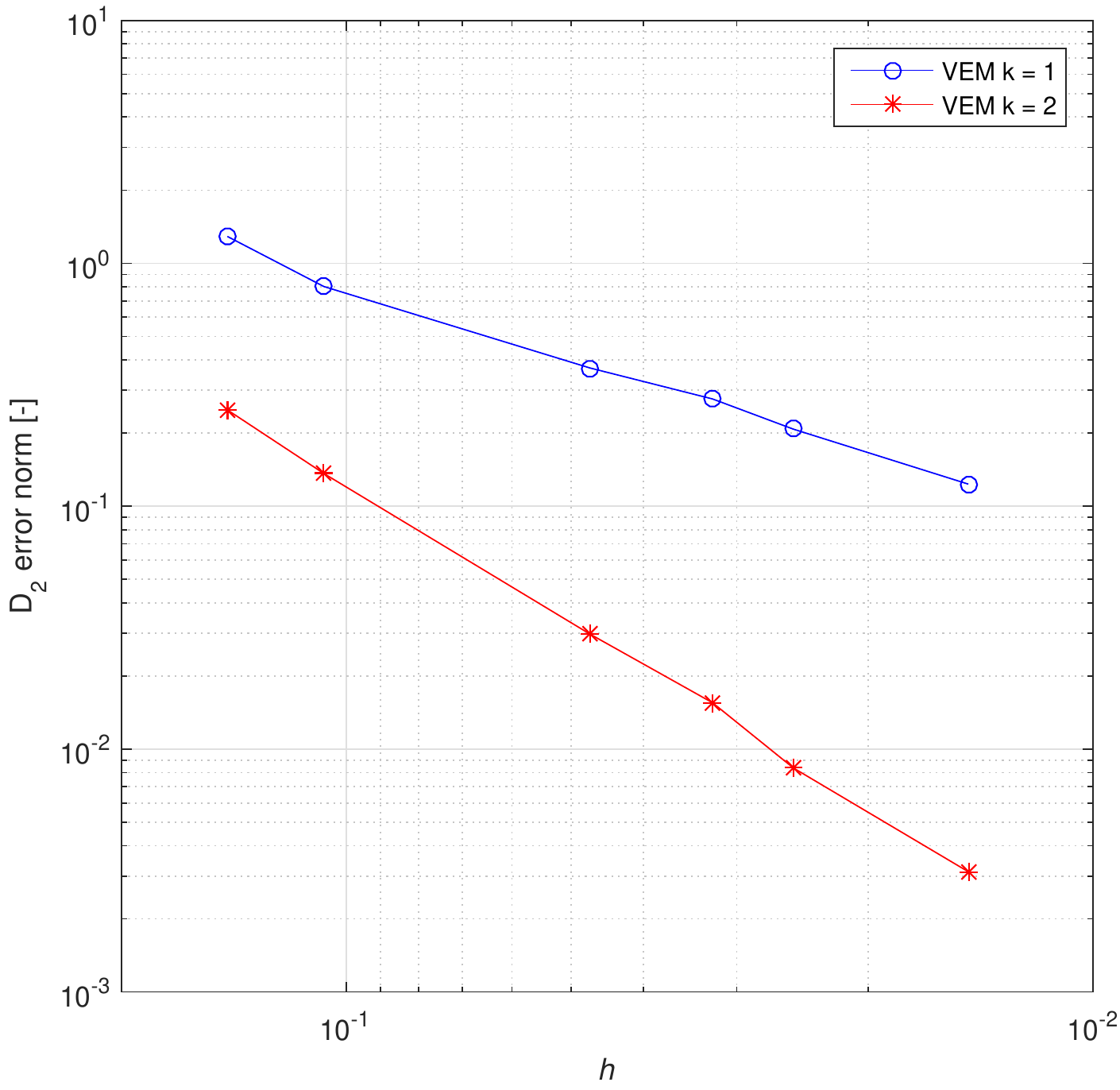}
\subcaption{}
\end{minipage}%
\begin{minipage}[b]{.5\linewidth}
\centering
\includegraphics[bb = 80 5 400 400, scale=0.41]
                {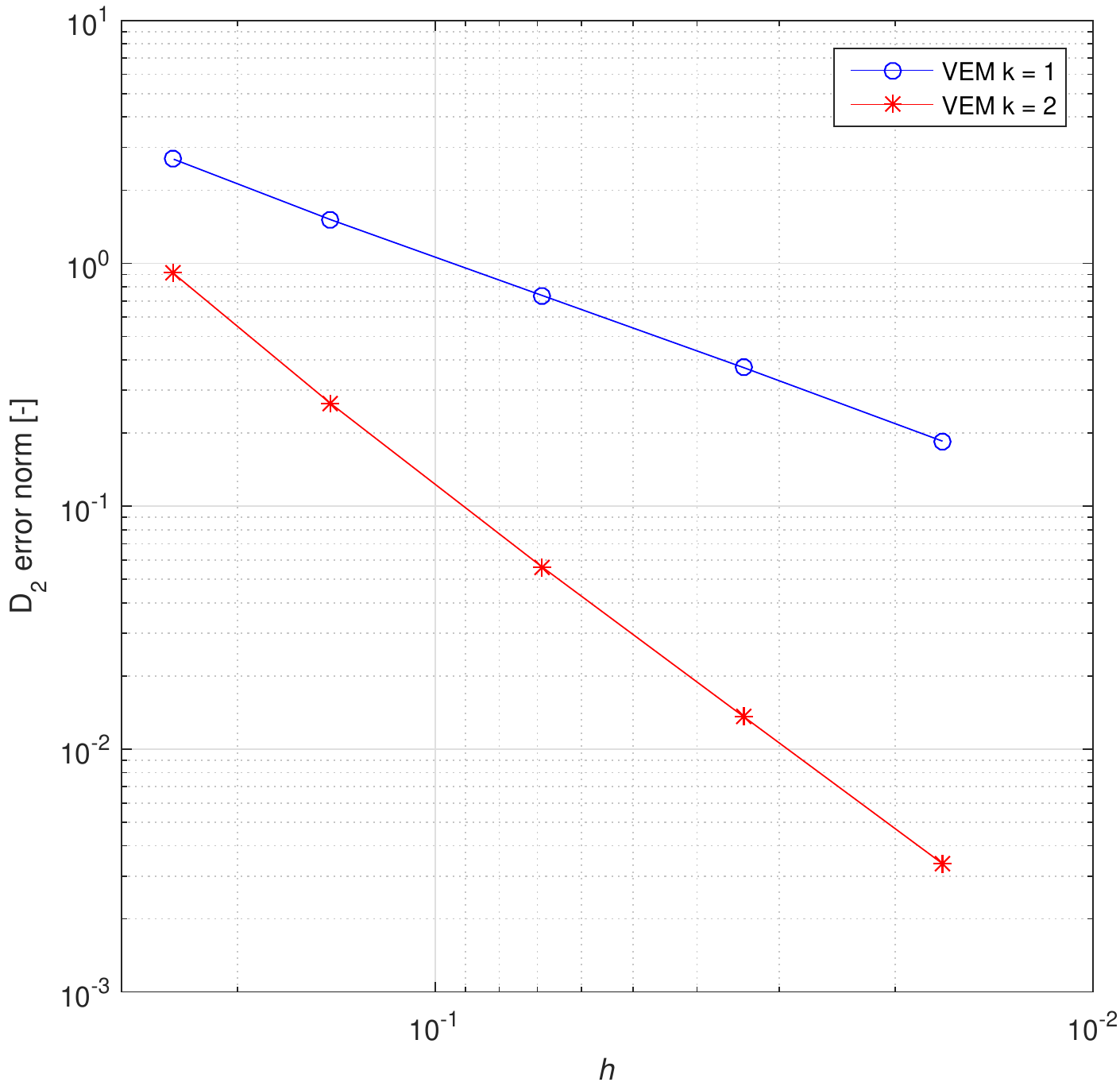}
\subcaption{}
\end{minipage}
\begin{minipage}[b]{.5\linewidth}
\centering
\includegraphics[bb = 80 5 400 400, scale=0.41]
                {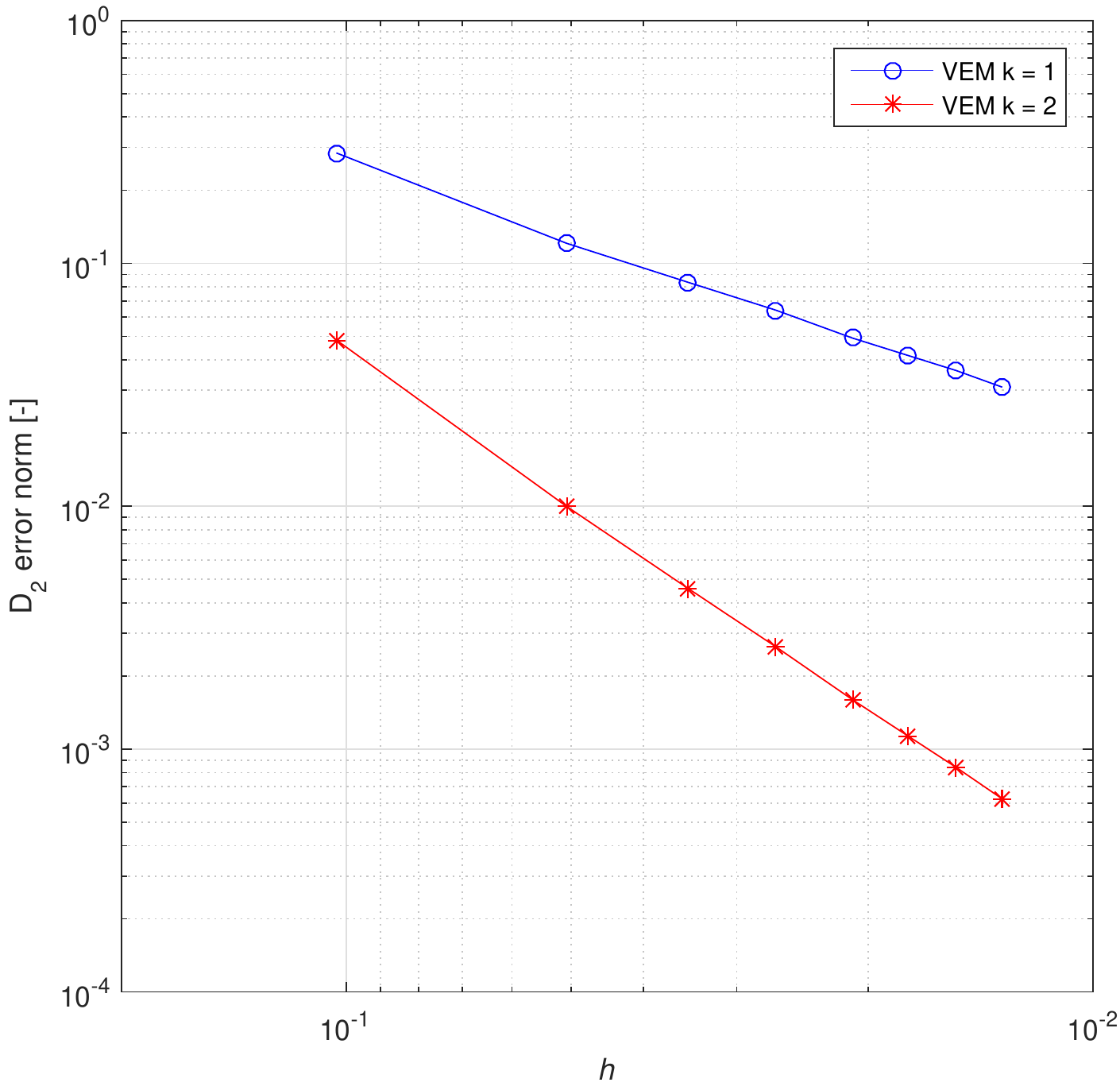}
\subcaption{}
\end{minipage}%
\begin{minipage}[b]{.5\linewidth}
\centering
\includegraphics[bb = 80 5 400 400, scale=0.41]
                {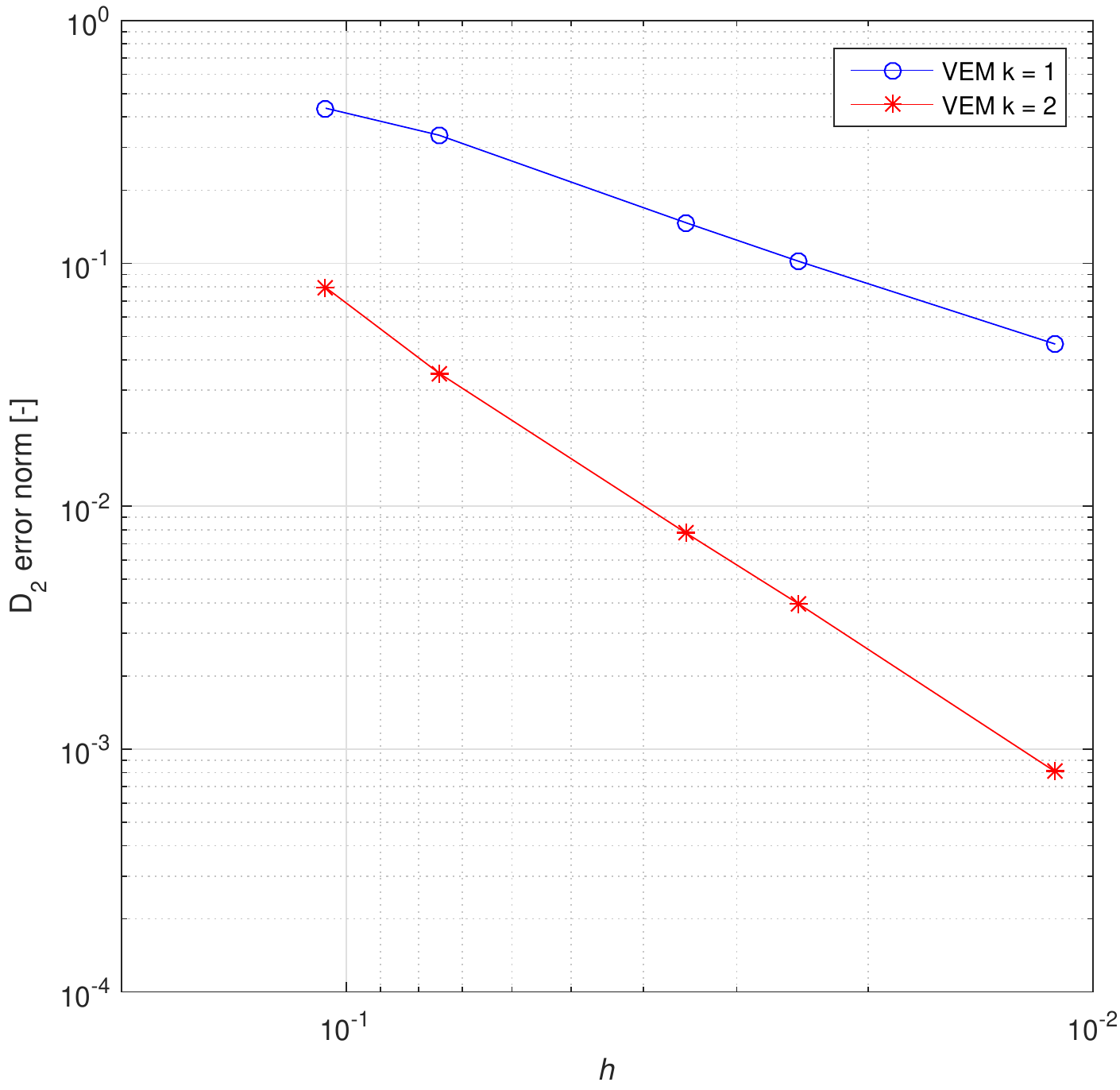}
\subcaption{}
\end{minipage}
\caption{
Test 2b: convergence plot. Error $D_2$ {\it vs} element size $h$. VEM $k=1$ and $k=2$. Non-uniform meshes of convex quadrilaterals (a); non-uniform meshes of triangles (b); uniform meshes of hexagons (c); Voronoi tessellations (d).}
\label{fig:Test_2_b_error_D2}
\end{figure}

Sample meshes for Test $2b$ are represented in Fig. \ref{fig:Test_2_b_mesh}, comprising non-uniform mesh of quadrilaterals (a), non-uniform mesh of triangles (b), uniform mesh of convex hexagons (c), and a centroid-based Voronoi tessellation (d). Correspondingly, Figures \ref{fig:Test_2_b_error_D1}(a)-(d) show the $h-$convergence plots in terms of error $D_1$ for the $k=1$, and $k=2$ VEM solutions. It is again observed that VEM show the expected convergence rates, which, as for Test $2a$, are maintained regardless of mesh non-uniformity and distortion. Similar considerations may be drawn in relation to the $h-$convergence pattern for error $D_2$, as can be appreciated in Fig. \ref{fig:Test_2_b_error_D2}.

\subsection{Test $3$}
\label{ss:test_3}
\begin{figure}
\vspace{20mm}
\begin{minipage}[b]{.5\linewidth}
\hspace{10mm}
\includegraphics[bb = 80 5 400 400, scale=0.42]
                {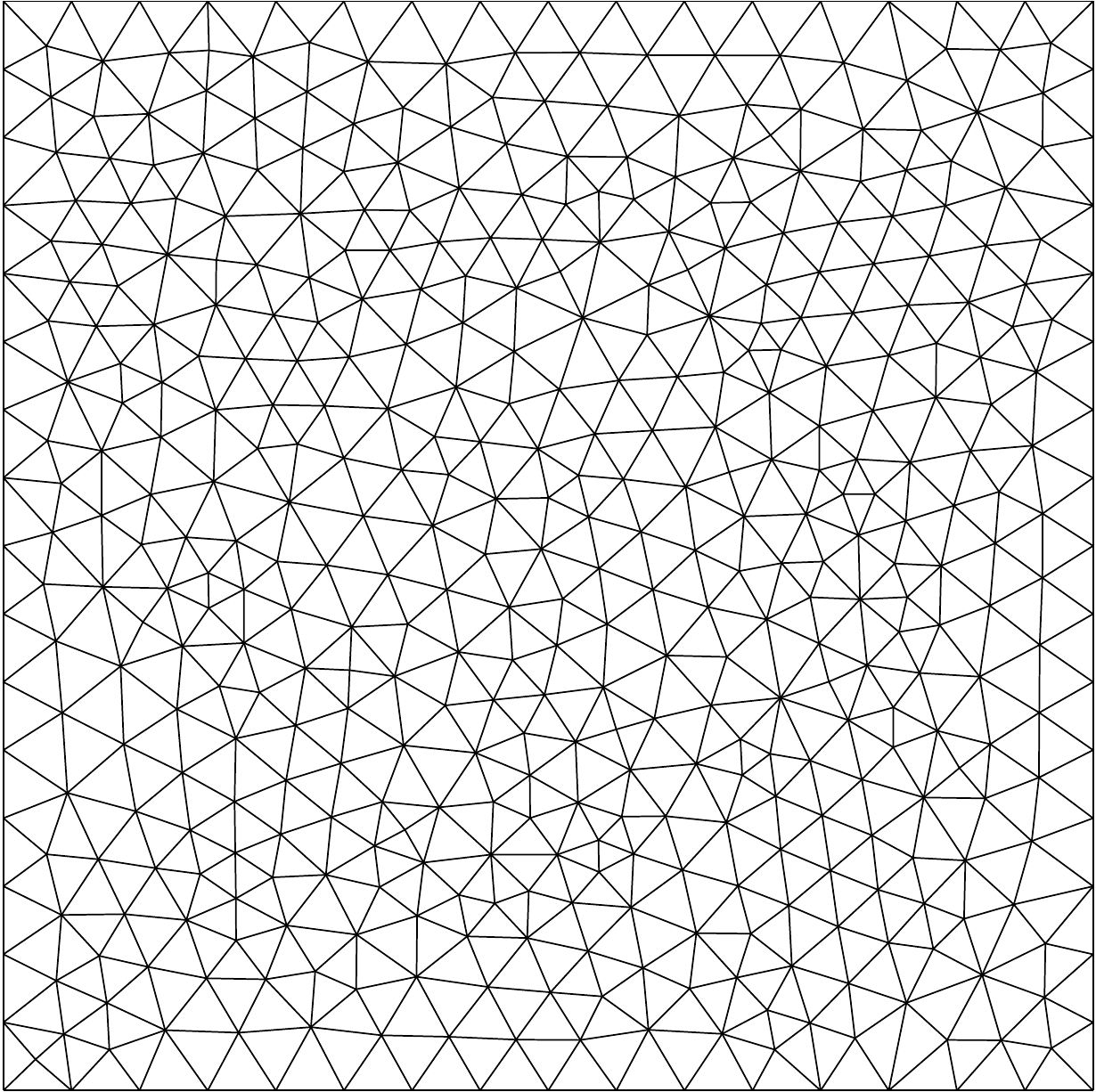}
\subcaption{}
\end{minipage}%
\begin{minipage}[b]{.5\linewidth}
\hspace{10mm}
\includegraphics[bb = 80 5 400 400, scale=0.42]
                {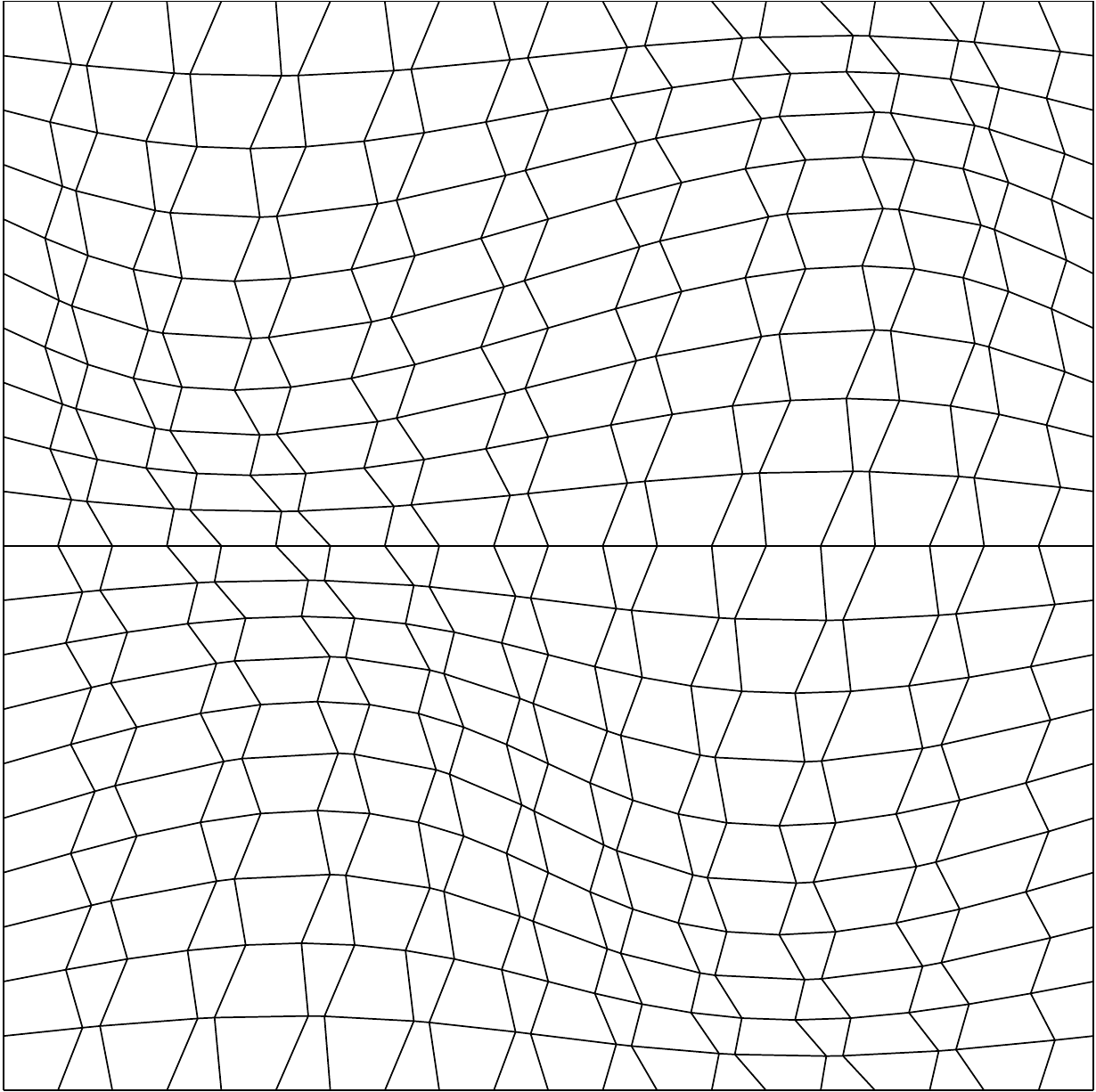}
\subcaption{}
\end{minipage}
\begin{minipage}[b]{.5\linewidth}
\hspace{10mm}
\includegraphics[bb = 80 5 400 400, scale=0.42]
                {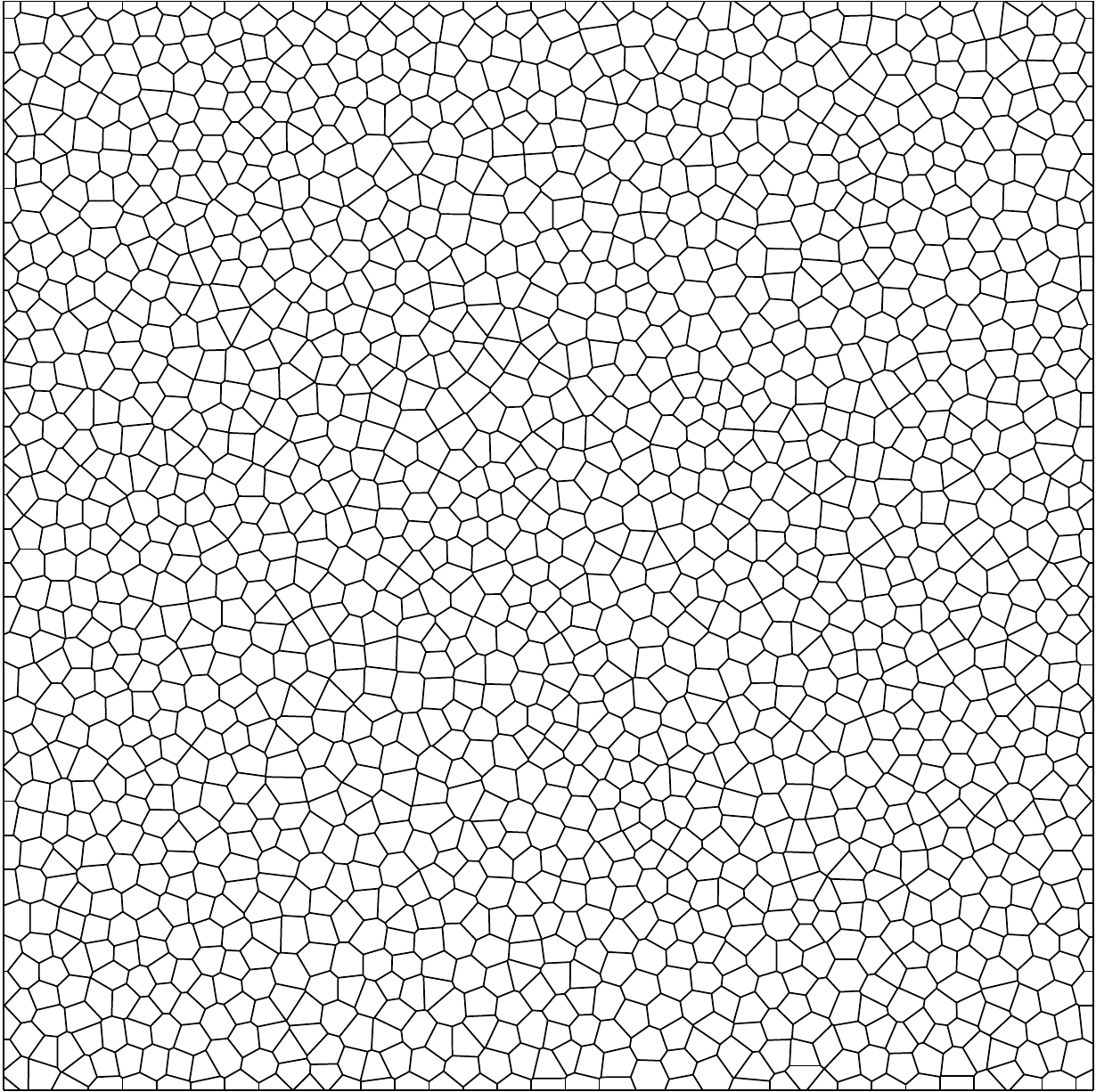}
\subcaption{}
\end{minipage}%
\begin{minipage}[b]{.5\linewidth}
\hspace{10mm}
\includegraphics[bb = 80 5 400 400, scale=0.42]
                {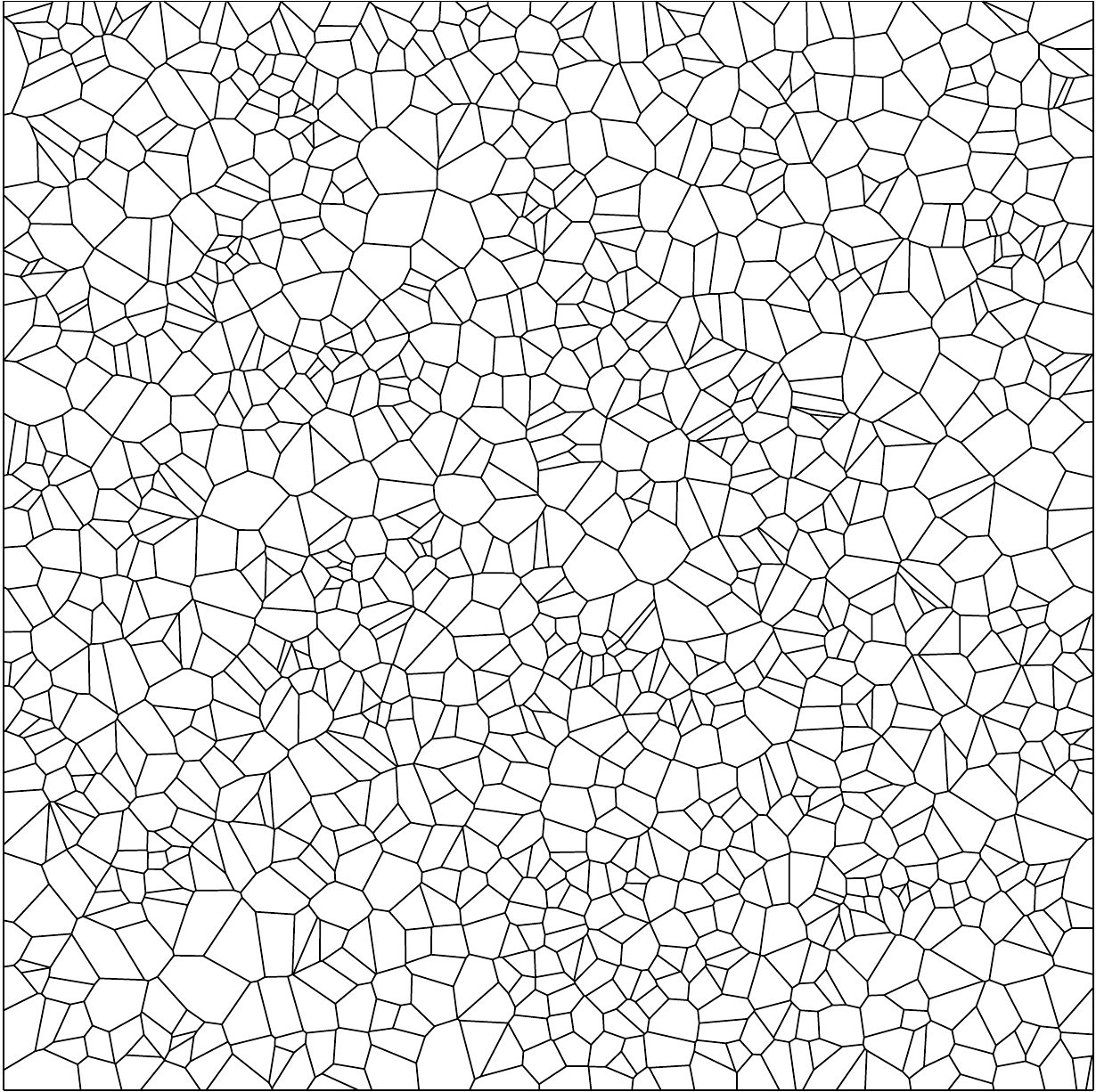}
\subcaption{}
\end{minipage}
\caption{
Test 3: Sample meshes. Non-uniform mesh of triangles (a); Distorted convex quadrilaterals (b); Centroid-based Voronoi tessellation (c); Random-based Voronoi tessellation (d).}
\label{fig:Test_3_mesh}
\end{figure}

\begin{figure}
\vspace{20mm}
\begin{minipage}[b]{.5\linewidth}
\centering
\includegraphics[bb = 80 5 400 400, scale=0.41]
                {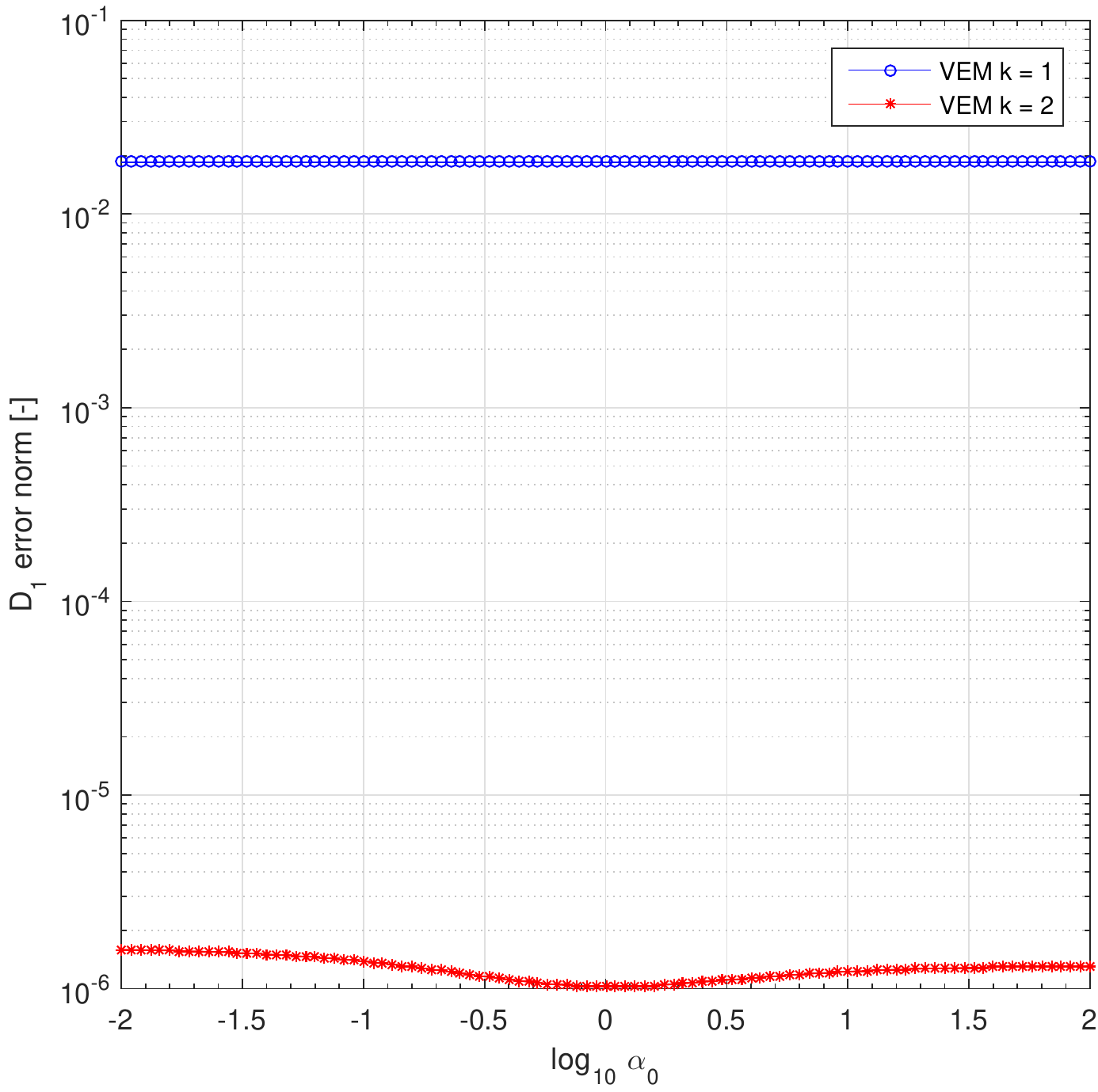}
\subcaption{}
\end{minipage}%
\begin{minipage}[b]{.5\linewidth}
\centering
\includegraphics[bb = 80 5 400 400, scale=0.41]
                {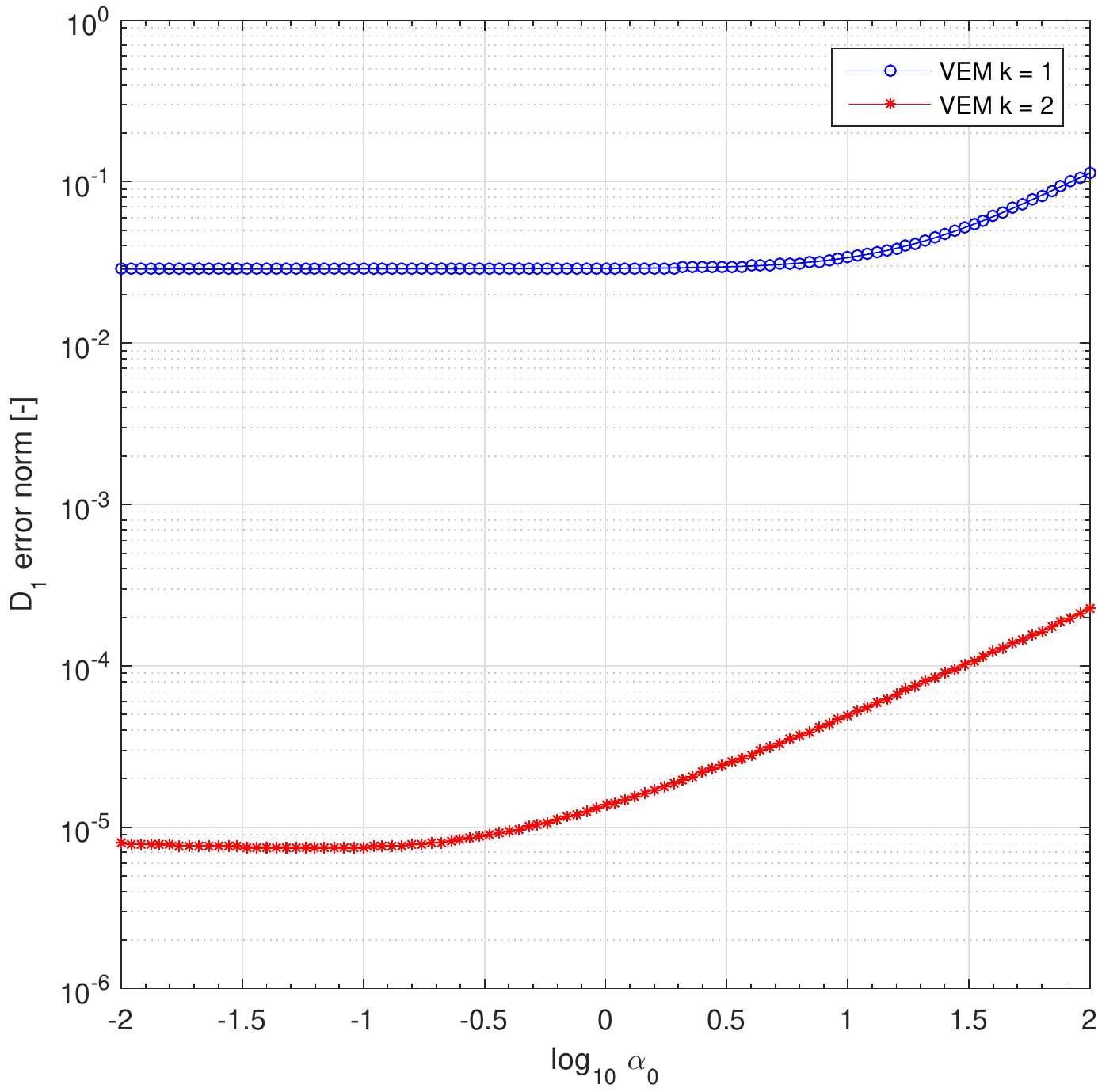}
\subcaption{}
\end{minipage}
\begin{minipage}[b]{.5\linewidth}
\centering
\includegraphics[bb = 80 5 400 400, scale=0.41]
                {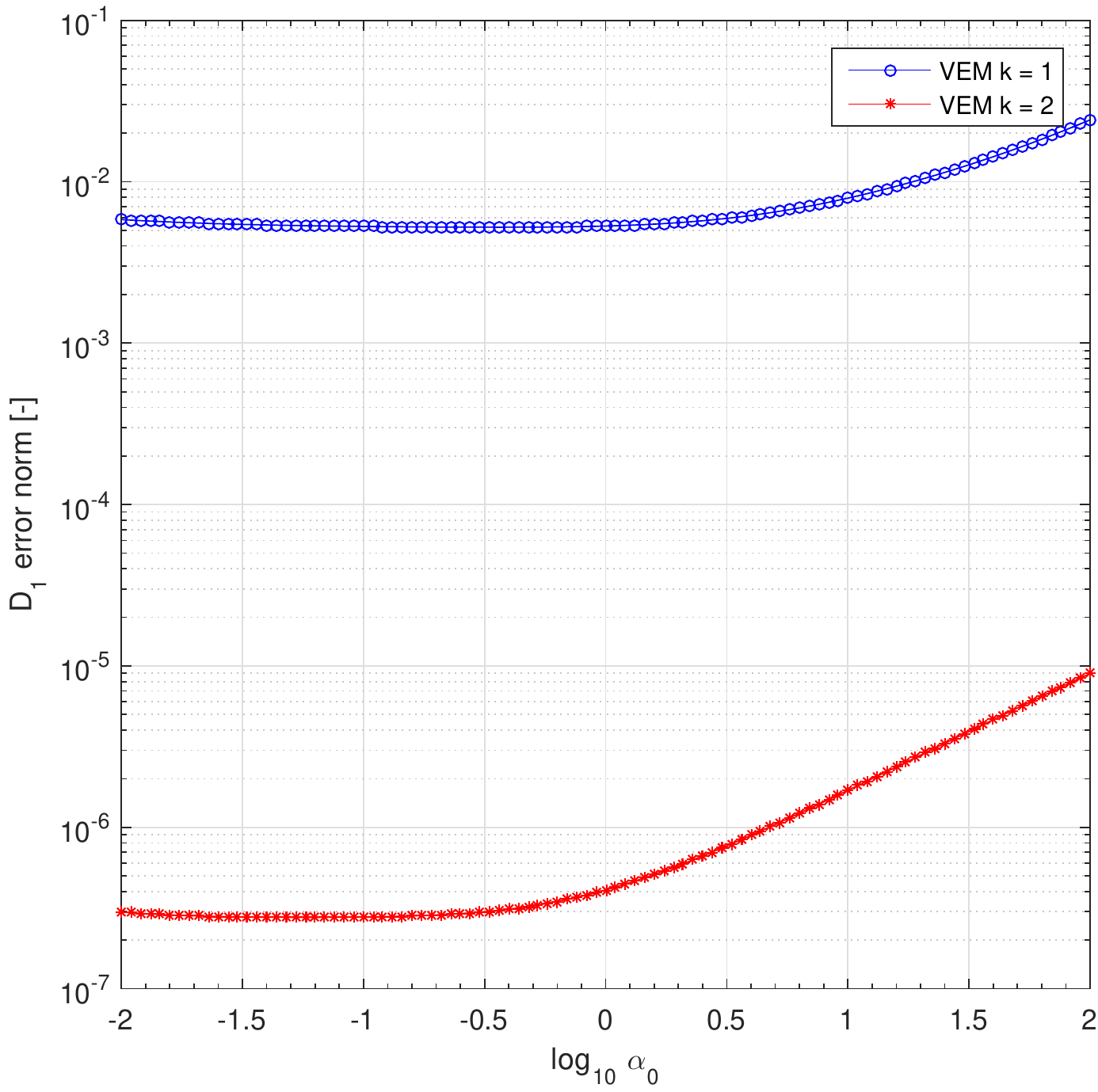}
\subcaption{}
\end{minipage}%
\begin{minipage}[b]{.5\linewidth}
\centering
\includegraphics[bb = 80 5 400 400, scale=0.41]
                {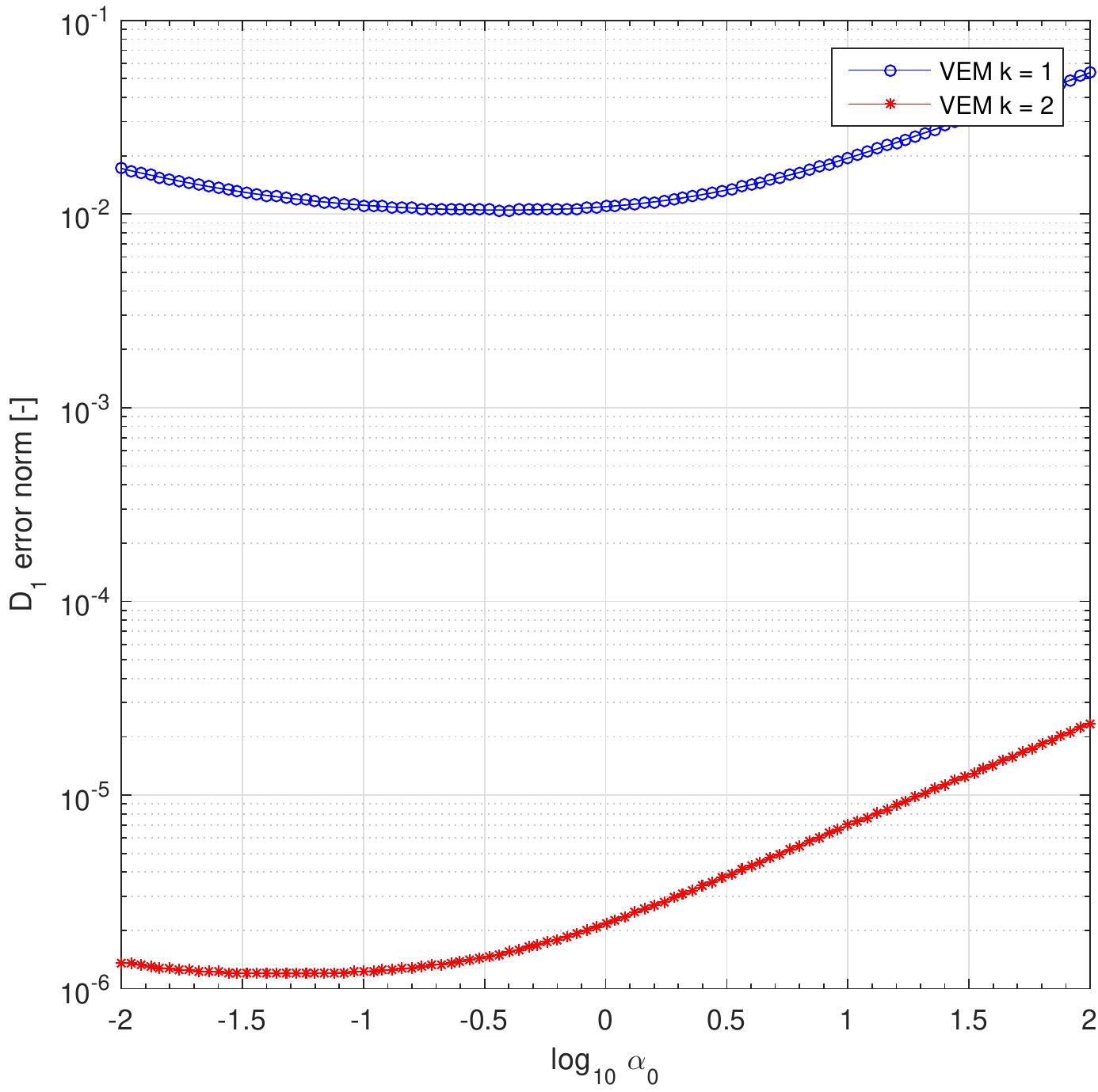}
\subcaption{}
\end{minipage}
\caption{
Test 3: sensitivity curves. Error $D_1$ {\it vs} stabilizing parameter factor $\alpha_0$, with $\tau=\alpha_0/2$. Non-uniform meshes of triangles (a); Distorted convex quadrilaterals (b); Centroid-based Voronoi tessellation (c); Random-based Voronoi tessellations (d).}
\label{fig:Test_3_sns_crv}
\end{figure}

The present section is devoted to the analysis of the influence on solution accuracy of the stabilization parameter, appearing in \eqref{eq:stab_term} pre-multiplying the stabilization stiffness, taken by default as $\tau=1/2$. Reference is made to Test $2a$, solved on the four different meshes reported in Fig. \ref{fig:Test_3_mesh}, comprising: non-uniform mesh of triangles (a), distorted convex quadrilaterals (b), a centroid-based Voronoi tessellation, and a random-based Voronoi tessellation. Correspondingly, we plot error $D_1$ (cf. \eqref{eq:energy_err_norm}) in Fig. \ref{fig:Test_3_sns_crv} against a parameter $\alpha_0$ ranging from $10^{-2}$ to $10^{2}$, with the corresponding stabilization parameter taken as $\tau=\alpha_0/2$, cf. equation \eqref{eq:stab_term}.  Note that for $\alpha_0=1$ the standard choice $\tau=1/2$ adopted in the tests of this paper is recovered, while for different values of $\alpha_0$ stronger or weaker stabilizations are introduced. 
We consider both $k = 1$ and $k=2$ VEMs. 
It can be observed that each examined case seems to be quite insensitive to the variation of $\tau$: the curves are quite flat and the overall variation of error with respect to the scaling is limited to within less than one order of magnitude. The proposed method thus seems to be quite robust with respect to the stabilization parameter.

\subsection{Test $4$}
\label{ss:test_4}
\begin{figure}
\begin{center}
\includegraphics[bb = 160 18 600 580, scale=0.25]
                {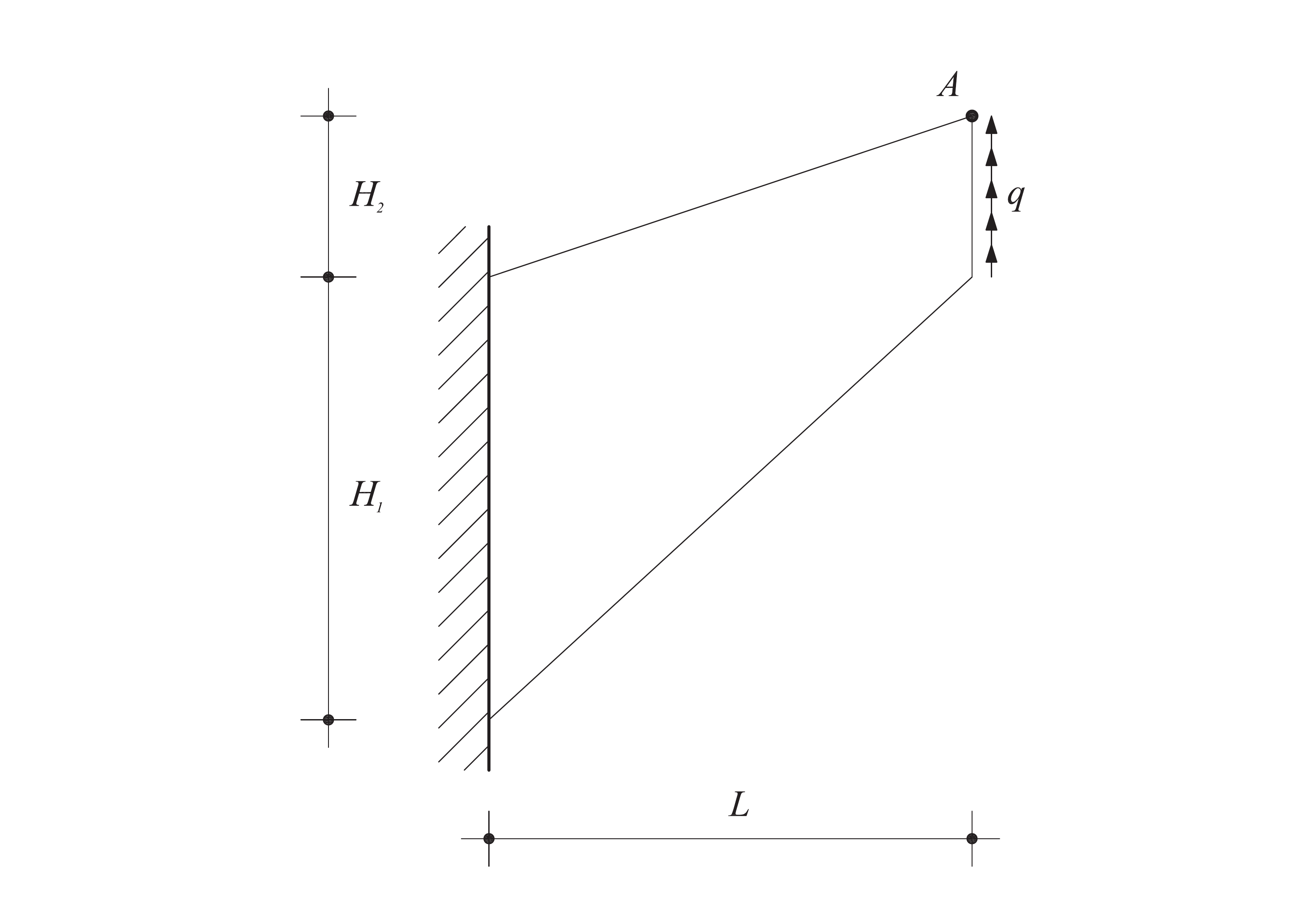}
\end{center}
\caption{
Test 4: Cook's membrane. Geometry, loading and boundary conditions.}
\label{fig:Test_4_geom}
\end{figure}

\begin{figure}

\begin{minipage}[b]{.5\linewidth}
\centering
\includegraphics[bb = 80 5 400 400, scale=0.42]
                {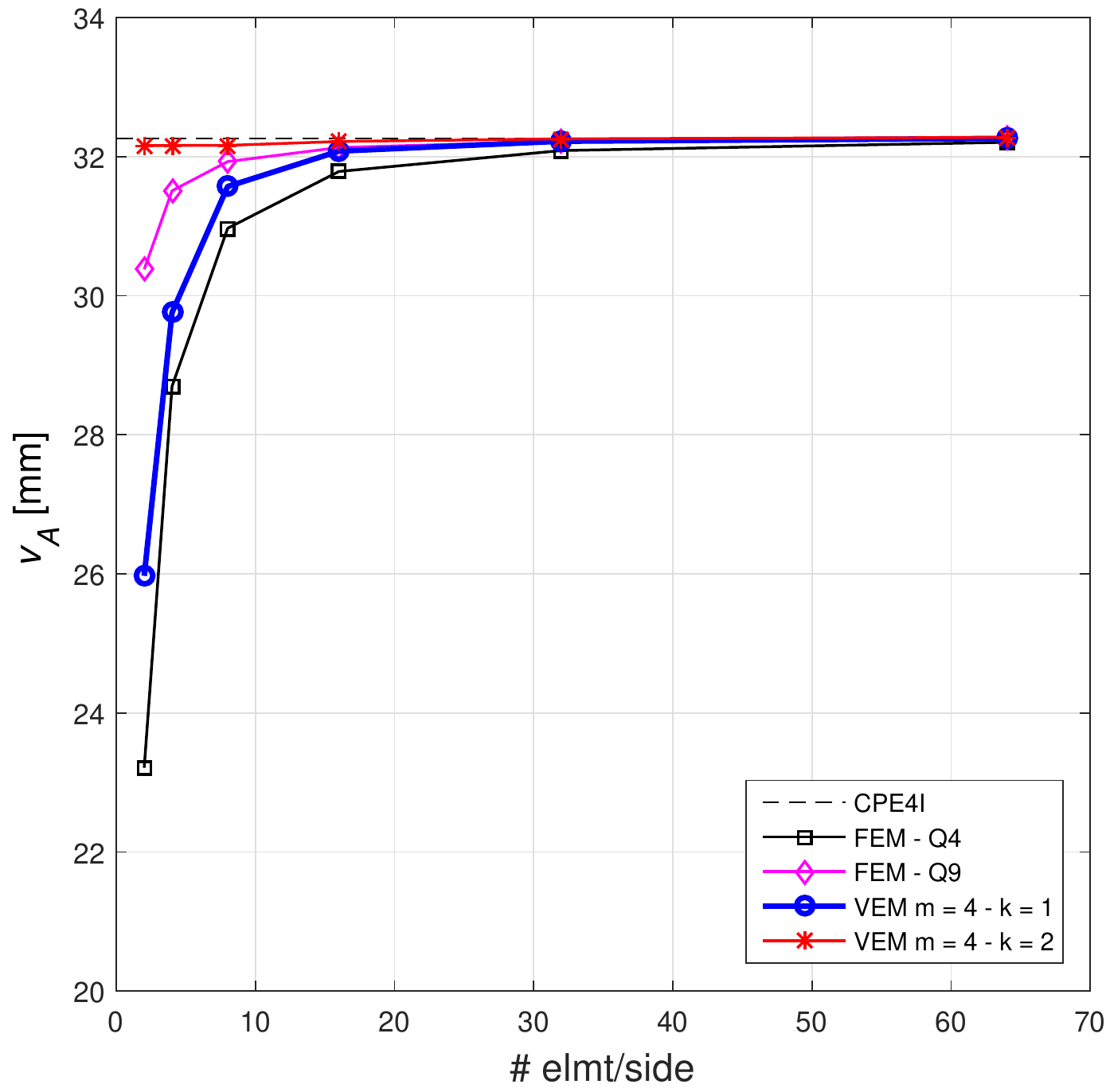}
\subcaption{}
\end{minipage}
\begin{minipage}[b]{.5\linewidth}
\centering
\includegraphics[bb = 80 5 400 400, scale=0.42]
                {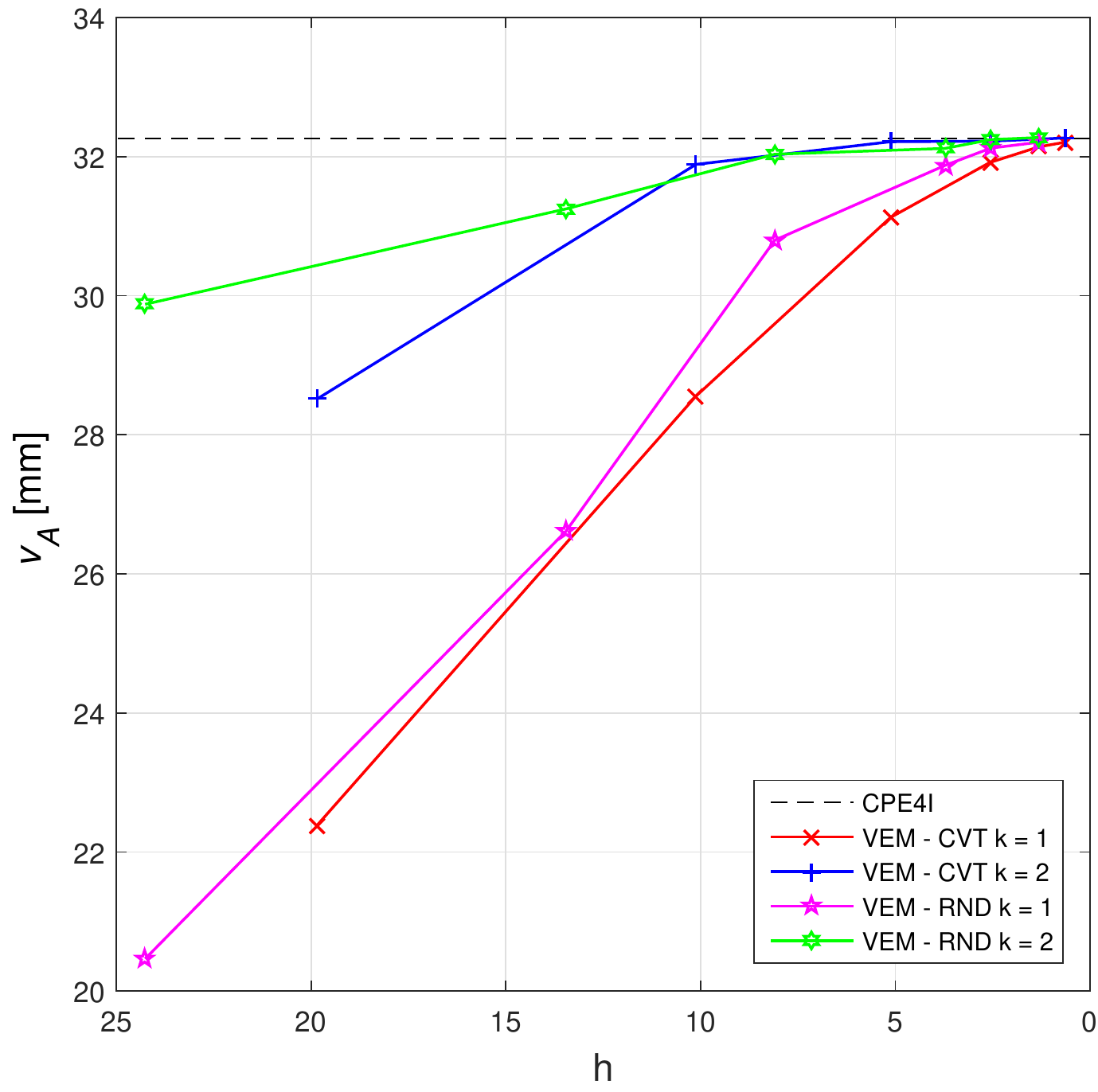}
\subcaption{}
\end{minipage}
\caption{
Test 4: Cook's membrane. Convergence results for vertical displacement of point $A$. Quad mesh with FEM $Q4$ and $Q9$; VEM $k=1$ and $k=2$, $m=4$ (a); Centroid-based and random-based Voronoi tessellations VEM $k=1$ and $k=2$ (b). Reference solution based on overkilling size mesh of quadrilateral CPE4I hybrid-mixed finite elements.}
\label{fig:Test_4_conv}
\end{figure}

The present section deals with the classical Cook's membrane 2D problem \cite{Zienckiewicz_Taylor2000}. Geometry of the domain $\O$ is presented in Fig. \ref{fig:Test_4_geom} with length data $H_1 = 44$, $H_2 = 16$, $L = 48$. Material parameters are Young modulus $E = 70$, and Poisson ration $\nu = 1/3$. The loading is given by a constant tangential traction $q = 6.25$ on the right edge of the domain. The problem is solved using two types of meshes: an evenly distributed quad mesh (see fig. \ref{fig:Test_4_geom}(a)), and centroid based [resp. random-based] Voronoi tessellations (see fig. \ref{fig:Test_4_geom}(b)). The former type is used for linear and quadratic VEMs in conjunction with $Q4$ and $Q9$ FEM for comparison purposes; whereas the latter types are adopted in conjunction to VEMs only. Convergence results are reported in terms of mesh refinement monitoring $v_A$, the vertical displacement of point A (see Fig. \ref{fig:Test_4_geom}), as can be appreciated in Fig. \ref{fig:Test_4_conv}(a)-(b), where a reference solution is indicated with a dotted black line corresponding to an overkilling accurate solution obtained with the hybrid-mixed CPE4I element \cite{ABAQUS:2011}. It is observed that quadrilateral VEMs have a slight edge in terms of accuracy with respect to Lagrangian finite elements of the corresponding order, while the polygonal VEMs from the centroid-based Voronoi tessellation slightly outperform the random-based VEMs in terms of accuracy.

\begin{figure}

\begin{minipage}[b]{.5\linewidth}
\centering
\includegraphics[bb = 160 18 550 480, scale=0.5]
                {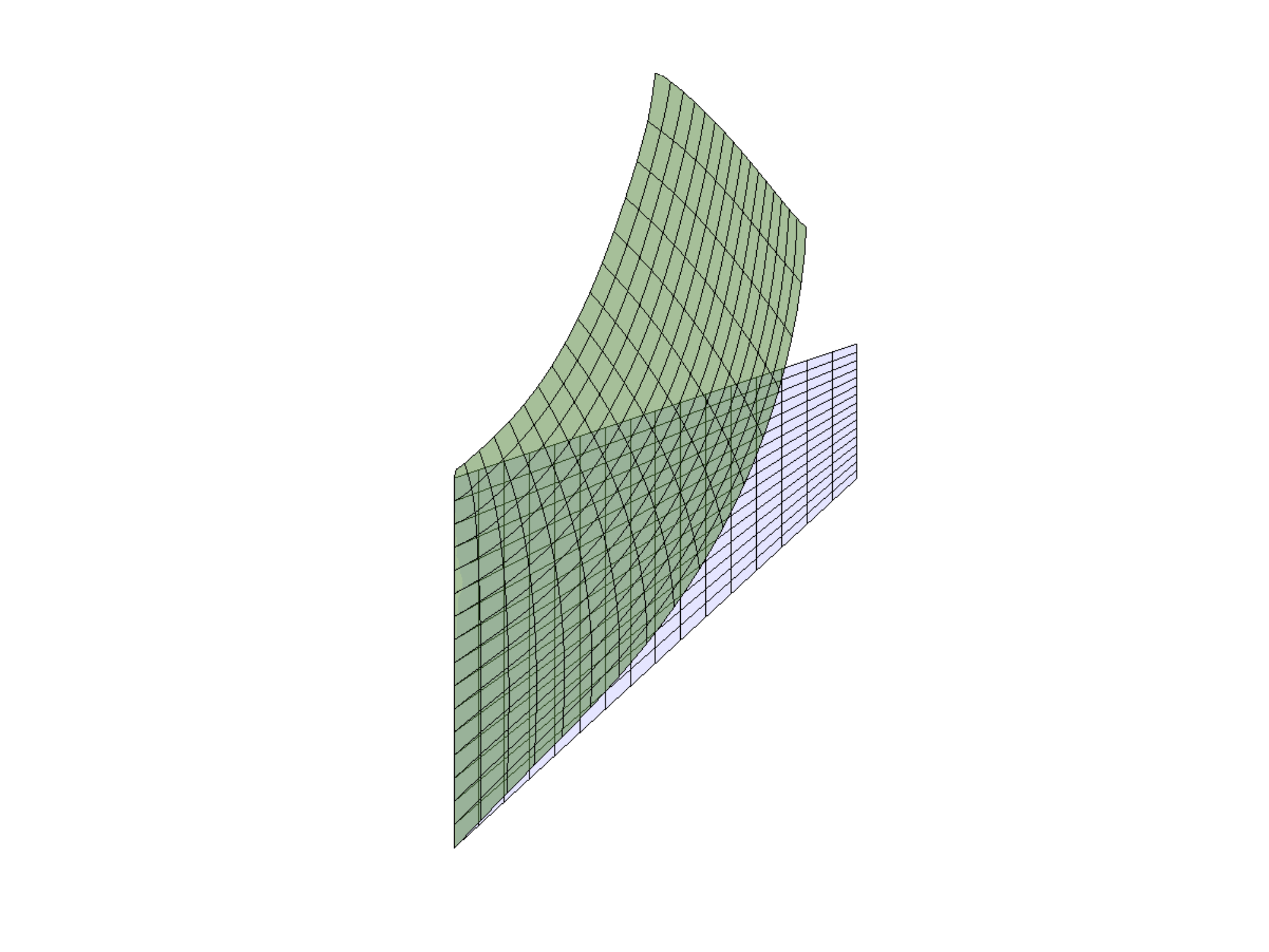}
\subcaption{}
\end{minipage}
\begin{minipage}[b]{.5\linewidth}
\centering
\includegraphics[bb = 250 18 300 480, scale=0.5]
                {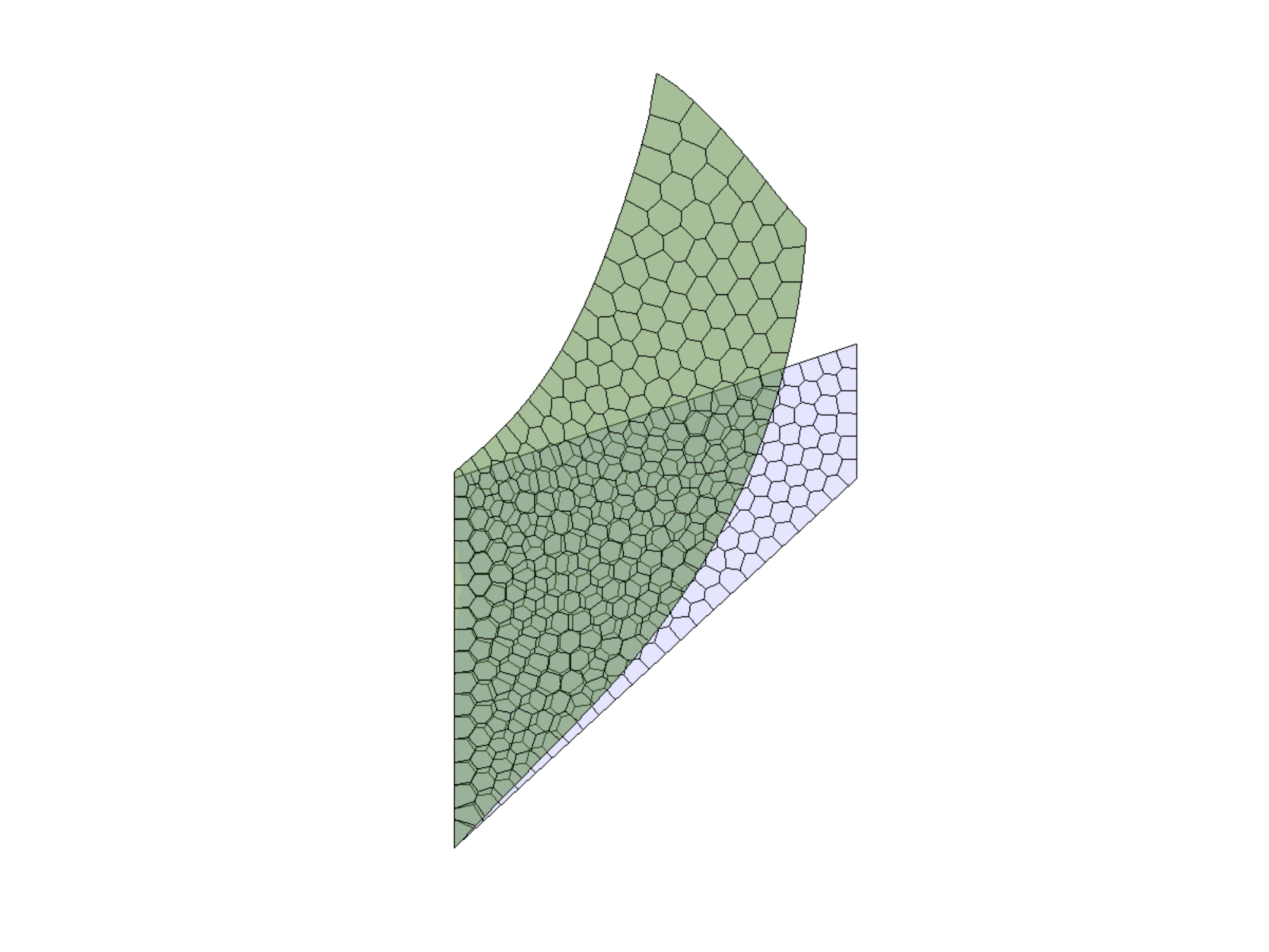}
\subcaption{}
\end{minipage}
\caption{
Test 4: Cook's membrane: undeformed and deformed configurations. Quad mesh, VEM $k=2$ (a); Centroid-based Voronoi tessellations, VEM $k=2$ (b).}
\label{fig:Test_4_def_con}
\end{figure}

\begin{figure}
\begin{minipage}[b]{.5\linewidth}
\hspace{13mm}
\includegraphics[bb = 80 5 400 400, scale=0.42]
                {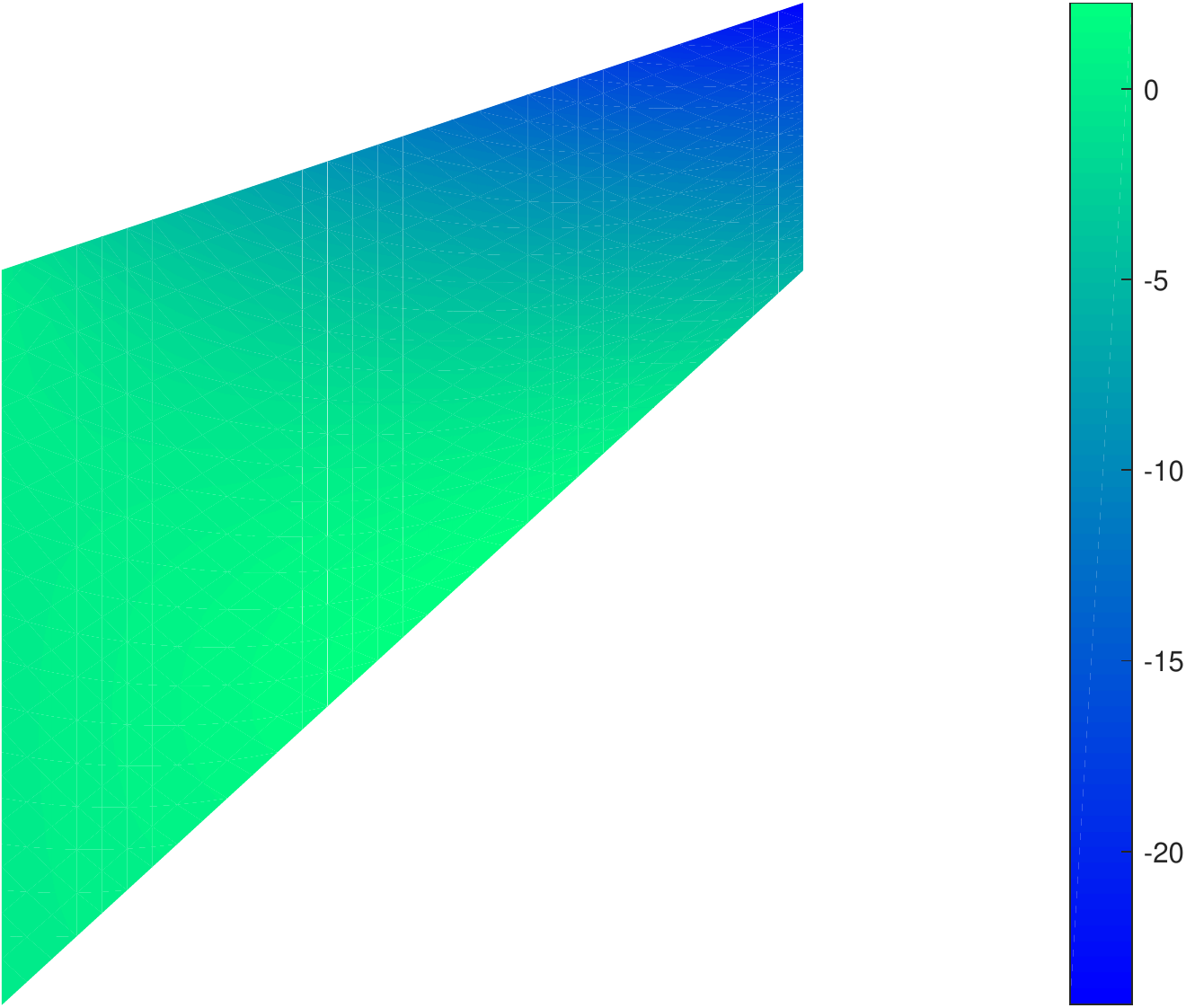}
\subcaption{}
\end{minipage}
\begin{minipage}[b]{.5\linewidth}
\hspace{13mm}
\includegraphics[bb = 80 5 400 400, scale=0.42]
                {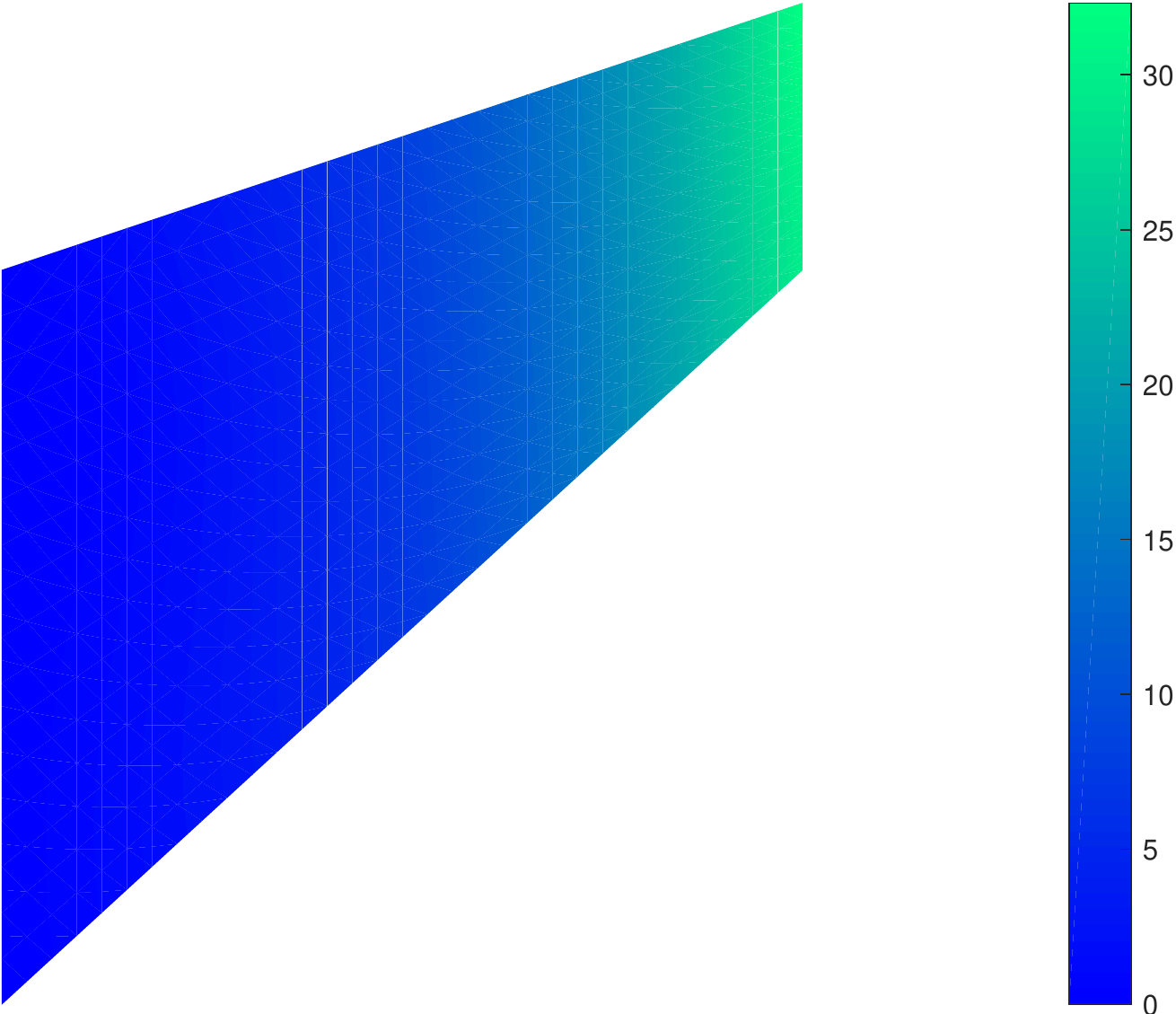}
\subcaption{}
\end{minipage}
\caption{
Test 4: Cook's membrane. Color plots of displacement field components $u$-(a), $v$-(b) on undeformed configuration. Centroid-based Voronoi tessellation for VEM $k=2$ solution.}
\label{fig:Test_4_disp}
\end{figure}

In Fig. \ref{fig:Test_4_def_con} the deformed configuration corresponding to quadratic quadrilaterals (a) and centroidal Voronoi polygons (b) is reported, showing a perfect agreement between the two solutions and hence no sensitivity to mesh type, distortion and non-uniformity. Finally, in Fig. \ref{fig:Test_4_disp}, the qualitative color plots of the displacement components stemming from the latter solution are displayed showing the expected smoothness of the displacement field over the mesh.

\section{Conclusion}
\label{s:conclusion}
A Virtual Element method for plane elasticity problems has been presented. The newly developed methodology, together with implementation details, has been illustrated for the case of quadrilateral elements with linear approximation of the displacement field, and for general polygonal elements with quadratic and higher order interpolation. The stiffness matrix and loading vector have been explicitly derived proving the straightforward implementability of the method into into a computer code for structural analysis.
An extensive campaign of numerical applications illustrate accuracy and convergence patterns of the method, and its performances in comparison with a classical finite element approach characterized by the same order of approximation of the unknown fields.
An accurate investigation on the convergence for structured and non structured discretizations has been performed, showing the superior robustness of the VEM method with respect to a FEM approach. It can be stated that VEM is almost insensitive to the mesh distortion.
The assessed reliability and the good performances of the VEM formulation indicate that the virtual element approach can rapidly become a numerical strategy largely diffused and used for structural computations.

\clearpage
\newpage

\clearpage
\newpage

\clearpage
\newpage

\clearpage
\newpage

\clearpage
\newpage
\clearpage
\newpage
\clearpage
\newpage
\newpage
\newpage
\newpage
\newpage
\newpage
\newpage
\newpage
\newpage
\begin{figure}

\begin{minipage}[b]{.5\linewidth}
\hspace{15mm}
\includegraphics[bb = 80 5 400 400, scale=0.42]
                {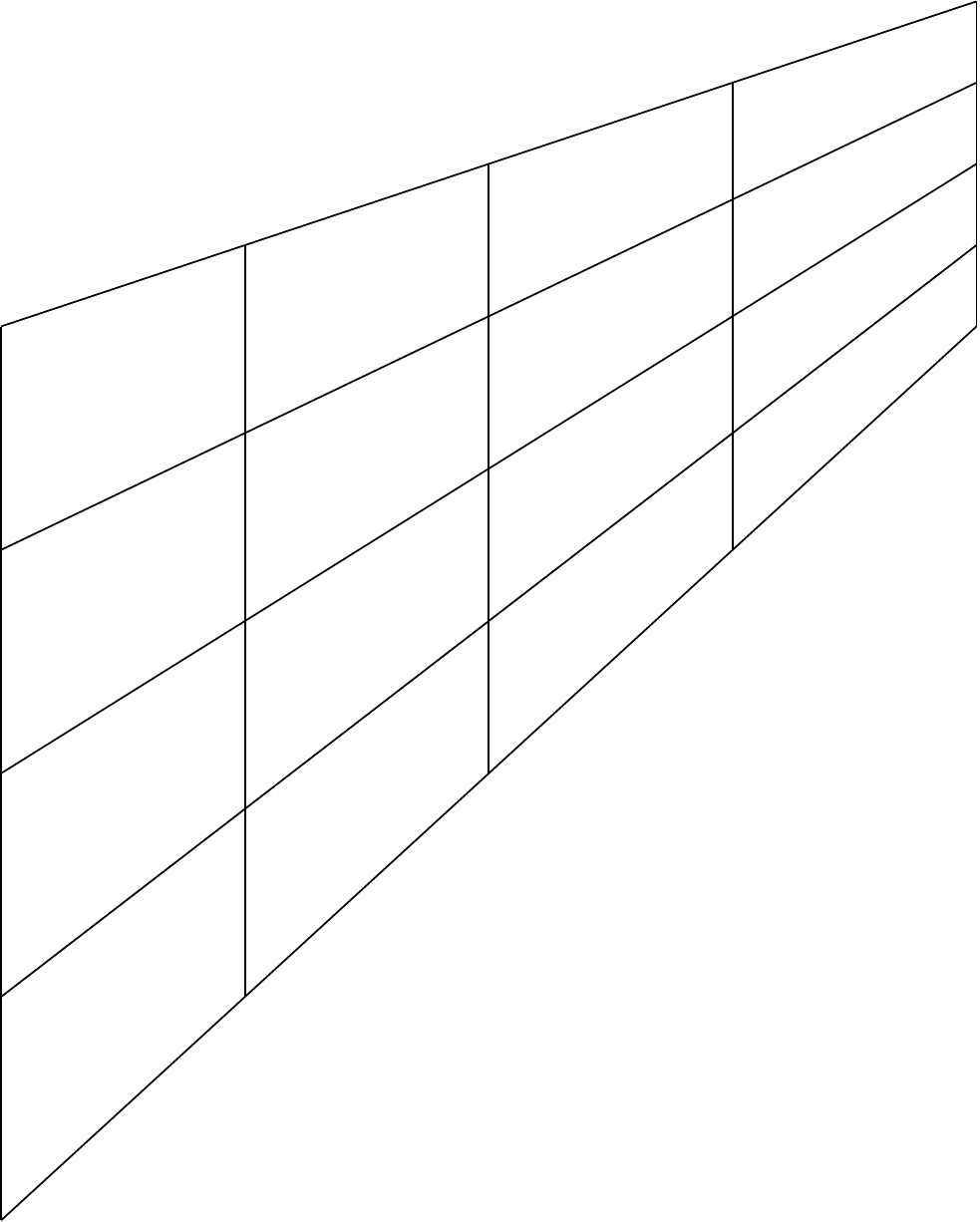}
\subcaption{}
\end{minipage}
\hspace{5mm}%
\begin{minipage}[b]{.5\linewidth}
\hspace{15mm}
\includegraphics[bb = 80 5 400 400, scale=0.42]
                {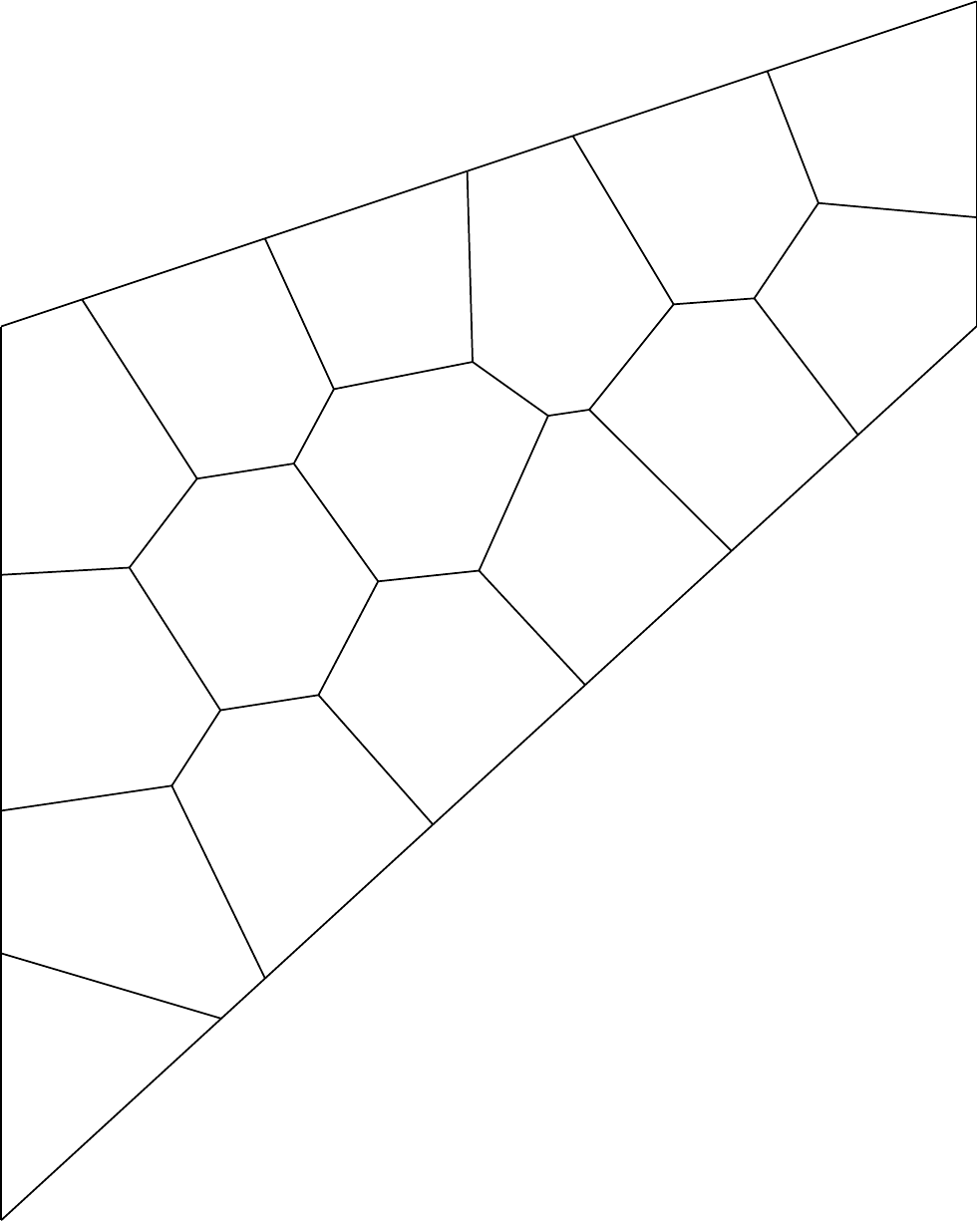}
\subcaption{}
\end{minipage}
\caption{
Test 4: Cook's membrane. Sample meshes: quad mesh (a); centroid-based Voronoi tessellation (b).}
\label{fig:Test_4_mesh}
\end{figure}
\newpage
\newpage
\newpage
\clearpage
\newpage
\bibliographystyle{amsplain}

\bibliography{general-bibliography}

\end{document}